\documentclass[aos,preprint]{imsart}

\RequirePackage[OT1]{fontenc}
\RequirePackage{amsthm,amsmath,natbib}
\RequirePackage[colorlinks,citecolor=blue,urlcolor=blue]{hyperref}

\usepackage{graphicx}

\setattribute{journal}{name}{}

\startlocaldefs
\numberwithin{equation}{section}
\theoremstyle{plain}

\endlocaldefs 
 
\defcitealias{Cut2015}{CPV15}

\usepackage{color}   
\definecolor{AliceBlue}{rgb}{0.94,0.97,1.00}
\definecolor{AntiqueWhite1}{rgb}{1.00,0.94,0.86}
\definecolor{AntiqueWhite2}{rgb}{0.93,0.87,0.80}
\definecolor{AntiqueWhite3}{rgb}{0.80,0.75,0.69}
\definecolor{AntiqueWhite4}{rgb}{0.55,0.51,0.47}
\definecolor{AntiqueWhite}{rgb}{0.98,0.92,0.84}
\definecolor{BlanchedAlmond}{rgb}{1.00,0.92,0.80}
\definecolor{BlueViolet}{rgb}{0.54,0.17,0.89}
\definecolor{CadetBlue1}{rgb}{0.60,0.96,1.00}
\definecolor{CadetBlue2}{rgb}{0.56,0.90,0.93}
\definecolor{CadetBlue3}{rgb}{0.48,0.77,0.80}
\definecolor{CadetBlue4}{rgb}{0.33,0.53,0.55}
\definecolor{CadetBlue}{rgb}{0.37,0.62,0.63}
\definecolor{CornflowerBlue}{rgb}{0.39,0.58,0.93}
\definecolor{DarkBlue}{rgb}{0.00,0.00,0.55}
\definecolor{DarkCyan}{rgb}{0.00,0.55,0.55}
\definecolor{DarkGoldenrod1}{rgb}{1.00,0.73,0.06}
\definecolor{DarkGoldenrod2}{rgb}{0.93,0.68,0.05}
\definecolor{DarkGoldenrod3}{rgb}{0.80,0.58,0.05}
\definecolor{DarkGoldenrod4}{rgb}{0.55,0.40,0.03}
\definecolor{DarkGoldenrod}{rgb}{0.72,0.53,0.04}
\definecolor{DarkGray}{rgb}{0.66,0.66,0.66}
\definecolor{DarkGreen}{rgb}{0.00,0.39,0.00}
\definecolor{DarkGrey}{rgb}{0.66,0.66,0.66}
\definecolor{DarkKhaki}{rgb}{0.74,0.72,0.42}
\definecolor{DarkMagenta}{rgb}{0.55,0.00,0.55}
\definecolor{DarkOliveGreen1}{rgb}{0.79,1.00,0.44}
\definecolor{DarkOliveGreen2}{rgb}{0.74,0.93,0.41}
\definecolor{DarkOliveGreen3}{rgb}{0.64,0.80,0.35}
\definecolor{DarkOliveGreen4}{rgb}{0.43,0.55,0.24}
\definecolor{DarkOliveGreen}{rgb}{0.33,0.42,0.18}
\definecolor{DarkOrange1}{rgb}{1.00,0.50,0.00}
\definecolor{DarkOrange2}{rgb}{0.93,0.46,0.00}
\definecolor{DarkOrange3}{rgb}{0.80,0.40,0.00}
\definecolor{DarkOrange4}{rgb}{0.55,0.27,0.00}
\definecolor{DarkOrange}{rgb}{1.00,0.55,0.00}
\definecolor{DarkOrchid1}{rgb}{0.75,0.24,1.00}
\definecolor{DarkOrchid2}{rgb}{0.70,0.23,0.93}
\definecolor{DarkOrchid3}{rgb}{0.60,0.20,0.80}
\definecolor{DarkOrchid4}{rgb}{0.41,0.13,0.55}
\definecolor{DarkOrchid}{rgb}{0.60,0.20,0.80}
\definecolor{DarkRed}{rgb}{0.55,0.00,0.00}
\definecolor{DarkSalmon}{rgb}{0.91,0.59,0.48}
\definecolor{DarkSeaGreen1}{rgb}{0.76,1.00,0.76}
\definecolor{DarkSeaGreen2}{rgb}{0.71,0.93,0.71}
\definecolor{DarkSeaGreen3}{rgb}{0.61,0.80,0.61}
\definecolor{DarkSeaGreen4}{rgb}{0.41,0.55,0.41}
\definecolor{DarkSeaGreen}{rgb}{0.56,0.74,0.56}
\definecolor{DarkSlateBlue}{rgb}{0.28,0.24,0.55}
\definecolor{DarkSlateGray1}{rgb}{0.59,1.00,1.00}
\definecolor{DarkSlateGray2}{rgb}{0.55,0.93,0.93}
\definecolor{DarkSlateGray3}{rgb}{0.47,0.80,0.80}
\definecolor{DarkSlateGray4}{rgb}{0.32,0.55,0.55}
\definecolor{DarkSlateGray}{rgb}{0.18,0.31,0.31}
\definecolor{DarkSlateGrey}{rgb}{0.18,0.31,0.31}
\definecolor{DarkTurquoise}{rgb}{0.00,0.81,0.82}
\definecolor{DarkViolet}{rgb}{0.58,0.00,0.83}
\definecolor{DeepPink1}{rgb}{1.00,0.08,0.58}
\definecolor{DeepPink2}{rgb}{0.93,0.07,0.54}
\definecolor{DeepPink3}{rgb}{0.80,0.06,0.46}
\definecolor{DeepPink4}{rgb}{0.55,0.04,0.31}
\definecolor{DeepPink}{rgb}{1.00,0.08,0.58}
\definecolor{DeepSkyBlue1}{rgb}{0.00,0.75,1.00}
\definecolor{DeepSkyBlue2}{rgb}{0.00,0.70,0.93}
\definecolor{DeepSkyBlue3}{rgb}{0.00,0.60,0.80}
\definecolor{DeepSkyBlue4}{rgb}{0.00,0.41,0.55}
\definecolor{DeepSkyBlue}{rgb}{0.00,0.75,1.00}
\definecolor{DimGray}{rgb}{0.41,0.41,0.41}
\definecolor{DimGrey}{rgb}{0.41,0.41,0.41}
\definecolor{DodgerBlue1}{rgb}{0.12,0.56,1.00}
\definecolor{DodgerBlue2}{rgb}{0.11,0.53,0.93}
\definecolor{DodgerBlue3}{rgb}{0.09,0.45,0.80}
\definecolor{DodgerBlue4}{rgb}{0.06,0.31,0.55}
\definecolor{DodgerBlue}{rgb}{0.12,0.56,1.00}
\definecolor{FloralWhite}{rgb}{1.00,0.98,0.94}
\definecolor{ForestGreen}{rgb}{0.13,0.55,0.13}
\definecolor{GhostWhite}{rgb}{0.97,0.97,1.00}
\definecolor{GreenYellow}{rgb}{0.68,1.00,0.18}
\definecolor{HotPink1}{rgb}{1.00,0.43,0.71}
\definecolor{HotPink2}{rgb}{0.93,0.42,0.65}
\definecolor{HotPink3}{rgb}{0.80,0.38,0.56}
\definecolor{HotPink4}{rgb}{0.55,0.23,0.38}
\definecolor{HotPink}{rgb}{1.00,0.41,0.71}
\definecolor{IndianRed1}{rgb}{1.00,0.42,0.42}
\definecolor{IndianRed2}{rgb}{0.93,0.39,0.39}
\definecolor{IndianRed3}{rgb}{0.80,0.33,0.33}
\definecolor{IndianRed4}{rgb}{0.55,0.23,0.23}
\definecolor{IndianRed}{rgb}{0.80,0.36,0.36}
\definecolor{LavenderBlush1}{rgb}{1.00,0.94,0.96}
\definecolor{LavenderBlush2}{rgb}{0.93,0.88,0.90}
\definecolor{LavenderBlush3}{rgb}{0.80,0.76,0.77}
\definecolor{LavenderBlush4}{rgb}{0.55,0.51,0.53}
\definecolor{LavenderBlush}{rgb}{1.00,0.94,0.96}
\definecolor{LawnGreen}{rgb}{0.49,0.99,0.00}
\definecolor{LemonChiffon1}{rgb}{1.00,0.98,0.80}
\definecolor{LemonChiffon2}{rgb}{0.93,0.91,0.75}
\definecolor{LemonChiffon3}{rgb}{0.80,0.79,0.65}
\definecolor{LemonChiffon4}{rgb}{0.55,0.54,0.44}
\definecolor{LemonChiffon}{rgb}{1.00,0.98,0.80}
\definecolor{LightBlue1}{rgb}{0.75,0.94,1.00}
\definecolor{LightBlue2}{rgb}{0.70,0.87,0.93}
\definecolor{LightBlue3}{rgb}{0.60,0.75,0.80}
\definecolor{LightBlue4}{rgb}{0.41,0.51,0.55}
\definecolor{LightBlue}{rgb}{0.68,0.85,0.90}
\definecolor{LightCoral}{rgb}{0.94,0.50,0.50}
\definecolor{LightCyan1}{rgb}{0.88,1.00,1.00}
\definecolor{LightCyan2}{rgb}{0.82,0.93,0.93}
\definecolor{LightCyan3}{rgb}{0.71,0.80,0.80}
\definecolor{LightCyan4}{rgb}{0.48,0.55,0.55}
\definecolor{LightCyan}{rgb}{0.88,1.00,1.00}
\definecolor{LightGoldenrod1}{rgb}{1.00,0.93,0.55}
\definecolor{LightGoldenrod2}{rgb}{0.93,0.86,0.51}
\definecolor{LightGoldenrod3}{rgb}{0.80,0.75,0.44}
\definecolor{LightGoldenrod4}{rgb}{0.55,0.51,0.30}
\definecolor{LightGoldenrodYellow}{rgb}{0.98,0.98,0.82}
\definecolor{LightGoldenrod}{rgb}{0.93,0.87,0.51}
\definecolor{LightGray}{rgb}{0.83,0.83,0.83}
\definecolor{LightGreen}{rgb}{0.56,0.93,0.56}
\definecolor{LightGrey}{rgb}{0.83,0.83,0.83}
\definecolor{LightPink1}{rgb}{1.00,0.68,0.73}
\definecolor{LightPink2}{rgb}{0.93,0.64,0.68}
\definecolor{LightPink3}{rgb}{0.80,0.55,0.58}
\definecolor{LightPink4}{rgb}{0.55,0.37,0.40}
\definecolor{LightPink}{rgb}{1.00,0.71,0.76}
\definecolor{LightSalmon1}{rgb}{1.00,0.63,0.48}
\definecolor{LightSalmon2}{rgb}{0.93,0.58,0.45}
\definecolor{LightSalmon3}{rgb}{0.80,0.51,0.38}
\definecolor{LightSalmon4}{rgb}{0.55,0.34,0.26}
\definecolor{LightSalmon}{rgb}{1.00,0.63,0.48}
\definecolor{LightSeaGreen}{rgb}{0.13,0.70,0.67}
\definecolor{LightSkyBlue1}{rgb}{0.69,0.89,1.00}
\definecolor{LightSkyBlue2}{rgb}{0.64,0.83,0.93}
\definecolor{LightSkyBlue3}{rgb}{0.55,0.71,0.80}
\definecolor{LightSkyBlue4}{rgb}{0.38,0.48,0.55}
\definecolor{LightSkyBlue}{rgb}{0.53,0.81,0.98}
\definecolor{LightSlateBlue}{rgb}{0.52,0.44,1.00}
\definecolor{LightSlateGray}{rgb}{0.47,0.53,0.60}
\definecolor{LightSlateGrey}{rgb}{0.47,0.53,0.60}
\definecolor{LightSteelBlue1}{rgb}{0.79,0.88,1.00}
\definecolor{LightSteelBlue2}{rgb}{0.74,0.82,0.93}
\definecolor{LightSteelBlue3}{rgb}{0.64,0.71,0.80}
\definecolor{LightSteelBlue4}{rgb}{0.43,0.48,0.55}
\definecolor{LightSteelBlue}{rgb}{0.69,0.77,0.87}
\definecolor{LightYellow1}{rgb}{1.00,1.00,0.88}
\definecolor{LightYellow2}{rgb}{0.93,0.93,0.82}
\definecolor{LightYellow3}{rgb}{0.80,0.80,0.71}
\definecolor{LightYellow4}{rgb}{0.55,0.55,0.48}
\definecolor{LightYellow}{rgb}{1.00,1.00,0.88}
\definecolor{LimeGreen}{rgb}{0.20,0.80,0.20}
\definecolor{MediumAquamarine}{rgb}{0.40,0.80,0.67}
\definecolor{MediumBlue}{rgb}{0.00,0.00,0.80}
\definecolor{MediumOrchid1}{rgb}{0.88,0.40,1.00}
\definecolor{MediumOrchid2}{rgb}{0.82,0.37,0.93}
\definecolor{MediumOrchid3}{rgb}{0.71,0.32,0.80}
\definecolor{MediumOrchid4}{rgb}{0.48,0.22,0.55}
\definecolor{MediumOrchid}{rgb}{0.73,0.33,0.83}
\definecolor{MediumPurple1}{rgb}{0.67,0.51,1.00}
\definecolor{MediumPurple2}{rgb}{0.62,0.47,0.93}
\definecolor{MediumPurple3}{rgb}{0.54,0.41,0.80}
\definecolor{MediumPurple4}{rgb}{0.36,0.28,0.55}
\definecolor{MediumPurple}{rgb}{0.58,0.44,0.86}
\definecolor{MediumSeaGreen}{rgb}{0.24,0.70,0.44}
\definecolor{MediumSlateBlue}{rgb}{0.48,0.41,0.93}
\definecolor{MediumSpringGreen}{rgb}{0.00,0.98,0.60}
\definecolor{MediumTurquoise}{rgb}{0.28,0.82,0.80}
\definecolor{MediumVioletRed}{rgb}{0.78,0.08,0.52}
\definecolor{MidnightBlue}{rgb}{0.10,0.10,0.44}
\definecolor{MintCream}{rgb}{0.96,1.00,0.98}
\definecolor{MistyRose1}{rgb}{1.00,0.89,0.88}
\definecolor{MistyRose2}{rgb}{0.93,0.84,0.82}
\definecolor{MistyRose3}{rgb}{0.80,0.72,0.71}
\definecolor{MistyRose4}{rgb}{0.55,0.49,0.48}
\definecolor{MistyRose}{rgb}{1.00,0.89,0.88}
\definecolor{NavajoWhite1}{rgb}{1.00,0.87,0.68}
\definecolor{NavajoWhite2}{rgb}{0.93,0.81,0.63}
\definecolor{NavajoWhite3}{rgb}{0.80,0.70,0.55}
\definecolor{NavajoWhite4}{rgb}{0.55,0.47,0.37}
\definecolor{NavajoWhite}{rgb}{1.00,0.87,0.68}
\definecolor{NavyBlue}{rgb}{0.00,0.00,0.50}
\definecolor{OldLace}{rgb}{0.99,0.96,0.90}
\definecolor{OliveDrab1}{rgb}{0.75,1.00,0.24}
\definecolor{OliveDrab2}{rgb}{0.70,0.93,0.23}
\definecolor{OliveDrab3}{rgb}{0.60,0.80,0.20}
\definecolor{OliveDrab4}{rgb}{0.41,0.55,0.13}
\definecolor{OliveDrab}{rgb}{0.42,0.56,0.14}
\definecolor{OrangeRed1}{rgb}{1.00,0.27,0.00}
\definecolor{OrangeRed2}{rgb}{0.93,0.25,0.00}
\definecolor{OrangeRed3}{rgb}{0.80,0.22,0.00}
\definecolor{OrangeRed4}{rgb}{0.55,0.15,0.00}
\definecolor{OrangeRed}{rgb}{1.00,0.27,0.00}
\definecolor{PaleGoldenrod}{rgb}{0.93,0.91,0.67}
\definecolor{PaleGreen1}{rgb}{0.60,1.00,0.60}
\definecolor{PaleGreen2}{rgb}{0.56,0.93,0.56}
\definecolor{PaleGreen3}{rgb}{0.49,0.80,0.49}
\definecolor{PaleGreen4}{rgb}{0.33,0.55,0.33}
\definecolor{PaleGreen}{rgb}{0.60,0.98,0.60}
\definecolor{PaleTurquoise1}{rgb}{0.73,1.00,1.00}
\definecolor{PaleTurquoise2}{rgb}{0.68,0.93,0.93}
\definecolor{PaleTurquoise3}{rgb}{0.59,0.80,0.80}
\definecolor{PaleTurquoise4}{rgb}{0.40,0.55,0.55}
\definecolor{PaleTurquoise}{rgb}{0.69,0.93,0.93}
\definecolor{PaleVioletRed1}{rgb}{1.00,0.51,0.67}
\definecolor{PaleVioletRed2}{rgb}{0.93,0.47,0.62}
\definecolor{PaleVioletRed3}{rgb}{0.80,0.41,0.54}
\definecolor{PaleVioletRed4}{rgb}{0.55,0.28,0.36}
\definecolor{PaleVioletRed}{rgb}{0.86,0.44,0.58}
\definecolor{PapayaWhip}{rgb}{1.00,0.94,0.84}
\definecolor{PeachPuff1}{rgb}{1.00,0.85,0.73}
\definecolor{PeachPuff2}{rgb}{0.93,0.80,0.68}
\definecolor{PeachPuff3}{rgb}{0.80,0.69,0.58}
\definecolor{PeachPuff4}{rgb}{0.55,0.47,0.40}
\definecolor{PeachPuff}{rgb}{1.00,0.85,0.73}
\definecolor{PowderBlue}{rgb}{0.69,0.88,0.90}
\definecolor{RosyBrown1}{rgb}{1.00,0.76,0.76}
\definecolor{RosyBrown2}{rgb}{0.93,0.71,0.71}
\definecolor{RosyBrown3}{rgb}{0.80,0.61,0.61}
\definecolor{RosyBrown4}{rgb}{0.55,0.41,0.41}
\definecolor{RosyBrown}{rgb}{0.74,0.56,0.56}
\definecolor{RoyalBlue1}{rgb}{0.28,0.46,1.00}
\definecolor{RoyalBlue2}{rgb}{0.26,0.43,0.93}
\definecolor{RoyalBlue3}{rgb}{0.23,0.37,0.80}
\definecolor{RoyalBlue4}{rgb}{0.15,0.25,0.55}
\definecolor{RoyalBlue}{rgb}{0.25,0.41,0.88}
\definecolor{SaddleBrown}{rgb}{0.55,0.27,0.07}
\definecolor{SandyBrown}{rgb}{0.96,0.64,0.38}
\definecolor{SeaGreen1}{rgb}{0.33,1.00,0.62}
\definecolor{SeaGreen2}{rgb}{0.31,0.93,0.58}
\definecolor{SeaGreen3}{rgb}{0.26,0.80,0.50}
\definecolor{SeaGreen4}{rgb}{0.18,0.55,0.34}
\definecolor{SeaGreen}{rgb}{0.18,0.55,0.34}
\definecolor{SkyBlue1}{rgb}{0.53,0.81,1.00}
\definecolor{SkyBlue2}{rgb}{0.49,0.75,0.93}
\definecolor{SkyBlue3}{rgb}{0.42,0.65,0.80}
\definecolor{SkyBlue4}{rgb}{0.29,0.44,0.55}
\definecolor{SkyBlue}{rgb}{0.53,0.81,0.92}
\definecolor{SlateBlue1}{rgb}{0.51,0.44,1.00}
\definecolor{SlateBlue2}{rgb}{0.48,0.40,0.93}
\definecolor{SlateBlue3}{rgb}{0.41,0.35,0.80}
\definecolor{SlateBlue4}{rgb}{0.28,0.24,0.55}
\definecolor{SlateBlue}{rgb}{0.42,0.35,0.80}
\definecolor{SlateGray1}{rgb}{0.78,0.89,1.00}
\definecolor{SlateGray2}{rgb}{0.73,0.83,0.93}
\definecolor{SlateGray3}{rgb}{0.62,0.71,0.80}
\definecolor{SlateGray4}{rgb}{0.42,0.48,0.55}
\definecolor{SlateGray}{rgb}{0.44,0.50,0.56}
\definecolor{SlateGrey}{rgb}{0.44,0.50,0.56}
\definecolor{SpringGreen1}{rgb}{0.00,1.00,0.50}
\definecolor{SpringGreen2}{rgb}{0.00,0.93,0.46}
\definecolor{SpringGreen3}{rgb}{0.00,0.80,0.40}
\definecolor{SpringGreen4}{rgb}{0.00,0.55,0.27}
\definecolor{SpringGreen}{rgb}{0.00,1.00,0.50}
\definecolor{SteelBlue1}{rgb}{0.39,0.72,1.00}
\definecolor{SteelBlue2}{rgb}{0.36,0.67,0.93}
\definecolor{SteelBlue3}{rgb}{0.31,0.58,0.80}
\definecolor{SteelBlue4}{rgb}{0.21,0.39,0.55}
\definecolor{SteelBlue}{rgb}{0.27,0.51,0.71}
\definecolor{VioletRed1}{rgb}{1.00,0.24,0.59}
\definecolor{VioletRed2}{rgb}{0.93,0.23,0.55}
\definecolor{VioletRed3}{rgb}{0.80,0.20,0.47}
\definecolor{VioletRed4}{rgb}{0.55,0.13,0.32}
\definecolor{VioletRed}{rgb}{0.82,0.13,0.56}
\definecolor{WhiteSmoke}{rgb}{0.96,0.96,0.96}
\definecolor{YellowGreen}{rgb}{0.60,0.80,0.20}
\definecolor{aliceblue}{rgb}{0.94,0.97,1.00}
\definecolor{antiquewhite}{rgb}{0.98,0.92,0.84}
\definecolor{aquamarine1}{rgb}{0.50,1.00,0.83}
\definecolor{aquamarine2}{rgb}{0.46,0.93,0.78}
\definecolor{aquamarine3}{rgb}{0.40,0.80,0.67}
\definecolor{aquamarine4}{rgb}{0.27,0.55,0.45}
\definecolor{aquamarine}{rgb}{0.50,1.00,0.83}
\definecolor{azure1}{rgb}{0.94,1.00,1.00}
\definecolor{azure2}{rgb}{0.88,0.93,0.93}
\definecolor{azure3}{rgb}{0.76,0.80,0.80}
\definecolor{azure4}{rgb}{0.51,0.55,0.55}
\definecolor{azure}{rgb}{0.94,1.00,1.00}
\definecolor{beige}{rgb}{0.96,0.96,0.86}
\definecolor{bisque1}{rgb}{1.00,0.89,0.77}
\definecolor{bisque2}{rgb}{0.93,0.84,0.72}
\definecolor{bisque3}{rgb}{0.80,0.72,0.62}
\definecolor{bisque4}{rgb}{0.55,0.49,0.42}
\definecolor{bisque}{rgb}{1.00,0.89,0.77}
\definecolor{black}{rgb}{0.00,0.00,0.00}
\definecolor{blanchedalmond}{rgb}{1.00,0.92,0.80}
\definecolor{blue1}{rgb}{0.00,0.00,1.00}
\definecolor{blue2}{rgb}{0.00,0.00,0.93}
\definecolor{blue3}{rgb}{0.00,0.00,0.80}
\definecolor{blue4}{rgb}{0.00,0.00,0.55}
\definecolor{blueviolet}{rgb}{0.54,0.17,0.89}
\definecolor{blue}{rgb}{0.00,0.00,1.00}
\definecolor{brown1}{rgb}{1.00,0.25,0.25}
\definecolor{brown2}{rgb}{0.93,0.23,0.23}
\definecolor{brown3}{rgb}{0.80,0.20,0.20}
\definecolor{brown4}{rgb}{0.55,0.14,0.14}
\definecolor{brown}{rgb}{0.65,0.16,0.16}
\definecolor{burlywood1}{rgb}{1.00,0.83,0.61}
\definecolor{burlywood2}{rgb}{0.93,0.77,0.57}
\definecolor{burlywood3}{rgb}{0.80,0.67,0.49}
\definecolor{burlywood4}{rgb}{0.55,0.45,0.33}
\definecolor{burlywood}{rgb}{0.87,0.72,0.53}
\definecolor{cadetblue}{rgb}{0.37,0.62,0.63}
\definecolor{chartreuse1}{rgb}{0.50,1.00,0.00}
\definecolor{chartreuse2}{rgb}{0.46,0.93,0.00}
\definecolor{chartreuse3}{rgb}{0.40,0.80,0.00}
\definecolor{chartreuse4}{rgb}{0.27,0.55,0.00}
\definecolor{chartreuse}{rgb}{0.50,1.00,0.00}
\definecolor{chocolate1}{rgb}{1.00,0.50,0.14}
\definecolor{chocolate2}{rgb}{0.93,0.46,0.13}
\definecolor{chocolate3}{rgb}{0.80,0.40,0.11}
\definecolor{chocolate4}{rgb}{0.55,0.27,0.07}
\definecolor{chocolate}{rgb}{0.82,0.41,0.12}
\definecolor{coral1}{rgb}{1.00,0.45,0.34}
\definecolor{coral2}{rgb}{0.93,0.42,0.31}
\definecolor{coral3}{rgb}{0.80,0.36,0.27}
\definecolor{coral4}{rgb}{0.55,0.24,0.18}
\definecolor{coral}{rgb}{1.00,0.50,0.31}
\definecolor{cornflowerblue}{rgb}{0.39,0.58,0.93}
\definecolor{cornsilk1}{rgb}{1.00,0.97,0.86}
\definecolor{cornsilk2}{rgb}{0.93,0.91,0.80}
\definecolor{cornsilk3}{rgb}{0.80,0.78,0.69}
\definecolor{cornsilk4}{rgb}{0.55,0.53,0.47}
\definecolor{cornsilk}{rgb}{1.00,0.97,0.86}
\definecolor{cyan1}{rgb}{0.00,1.00,1.00}
\definecolor{cyan2}{rgb}{0.00,0.93,0.93}
\definecolor{cyan3}{rgb}{0.00,0.80,0.80}
\definecolor{cyan4}{rgb}{0.00,0.55,0.55}
\definecolor{cyan}{rgb}{0.00,1.00,1.00}
\definecolor{darkblue}{rgb}{0.00,0.00,0.55}
\definecolor{darkcyan}{rgb}{0.00,0.55,0.55}
\definecolor{darkgoldenrod}{rgb}{0.72,0.53,0.04}
\definecolor{darkgray}{rgb}{0.66,0.66,0.66}
\definecolor{darkgreen}{rgb}{0.00,0.39,0.00}
\definecolor{darkgrey}{rgb}{0.66,0.66,0.66}
\definecolor{darkkhaki}{rgb}{0.74,0.72,0.42}
\definecolor{darkmagenta}{rgb}{0.55,0.00,0.55}
\definecolor{darkolive}{rgb}{0.33,0.42,0.18}
\definecolor{darkorange}{rgb}{1.00,0.55,0.00}
\definecolor{darkorchid}{rgb}{0.60,0.20,0.80}
\definecolor{darkred}{rgb}{0.55,0.00,0.00}
\definecolor{darksalmon}{rgb}{0.91,0.59,0.48}
\definecolor{darksea}{rgb}{0.56,0.74,0.56}
\definecolor{darkslate}{rgb}{0.18,0.31,0.31}
\definecolor{darkslate}{rgb}{0.18,0.31,0.31}
\definecolor{darkslate}{rgb}{0.28,0.24,0.55}
\definecolor{darkturquoise}{rgb}{0.00,0.81,0.82}
\definecolor{darkviolet}{rgb}{0.58,0.00,0.83}
\definecolor{deeppink}{rgb}{1.00,0.08,0.58}
\definecolor{deepsky}{rgb}{0.00,0.75,1.00}
\definecolor{dimgray}{rgb}{0.41,0.41,0.41}
\definecolor{dimgrey}{rgb}{0.41,0.41,0.41}
\definecolor{dodgerblue}{rgb}{0.12,0.56,1.00}
\definecolor{firebrick1}{rgb}{1.00,0.19,0.19}
\definecolor{firebrick2}{rgb}{0.93,0.17,0.17}
\definecolor{firebrick3}{rgb}{0.80,0.15,0.15}
\definecolor{firebrick4}{rgb}{0.55,0.10,0.10}
\definecolor{firebrick}{rgb}{0.70,0.13,0.13}
\definecolor{floralwhite}{rgb}{1.00,0.98,0.94}
\definecolor{forestgreen}{rgb}{0.13,0.55,0.13}
\definecolor{gainsboro}{rgb}{0.86,0.86,0.86}
\definecolor{ghostwhite}{rgb}{0.97,0.97,1.00}
\definecolor{gold1}{rgb}{1.00,0.84,0.00}
\definecolor{gold2}{rgb}{0.93,0.79,0.00}
\definecolor{gold3}{rgb}{0.80,0.68,0.00}
\definecolor{gold4}{rgb}{0.55,0.46,0.00}
\definecolor{goldenrod1}{rgb}{1.00,0.76,0.15}
\definecolor{goldenrod2}{rgb}{0.93,0.71,0.13}
\definecolor{goldenrod3}{rgb}{0.80,0.61,0.11}
\definecolor{goldenrod4}{rgb}{0.55,0.41,0.08}
\definecolor{goldenrod}{rgb}{0.85,0.65,0.13}
\definecolor{gold}{rgb}{1.00,0.84,0.00}
\definecolor{gray0}{rgb}{0.00,0.00,0.00}
\definecolor{gray100}{rgb}{1.00,1.00,1.00}
\definecolor{gray10}{rgb}{0.10,0.10,0.10}
\definecolor{gray11}{rgb}{0.11,0.11,0.11}
\definecolor{gray12}{rgb}{0.12,0.12,0.12}
\definecolor{gray13}{rgb}{0.13,0.13,0.13}
\definecolor{gray14}{rgb}{0.14,0.14,0.14}
\definecolor{gray15}{rgb}{0.15,0.15,0.15}
\definecolor{gray16}{rgb}{0.16,0.16,0.16}
\definecolor{gray17}{rgb}{0.17,0.17,0.17}
\definecolor{gray18}{rgb}{0.18,0.18,0.18}
\definecolor{gray19}{rgb}{0.19,0.19,0.19}
\definecolor{gray1}{rgb}{0.01,0.01,0.01}
\definecolor{gray20}{rgb}{0.20,0.20,0.20}
\definecolor{gray21}{rgb}{0.21,0.21,0.21}
\definecolor{gray22}{rgb}{0.22,0.22,0.22}
\definecolor{gray23}{rgb}{0.23,0.23,0.23}
\definecolor{gray24}{rgb}{0.24,0.24,0.24}
\definecolor{gray25}{rgb}{0.25,0.25,0.25}
\definecolor{gray26}{rgb}{0.26,0.26,0.26}
\definecolor{gray27}{rgb}{0.27,0.27,0.27}
\definecolor{gray28}{rgb}{0.28,0.28,0.28}
\definecolor{gray29}{rgb}{0.29,0.29,0.29}
\definecolor{gray2}{rgb}{0.02,0.02,0.02}
\definecolor{gray30}{rgb}{0.30,0.30,0.30}
\definecolor{gray31}{rgb}{0.31,0.31,0.31}
\definecolor{gray32}{rgb}{0.32,0.32,0.32}
\definecolor{gray33}{rgb}{0.33,0.33,0.33}
\definecolor{gray34}{rgb}{0.34,0.34,0.34}
\definecolor{gray35}{rgb}{0.35,0.35,0.35}
\definecolor{gray36}{rgb}{0.36,0.36,0.36}
\definecolor{gray37}{rgb}{0.37,0.37,0.37}
\definecolor{gray38}{rgb}{0.38,0.38,0.38}
\definecolor{gray39}{rgb}{0.39,0.39,0.39}
\definecolor{gray3}{rgb}{0.03,0.03,0.03}
\definecolor{gray40}{rgb}{0.40,0.40,0.40}
\definecolor{gray41}{rgb}{0.41,0.41,0.41}
\definecolor{gray42}{rgb}{0.42,0.42,0.42}
\definecolor{gray43}{rgb}{0.43,0.43,0.43}
\definecolor{gray44}{rgb}{0.44,0.44,0.44}
\definecolor{gray45}{rgb}{0.45,0.45,0.45}
\definecolor{gray46}{rgb}{0.46,0.46,0.46}
\definecolor{gray47}{rgb}{0.47,0.47,0.47}
\definecolor{gray48}{rgb}{0.48,0.48,0.48}
\definecolor{gray49}{rgb}{0.49,0.49,0.49}
\definecolor{gray4}{rgb}{0.04,0.04,0.04}
\definecolor{gray50}{rgb}{0.50,0.50,0.50}
\definecolor{gray51}{rgb}{0.51,0.51,0.51}
\definecolor{gray52}{rgb}{0.52,0.52,0.52}
\definecolor{gray53}{rgb}{0.53,0.53,0.53}
\definecolor{gray54}{rgb}{0.54,0.54,0.54}
\definecolor{gray55}{rgb}{0.55,0.55,0.55}
\definecolor{gray56}{rgb}{0.56,0.56,0.56}
\definecolor{gray57}{rgb}{0.57,0.57,0.57}
\definecolor{gray58}{rgb}{0.58,0.58,0.58}
\definecolor{gray59}{rgb}{0.59,0.59,0.59}
\definecolor{gray5}{rgb}{0.05,0.05,0.05}
\definecolor{gray60}{rgb}{0.60,0.60,0.60}
\definecolor{gray61}{rgb}{0.61,0.61,0.61}
\definecolor{gray62}{rgb}{0.62,0.62,0.62}
\definecolor{gray63}{rgb}{0.63,0.63,0.63}
\definecolor{gray64}{rgb}{0.64,0.64,0.64}
\definecolor{gray65}{rgb}{0.65,0.65,0.65}
\definecolor{gray66}{rgb}{0.66,0.66,0.66}
\definecolor{gray67}{rgb}{0.67,0.67,0.67}
\definecolor{gray68}{rgb}{0.68,0.68,0.68}
\definecolor{gray69}{rgb}{0.69,0.69,0.69}
\definecolor{gray6}{rgb}{0.06,0.06,0.06}
\definecolor{gray70}{rgb}{0.70,0.70,0.70}
\definecolor{gray71}{rgb}{0.71,0.71,0.71}
\definecolor{gray72}{rgb}{0.72,0.72,0.72}
\definecolor{gray73}{rgb}{0.73,0.73,0.73}
\definecolor{gray74}{rgb}{0.74,0.74,0.74}
\definecolor{gray75}{rgb}{0.75,0.75,0.75}
\definecolor{gray76}{rgb}{0.76,0.76,0.76}
\definecolor{gray77}{rgb}{0.77,0.77,0.77}
\definecolor{gray78}{rgb}{0.78,0.78,0.78}
\definecolor{gray79}{rgb}{0.79,0.79,0.79}
\definecolor{gray7}{rgb}{0.07,0.07,0.07}
\definecolor{gray80}{rgb}{0.80,0.80,0.80}
\definecolor{gray81}{rgb}{0.81,0.81,0.81}
\definecolor{gray82}{rgb}{0.82,0.82,0.82}
\definecolor{gray83}{rgb}{0.83,0.83,0.83}
\definecolor{gray84}{rgb}{0.84,0.84,0.84}
\definecolor{gray85}{rgb}{0.85,0.85,0.85}
\definecolor{gray86}{rgb}{0.86,0.86,0.86}
\definecolor{gray87}{rgb}{0.87,0.87,0.87}
\definecolor{gray88}{rgb}{0.88,0.88,0.88}
\definecolor{gray89}{rgb}{0.89,0.89,0.89}
\definecolor{gray8}{rgb}{0.08,0.08,0.08}
\definecolor{gray90}{rgb}{0.90,0.90,0.90}
\definecolor{gray91}{rgb}{0.91,0.91,0.91}
\definecolor{gray92}{rgb}{0.92,0.92,0.92}
\definecolor{gray93}{rgb}{0.93,0.93,0.93}
\definecolor{gray94}{rgb}{0.94,0.94,0.94}
\definecolor{gray95}{rgb}{0.95,0.95,0.95}
\definecolor{gray96}{rgb}{0.96,0.96,0.96}
\definecolor{gray97}{rgb}{0.97,0.97,0.97}
\definecolor{gray98}{rgb}{0.98,0.98,0.98}
\definecolor{gray99}{rgb}{0.99,0.99,0.99}
\definecolor{gray9}{rgb}{0.09,0.09,0.09}
\definecolor{gray}{rgb}{0.75,0.75,0.75}
\definecolor{green1}{rgb}{0.00,1.00,0.00}
\definecolor{green2}{rgb}{0.00,0.93,0.00}
\definecolor{green3}{rgb}{0.00,0.80,0.00}
\definecolor{green4}{rgb}{0.00,0.55,0.00}
\definecolor{greenyellow}{rgb}{0.68,1.00,0.18}
\definecolor{green}{rgb}{0.00,1.00,0.00}
\definecolor{grey0}{rgb}{0.00,0.00,0.00}
\definecolor{grey100}{rgb}{1.00,1.00,1.00}
\definecolor{grey10}{rgb}{0.10,0.10,0.10}
\definecolor{grey11}{rgb}{0.11,0.11,0.11}
\definecolor{grey12}{rgb}{0.12,0.12,0.12}
\definecolor{grey13}{rgb}{0.13,0.13,0.13}
\definecolor{grey14}{rgb}{0.14,0.14,0.14}
\definecolor{grey15}{rgb}{0.15,0.15,0.15}
\definecolor{grey16}{rgb}{0.16,0.16,0.16}
\definecolor{grey17}{rgb}{0.17,0.17,0.17}
\definecolor{grey18}{rgb}{0.18,0.18,0.18}
\definecolor{grey19}{rgb}{0.19,0.19,0.19}
\definecolor{grey1}{rgb}{0.01,0.01,0.01}
\definecolor{grey20}{rgb}{0.20,0.20,0.20}
\definecolor{grey21}{rgb}{0.21,0.21,0.21}
\definecolor{grey22}{rgb}{0.22,0.22,0.22}
\definecolor{grey23}{rgb}{0.23,0.23,0.23}
\definecolor{grey24}{rgb}{0.24,0.24,0.24}
\definecolor{grey25}{rgb}{0.25,0.25,0.25}
\definecolor{grey26}{rgb}{0.26,0.26,0.26}
\definecolor{grey27}{rgb}{0.27,0.27,0.27}
\definecolor{grey28}{rgb}{0.28,0.28,0.28}
\definecolor{grey29}{rgb}{0.29,0.29,0.29}
\definecolor{grey2}{rgb}{0.02,0.02,0.02}
\definecolor{grey30}{rgb}{0.30,0.30,0.30}
\definecolor{grey31}{rgb}{0.31,0.31,0.31}
\definecolor{grey32}{rgb}{0.32,0.32,0.32}
\definecolor{grey33}{rgb}{0.33,0.33,0.33}
\definecolor{grey34}{rgb}{0.34,0.34,0.34}
\definecolor{grey35}{rgb}{0.35,0.35,0.35}
\definecolor{grey36}{rgb}{0.36,0.36,0.36}
\definecolor{grey37}{rgb}{0.37,0.37,0.37}
\definecolor{grey38}{rgb}{0.38,0.38,0.38}
\definecolor{grey39}{rgb}{0.39,0.39,0.39}
\definecolor{grey3}{rgb}{0.03,0.03,0.03}
\definecolor{grey40}{rgb}{0.40,0.40,0.40}
\definecolor{grey41}{rgb}{0.41,0.41,0.41}
\definecolor{grey42}{rgb}{0.42,0.42,0.42}
\definecolor{grey43}{rgb}{0.43,0.43,0.43}
\definecolor{grey44}{rgb}{0.44,0.44,0.44}
\definecolor{grey45}{rgb}{0.45,0.45,0.45}
\definecolor{grey46}{rgb}{0.46,0.46,0.46}
\definecolor{grey47}{rgb}{0.47,0.47,0.47}
\definecolor{grey48}{rgb}{0.48,0.48,0.48}
\definecolor{grey49}{rgb}{0.49,0.49,0.49}
\definecolor{grey4}{rgb}{0.04,0.04,0.04}
\definecolor{grey50}{rgb}{0.50,0.50,0.50}
\definecolor{grey51}{rgb}{0.51,0.51,0.51}
\definecolor{grey52}{rgb}{0.52,0.52,0.52}
\definecolor{grey53}{rgb}{0.53,0.53,0.53}
\definecolor{grey54}{rgb}{0.54,0.54,0.54}
\definecolor{grey55}{rgb}{0.55,0.55,0.55}
\definecolor{grey56}{rgb}{0.56,0.56,0.56}
\definecolor{grey57}{rgb}{0.57,0.57,0.57}
\definecolor{grey58}{rgb}{0.58,0.58,0.58}
\definecolor{grey59}{rgb}{0.59,0.59,0.59}
\definecolor{grey5}{rgb}{0.05,0.05,0.05}
\definecolor{grey60}{rgb}{0.60,0.60,0.60}
\definecolor{grey61}{rgb}{0.61,0.61,0.61}
\definecolor{grey62}{rgb}{0.62,0.62,0.62}
\definecolor{grey63}{rgb}{0.63,0.63,0.63}
\definecolor{grey64}{rgb}{0.64,0.64,0.64}
\definecolor{grey65}{rgb}{0.65,0.65,0.65}
\definecolor{grey66}{rgb}{0.66,0.66,0.66}
\definecolor{grey67}{rgb}{0.67,0.67,0.67}
\definecolor{grey68}{rgb}{0.68,0.68,0.68}
\definecolor{grey69}{rgb}{0.69,0.69,0.69}
\definecolor{grey6}{rgb}{0.06,0.06,0.06}
\definecolor{grey70}{rgb}{0.70,0.70,0.70}
\definecolor{grey71}{rgb}{0.71,0.71,0.71}
\definecolor{grey72}{rgb}{0.72,0.72,0.72}
\definecolor{grey73}{rgb}{0.73,0.73,0.73}
\definecolor{grey74}{rgb}{0.74,0.74,0.74}
\definecolor{grey75}{rgb}{0.75,0.75,0.75}
\definecolor{grey76}{rgb}{0.76,0.76,0.76}
\definecolor{grey77}{rgb}{0.77,0.77,0.77}
\definecolor{grey78}{rgb}{0.78,0.78,0.78}
\definecolor{grey79}{rgb}{0.79,0.79,0.79}
\definecolor{grey7}{rgb}{0.07,0.07,0.07}
\definecolor{grey80}{rgb}{0.80,0.80,0.80}
\definecolor{grey81}{rgb}{0.81,0.81,0.81}
\definecolor{grey82}{rgb}{0.82,0.82,0.82}
\definecolor{grey83}{rgb}{0.83,0.83,0.83}
\definecolor{grey84}{rgb}{0.84,0.84,0.84}
\definecolor{grey85}{rgb}{0.85,0.85,0.85}
\definecolor{grey86}{rgb}{0.86,0.86,0.86}
\definecolor{grey87}{rgb}{0.87,0.87,0.87}
\definecolor{grey88}{rgb}{0.88,0.88,0.88}
\definecolor{grey89}{rgb}{0.89,0.89,0.89}
\definecolor{grey8}{rgb}{0.08,0.08,0.08}
\definecolor{grey90}{rgb}{0.90,0.90,0.90}
\definecolor{grey91}{rgb}{0.91,0.91,0.91}
\definecolor{grey92}{rgb}{0.92,0.92,0.92}
\definecolor{grey93}{rgb}{0.93,0.93,0.93}
\definecolor{grey94}{rgb}{0.94,0.94,0.94}
\definecolor{grey95}{rgb}{0.95,0.95,0.95}
\definecolor{grey96}{rgb}{0.96,0.96,0.96}
\definecolor{grey97}{rgb}{0.97,0.97,0.97}
\definecolor{grey98}{rgb}{0.98,0.98,0.98}
\definecolor{grey99}{rgb}{0.99,0.99,0.99}
\definecolor{grey9}{rgb}{0.09,0.09,0.09}
\definecolor{grey}{rgb}{0.75,0.75,0.75}
\definecolor{honeydew1}{rgb}{0.94,1.00,0.94}
\definecolor{honeydew2}{rgb}{0.88,0.93,0.88}
\definecolor{honeydew3}{rgb}{0.76,0.80,0.76}
\definecolor{honeydew4}{rgb}{0.51,0.55,0.51}
\definecolor{honeydew}{rgb}{0.94,1.00,0.94}
\definecolor{hotpink}{rgb}{1.00,0.41,0.71}
\definecolor{indianred}{rgb}{0.80,0.36,0.36}
\definecolor{ivory1}{rgb}{1.00,1.00,0.94}
\definecolor{ivory2}{rgb}{0.93,0.93,0.88}
\definecolor{ivory3}{rgb}{0.80,0.80,0.76}
\definecolor{ivory4}{rgb}{0.55,0.55,0.51}
\definecolor{ivory}{rgb}{1.00,1.00,0.94}
\definecolor{khaki1}{rgb}{1.00,0.96,0.56}
\definecolor{khaki2}{rgb}{0.93,0.90,0.52}
\definecolor{khaki3}{rgb}{0.80,0.78,0.45}
\definecolor{khaki4}{rgb}{0.55,0.53,0.31}
\definecolor{khaki}{rgb}{0.94,0.90,0.55}
\definecolor{lavenderblush}{rgb}{1.00,0.94,0.96}
\definecolor{lavender}{rgb}{0.90,0.90,0.98}
\definecolor{lawngreen}{rgb}{0.49,0.99,0.00}
\definecolor{lemonchiffon}{rgb}{1.00,0.98,0.80}
\definecolor{lightblue}{rgb}{0.68,0.85,0.90}
\definecolor{lightcoral}{rgb}{0.94,0.50,0.50}
\definecolor{lightcyan}{rgb}{0.88,1.00,1.00}
\definecolor{lightgoldenrod}{rgb}{0.93,0.87,0.51}
\definecolor{lightgoldenrod}{rgb}{0.98,0.98,0.82}
\definecolor{lightgray}{rgb}{0.83,0.83,0.83}
\definecolor{lightgreen}{rgb}{0.56,0.93,0.56}
\definecolor{lightgrey}{rgb}{0.83,0.83,0.83}
\definecolor{lightpink}{rgb}{1.00,0.71,0.76}
\definecolor{lightsalmon}{rgb}{1.00,0.63,0.48}
\definecolor{lightsea}{rgb}{0.13,0.70,0.67}
\definecolor{lightsky}{rgb}{0.53,0.81,0.98}
\definecolor{lightslate}{rgb}{0.47,0.53,0.60}
\definecolor{lightslate}{rgb}{0.47,0.53,0.60}
\definecolor{lightslate}{rgb}{0.52,0.44,1.00}
\definecolor{lightsteel}{rgb}{0.69,0.77,0.87}
\definecolor{lightyellow}{rgb}{1.00,1.00,0.88}
\definecolor{limegreen}{rgb}{0.20,0.80,0.20}
\definecolor{linen}{rgb}{0.98,0.94,0.90}
\definecolor{magenta1}{rgb}{1.00,0.00,1.00}
\definecolor{magenta2}{rgb}{0.93,0.00,0.93}
\definecolor{magenta3}{rgb}{0.80,0.00,0.80}
\definecolor{magenta4}{rgb}{0.55,0.00,0.55}
\definecolor{magenta}{rgb}{1.00,0.00,1.00}
\definecolor{maroon1}{rgb}{1.00,0.20,0.70}
\definecolor{maroon2}{rgb}{0.93,0.19,0.65}
\definecolor{maroon3}{rgb}{0.80,0.16,0.56}
\definecolor{maroon4}{rgb}{0.55,0.11,0.38}
\definecolor{maroon}{rgb}{0.69,0.19,0.38}
\definecolor{mediumaquamarine}{rgb}{0.40,0.80,0.67}
\definecolor{mediumblue}{rgb}{0.00,0.00,0.80}
\definecolor{mediumorchid}{rgb}{0.73,0.33,0.83}
\definecolor{mediumpurple}{rgb}{0.58,0.44,0.86}
\definecolor{mediumsea}{rgb}{0.24,0.70,0.44}
\definecolor{mediumslate}{rgb}{0.48,0.41,0.93}
\definecolor{mediumspring}{rgb}{0.00,0.98,0.60}
\definecolor{mediumturquoise}{rgb}{0.28,0.82,0.80}
\definecolor{mediumviolet}{rgb}{0.78,0.08,0.52}
\definecolor{midnightblue}{rgb}{0.10,0.10,0.44}
\definecolor{mintcream}{rgb}{0.96,1.00,0.98}
\definecolor{mistyrose}{rgb}{1.00,0.89,0.88}
\definecolor{moccasin}{rgb}{1.00,0.89,0.71}
\definecolor{navajowhite}{rgb}{1.00,0.87,0.68}
\definecolor{navyblue}{rgb}{0.00,0.00,0.50}
\definecolor{navy}{rgb}{0.00,0.00,0.50}
\definecolor{oldlace}{rgb}{0.99,0.96,0.90}
\definecolor{olivedrab}{rgb}{0.42,0.56,0.14}
\definecolor{orange1}{rgb}{1.00,0.65,0.00}
\definecolor{orange2}{rgb}{0.93,0.60,0.00}
\definecolor{orange3}{rgb}{0.80,0.52,0.00}
\definecolor{orange4}{rgb}{0.55,0.35,0.00}
\definecolor{orangered}{rgb}{1.00,0.27,0.00}
\definecolor{orange}{rgb}{1.00,0.65,0.00}
\definecolor{orchid1}{rgb}{1.00,0.51,0.98}
\definecolor{orchid2}{rgb}{0.93,0.48,0.91}
\definecolor{orchid3}{rgb}{0.80,0.41,0.79}
\definecolor{orchid4}{rgb}{0.55,0.28,0.54}
\definecolor{orchid}{rgb}{0.85,0.44,0.84}
\definecolor{palegoldenrod}{rgb}{0.93,0.91,0.67}
\definecolor{palegreen}{rgb}{0.60,0.98,0.60}
\definecolor{paleturquoise}{rgb}{0.69,0.93,0.93}
\definecolor{paleviolet}{rgb}{0.86,0.44,0.58}
\definecolor{papayawhip}{rgb}{1.00,0.94,0.84}
\definecolor{peachpuff}{rgb}{1.00,0.85,0.73}
\definecolor{peru}{rgb}{0.80,0.52,0.25}
\definecolor{pink1}{rgb}{1.00,0.71,0.77}
\definecolor{pink2}{rgb}{0.93,0.66,0.72}
\definecolor{pink3}{rgb}{0.80,0.57,0.62}
\definecolor{pink4}{rgb}{0.55,0.39,0.42}
\definecolor{pink}{rgb}{1.00,0.75,0.80}
\definecolor{plum1}{rgb}{1.00,0.73,1.00}
\definecolor{plum2}{rgb}{0.93,0.68,0.93}
\definecolor{plum3}{rgb}{0.80,0.59,0.80}
\definecolor{plum4}{rgb}{0.55,0.40,0.55}
\definecolor{plum}{rgb}{0.87,0.63,0.87}
\definecolor{powderblue}{rgb}{0.69,0.88,0.90}
\definecolor{purple1}{rgb}{0.61,0.19,1.00}
\definecolor{purple2}{rgb}{0.57,0.17,0.93}
\definecolor{purple3}{rgb}{0.49,0.15,0.80}
\definecolor{purple4}{rgb}{0.33,0.10,0.55}
\definecolor{purple}{rgb}{0.63,0.13,0.94}
\definecolor{red1}{rgb}{1.00,0.00,0.00}
\definecolor{red2}{rgb}{0.93,0.00,0.00}
\definecolor{red3}{rgb}{0.80,0.00,0.00}
\definecolor{red4}{rgb}{0.55,0.00,0.00}
\definecolor{red}{rgb}{1.00,0.00,0.00}
\definecolor{rosybrown}{rgb}{0.74,0.56,0.56}
\definecolor{royalblue}{rgb}{0.25,0.41,0.88}
\definecolor{saddlebrown}{rgb}{0.55,0.27,0.07}
\definecolor{salmon1}{rgb}{1.00,0.55,0.41}
\definecolor{salmon2}{rgb}{0.93,0.51,0.38}
\definecolor{salmon3}{rgb}{0.80,0.44,0.33}
\definecolor{salmon4}{rgb}{0.55,0.30,0.22}
\definecolor{salmon}{rgb}{0.98,0.50,0.45}
\definecolor{sandybrown}{rgb}{0.96,0.64,0.38}
\definecolor{seagreen}{rgb}{0.18,0.55,0.34}
\definecolor{seashell1}{rgb}{1.00,0.96,0.93}
\definecolor{seashell2}{rgb}{0.93,0.90,0.87}
\definecolor{seashell3}{rgb}{0.80,0.77,0.75}
\definecolor{seashell4}{rgb}{0.55,0.53,0.51}
\definecolor{seashell}{rgb}{1.00,0.96,0.93}
\definecolor{sienna1}{rgb}{1.00,0.51,0.28}
\definecolor{sienna2}{rgb}{0.93,0.47,0.26}
\definecolor{sienna3}{rgb}{0.80,0.41,0.22}
\definecolor{sienna4}{rgb}{0.55,0.28,0.15}
\definecolor{sienna}{rgb}{0.63,0.32,0.18}
\definecolor{skyblue}{rgb}{0.53,0.81,0.92}
\definecolor{slateblue}{rgb}{0.42,0.35,0.80}
\definecolor{slategray}{rgb}{0.44,0.50,0.56}
\definecolor{slategrey}{rgb}{0.44,0.50,0.56}
\definecolor{snow1}{rgb}{1.00,0.98,0.98}
\definecolor{snow2}{rgb}{0.93,0.91,0.91}
\definecolor{snow3}{rgb}{0.80,0.79,0.79}
\definecolor{snow4}{rgb}{0.55,0.54,0.54}
\definecolor{snow}{rgb}{1.00,0.98,0.98}
\definecolor{springgreen}{rgb}{0.00,1.00,0.50}
\definecolor{steelblue}{rgb}{0.27,0.51,0.71}
\definecolor{tan1}{rgb}{1.00,0.65,0.31}
\definecolor{tan2}{rgb}{0.93,0.60,0.29}
\definecolor{tan3}{rgb}{0.80,0.52,0.25}
\definecolor{tan4}{rgb}{0.55,0.35,0.17}
\definecolor{tan}{rgb}{0.82,0.71,0.55}
\definecolor{thistle1}{rgb}{1.00,0.88,1.00}
\definecolor{thistle2}{rgb}{0.93,0.82,0.93}
\definecolor{thistle3}{rgb}{0.80,0.71,0.80}
\definecolor{thistle4}{rgb}{0.55,0.48,0.55}
\definecolor{thistle}{rgb}{0.85,0.75,0.85}
\definecolor{tomato1}{rgb}{1.00,0.39,0.28}
\definecolor{tomato2}{rgb}{0.93,0.36,0.26}
\definecolor{tomato3}{rgb}{0.80,0.31,0.22}
\definecolor{tomato4}{rgb}{0.55,0.21,0.15}
\definecolor{tomato}{rgb}{1.00,0.39,0.28}
\definecolor{turquoise1}{rgb}{0.00,0.96,1.00}
\definecolor{turquoise2}{rgb}{0.00,0.90,0.93}
\definecolor{turquoise3}{rgb}{0.00,0.77,0.80}
\definecolor{turquoise4}{rgb}{0.00,0.53,0.55}
\definecolor{turquoise}{rgb}{0.25,0.88,0.82}
\definecolor{violetred}{rgb}{0.82,0.13,0.56}
\definecolor{violet}{rgb}{0.93,0.51,0.93}
\definecolor{wheat1}{rgb}{1.00,0.91,0.73}
\definecolor{wheat2}{rgb}{0.93,0.85,0.68}
\definecolor{wheat3}{rgb}{0.80,0.73,0.59}
\definecolor{wheat4}{rgb}{0.55,0.49,0.40}
\definecolor{wheat}{rgb}{0.96,0.87,0.70}
\definecolor{whitesmoke}{rgb}{0.96,0.96,0.96}
\definecolor{white}{rgb}{1.00,1.00,1.00}
\definecolor{yellow1}{rgb}{1.00,1.00,0.00}
\definecolor{yellow2}{rgb}{0.93,0.93,0.00}
\definecolor{yellow3}{rgb}{0.80,0.80,0.00}
\definecolor{yellow4}{rgb}{0.55,0.55,0.00}
\definecolor{yellowgreen}{rgb}{0.60,0.80,0.20}
\definecolor{yellow}{rgb}{1.00,1.00,0.00}

\usepackage[usenames,dvipsnames]{xcolor}

\RequirePackage[colorlinks]{hyperref}
\RequirePackage{hypernat}

\usepackage[psamsfonts]{amsfonts}
 
\usepackage{pdflscape}

\usepackage{dsfont}

\newtheorem{Lem}{Lemma}[section]

\newtheorem{Theor}{Theorem}[section]

\newtheorem{defin}{Definition}[section]

\usepackage{hyperref}
\hypersetup{pdfpagemode=FullScreen,  
colorlinks=true,
citecolor=Bittersweet,
linkcolor=Bittersweet,
urlcolor=Bittersweet
}

\newcommand{\cqfd}{\hfill $\square$}

\newcommand{\R}{\mathbb R}

\newcommand{\n}{^{(n)}}

\newcommand{\Xb}{\mathbf{X}}
\newcommand{\Sb}{\mathbf{S}}

\newcommand{\Zb}{\mathbf{Z}}

\newcommand{\ub}{\ensuremath{\mathbf{u}}}
\newcommand{\vb}{\ensuremath{\mathbf{v}}}
\newcommand{\xb}{\ensuremath{\mathbf{x}}}

\newcommand{\Ab}{\ensuremath{\mathbf{A}}}

\newcommand{\Tb}{\ensuremath{\mathbf{T}}}

\newcommand{\Yb}{\ensuremath{\mathbf{Y}}}

\newcommand{\thetab}{{\pmb \theta}}

\newcommand{\Sigb}{{\pmb \Sigma}}

\newcommand{\Deltab}{{\pmb \Delta}}
\newcommand{\taub}{{\pmb \tau}}

\newcommand{\Gammab}{{\pmb \Gamma}}

\newcommand{\pr}{^{\prime}}

\newcommand{\ny}{n\rightarrow\infty}

\begin{document}

\begin{frontmatter}
\title{Inference on the mode of weak directional signals\,:
\\
\mbox{\hspace{-15mm} A Le Cam perspective on hypothesis testing near singularities}}
\runtitle{Inference on the mode of weak directional signals}
\vspace{1mm}

\begin{aug}
\author{\fnms{Davy} \snm{Paindaveine}\thanksref{t1}\ead[label=e1]{dpaindav@ulb.ac.be}
\ead[label=u1,url]{http://homepages.ulb.ac.be/\textasciitilde dpaindav}}
\and
\author{\fnms{ Thomas} \snm{Verdebout}
\ead[label=e2]{tverdebo@ulb.ac.be}
\ead[label=u2,url]{http://tverdebo.ulb.ac.be}}

\thankstext{t1}{Research is supported by\vspace{-1mm} the IAP research network grant 
\vspace{.8mm}
\mbox{nr.} P7/06 of the Belgian government (Belgian Science Policy), the Cr\'{e}dit de Recherche  J.0113.16 of the FNRS (Fonds National pour la Recherche Scientifique), Communaut\'{e} 
\vspace{1mm}
Fran\c{c}aise de Belgique, and a grant from the National Bank of Belgium.}
\runauthor{D. Paindaveine and Th. Verdebout}
\vspace{2mm}

\affiliation{Universit\'{e} libre de Bruxelles}

\address{Universit\'{e} libre de Bruxelles\\
ECARES and D\'{e}partement de Math\'{e}matique\\
Avenue F.D. Roosevelt, 50, ECARES, CP114/04\\
B-1050, Brussels\\
Belgium\\
\printead{e1}\\
\printead{u1}}

\address{Universit\'{e} libre de Bruxelles\\
ECARES and D\'{e}partement de Math\'{e}matique\\
Bld du Triomphe, 
Campus Plaine, CP210\\
B-1050, Brussels\\
Belgium\\
E-mail:\ \printead*{e2}\\
\printead{u2}}

\end{aug}

\vspace{4mm}
\begin{abstract}
We revisit, in an original and challenging perspective, the problem of testing the null hypothesis that the mode of a directional signal is equal to a given value. Motivated by a real data example where the signal is weak, we consider this problem under asymptotic scenarios for which the signal strength goes to zero at an arbitrary rate~$\eta_n$. Both under the null and the alternative, we focus on rotationally symmetric distributions. We show that, while they are asymptotically equivalent under fixed signal strength, the classical Wald and Watson tests exhibit very different (null and non-null) behaviours when the signal becomes arbitrarily weak. To fully characterize how challenging the problem is as a function of~$\eta_n$, we adopt a Le Cam, convergence-of-statistical-experiments, point of view and show that the resulting limiting experiments crucially depend on~$\eta_n$. In the light of these results, the Watson test is shown to be \emph{adaptively} rate-consistent and essentially adaptively Le Cam optimal. Throughout, our theoretical findings are illustrated via Monte-Carlo simulations. The practical relevance of our results is also shown on the real data example that motivated the present work.
\vspace{-1mm}
\end{abstract}

\begin{keyword}[class=MSC]
\kwd[Primary ]{62G10}
\kwd{62G20}
\kwd[; secondary ]{62G35}
\kwd{62H11}
\end{keyword}

\begin{keyword}
\kwd{Contiguity}
\kwd{convergence of statistical experiments}
\kwd{directional statistics}
\kwd{robust tests}
\kwd{rotationally symmetric distributions}
\end{keyword}

\end{frontmatter}


\section{Introduction}
\label{introsec}

In applications involving multivariate data, it is not uncommon that practitioners observe directions only, rather than both directions and magnitudes. Such data are said to be \emph{directional} and are viewed as realizations of a random vector~$\Xb$ with a distribution~${\rm P}$ that only charges the unit sphere ${\cal S}^{p-1}:= \{ {\bf x} \in \R^p : \|\xb\|^2={\bf x}\pr {\bf x}=1\}$ of $\R^p$. Common examples include data related to wind, earth magnetic field or cosmology. In most applications, the primary focus is on location functionals, such as the \emph{spherical mean} ${\rm E}_{\rm P}[\Xb]/\|{\rm E}_{\rm P}[\Xb]\|$ or the mode of~${\rm P}$ (that is, the maximizer of the density of~${\rm P}$ with respect to an appropriate dominating measure on the unit sphere). In this introduction, we focus without loss of generality on the spherical mean,~$\thetab(=\thetab({\rm P}))$ say, since in the rest of the paper, distributional assumptions will ensure that the mode and the spherical mean do coincide. 

Inference on~$\thetab$ has been considered in many papers. Asymptotic score tests for the null hypothesis~${\cal H}_0: \thetab= \thetab_0$  have been studied in \citet[p.\!~140]{Wat1983b} and \cite{PaiVer2015a}, while Wald tests were considered in \cite{HayPur1985}, \cite{Hay1990} and \cite{LaJu03}. Robust M-estimation of $\thetab$ has been tackled in \cite{CR01} and rank-based procedures were proposed in \cite{TS07}, \cite{LSTV13} and \cite{PaiVer2015a}. The score test  for ${\cal H}_0: \thetab= \thetab_0$ has recently been shown to be robust to high-dimensionality in \cite{Leyetal2015}.

Clearly, performing inference on~$\thetab$ is a semiparametric problem whose difficulty depends on the underlying distribution~${\rm P}$~: if~${\rm P}$ is much concentrated about~$\thetab$, then it is in principle easy to, e.g., identify small confidence zones for~$\thetab$. On the contrary, if~${\rm P}$ is close to the uniform distribution~${\rm P}_0$ over the unit sphere, then performing inference on~$\thetab$ is much more delicate and the corresponding confidence zones will be very broad. In line with this, the Fisher information for~$\thetab$ obtained in Proposition~2.2 of \cite{LSTV13} (in the context of rotationally symmetric distributions) explicitly depends on~${\rm P}$ and goes to the zero matrix as~${\rm P}$ converges weakly to~${\rm P}_0$. This singularity, of course, is intimately related to the fact that~$\thetab$ is not identifiable at~${\rm P}_0$. More generally, performing inference on~$\thetab$ is expected to be non-standard and difficult when~$\lambda=\lambda({\rm P})=\|{\rm E}_{\rm P}[\Xb]\|$ is close to the zero value that makes~$\thetab={\rm E}_{\rm P}[\Xb]/\lambda({\rm P})$ undefined.

So far, asymptotic inference on~$\thetab$ has been conducted under the assumption that observations are randomly sampled from a distribution~${\rm P}$ that does not depend on the sample size~$n$. If, however, the directional signal is weak, meaning that~${\rm P}$ is close to~${\rm P}_0$ (or that the corresponding~$\lambda$ value is close to zero), then such a standard asymptotic scenario may be inappropriate for conducting inference on the signal direction~$\thetab$, in the sense that the resulting asymptotic distribution of some statistic of interest may fail, even if~$n$ is large, to provide satisfactory approximations of the corresponding fixed-$n$ distribution. One of the goals of this paper is to show that this may indeed be the case and that it may have dramatic implications on standard inference procedures. Such considerations are relevant as soon as the directional signal is weak, as it is the case for instance for the cosmic ray data set we will consider in Section~\ref{realsec}.
 

As a reaction, we consider in this paper asymptotic scenarios associated with triangular arrays of observations where, for each positive integer~$n$, $\Xb_{n1},\ldots,\Xb_{nn}$ are randomly sampled from a distribution~${\rm P}_n$ over the unit sphere. We will allow the strength of the signal, $\lambda_n=\lambda({\rm P}_n)$ say, to go to zero at an arbitrary rate~$\eta_n$. In the semiparametric model we will actually adopt (whose validity could be tested a priori in the spirit of \citealp{Dette013} or \citealp{Gon14}),  
this is equivalent to allowing the underlying distribution~${\rm P}_n$ to converge to the uniform distribution~${\rm P}_0$ at an arbitrary rate. In this context, we will mainly focus on the problem of testing~${\cal H}_0\n: \thetab_n= \thetab_0$ against~${\cal H}_1\n: \thetab_n\neq \thetab_0$, where~$\thetab_n=\thetab({\rm P}_n)$ is the parameter of interest and~$\thetab_0(\in\mathcal{S}^{p-1})$ is fixed. We first consider the two most famous tests for this problem, namely the score test of \citet[p.\!~140]{Wat1983b} (see also \citealp{Leyetal2015} and \citealp{PaiVer2015a}) and the traditional Wald test based on the sample spherical mean (see \citealp{HayPur1985}, \citealp{Hay1990} and \citealp{LaJu03}). We show that these tests exhibit very different asymptotic null behaviours in the vicinity of uniformity~: while the null behaviour of the Watson test (see~(\ref{Watsontest}) below) is robust to the possible convergence of~${\rm P}_n$ to~${\rm P}_0$, the null behaviour of the Wald test (see~(\ref{Waldtest}) below) is not and crucially depends on the rate~$\eta_n$. This is in sharp contrast with what happens \emph{away from uniformity}, that is for~${\rm P}_n\equiv {\rm P}$, where the Wald and Watson tests have been shown to be asymptotically equivalent under the null; see~\cite{Hay1990} in a specific parametric setup, or Theorem~\ref{LDnullasymp}(i) below in the broader semiparametric framework considered in the present paper. In view of this asymptotic equivalence, practitioners might be tempted to use indifferently the Wald or Watson tests in the vicinity of uniformity as well. However, our results show that, for data sets such as the cosmological one considered in Section~\ref{realsec}, this might have dramatic consequences for inference. 

Of course, robustness of the null behaviour in the vicinity of uniformity should not be obtained at the expense of efficiency. To investigate whether this is the case or not, we also study, as the signal strength goes to zero, the asymptotic distribution of the Watson test under appropriate local alternatives. We show that the weaker the signal (more precisely, the faster the rate~$\eta_n$ at which the signal strength goes to zero), the less severe the alternatives that can be detected by the Watson test (more precisely, the poorer its consistency rate), which is of course reasonable. 
Moreover, if the rate at which the signal vanishes exceeds some threshold, then the Watson test, like the Wald test, is blind to all alternatives, as severe as they may be. We show that this threshold rate, that is, the fastest rate~$\eta_n$ for which some alternatives can be detected by the Watson test, is the slowest rate for which the corresponding distributions~${\rm P}_n$ form a sequence of probability measures that is contiguous to the sequence associated with~${\rm P}_0$. Contiguity will therefore play an important role when quantifying what we call ``vicinity of uniformity". 

%
%
%
%
%
%
%
%
%

Finally, while it is of course nice to identify the alternatives that can be detected by the Watson test for any possible rate~$\eta_n$,  some important questions remain~: (i) for a given rate~$\eta_n$, does there exist a test that can see less severe alternatives than those detected by the Watson test? (ii) If not, does the Watson test maximize the asymptotic power against the least severe alternatives it can detect? To answer these questions, we adopt Le Cam's convergence-of-statistical-experiments approach and derive, for any given rate~$\eta_n$, the corresponding limiting experiments. Interestingly, these limiting experiments are locally asymptotically normal for any~$\eta_n$ yet depend crucially on~$\eta_n$. Our results reveal that (i) the Watson test is \emph{rate-adaptive}, in the sense that, irrespective of~$\eta_n$, no tests can show non-trivial asymptotic powers against less severe alternatives than those detected by the Watson test. We also show that (ii) the Watson test is essentially \emph{adaptively Le Cam optimal}: it is \emph{uniformly (in the underlying distribution) optimal}  whenever the underlying sequences of distributions~${\rm P}_n$ is not contiguous to~${\rm P}_0$ and \emph{uniformly (in the underlying distribution) locally optimal} under contiguity (see Section \ref{LANfixedpsec} for details).

The problem we consider is characterized by the fact that the parameter of interest~$\thetab$ becomes unidentified/undefined when a nuisance parameter takes some given value (here, e.g., when $\lambda=0$). Such situations have been considered in the literature in various frameworks and it has been recognized that performing inference on~$\thetab$ when the nuisance is close to this particular value is challenging. This is particularly true in the field of econometrics; we refer to, e.g., \cite{Duf1997}, \cite{Pot2002}, \cite{ForHil2003}, \cite{Duf2006}, or \cite{For2009}. To the best of our knowledge, the results of this paper are the first to discuss asymptotic optimality issues (through fine Le Cam-type results) in such close-to-singular cases. 
Incidentally, another setup that is of a similar nature is the one associated with Gaussian mixtures of the form
$(1- \lambda_n) \, {\cal N}(0,1)+ \lambda_n \, {\cal N}(\theta_n,1)$.
Many works considered the problem of testing~$\mathcal{H}_0\n:\lambda_n=0$ against alternatives under which~$\lambda_n$ goes to zero and~$\theta_n$ diverges to infinity in an appropriate way; see~\cite{Cai07} and the references therein. If the null is rejected, then it becomes of interest to identify the signal, that is, to perform inference on~$\theta_n$, which is close to being unidentified in the setup considered where~$\lambda_n$ is close to zero.  Our investigation brings precise results in a framework that is very similar to those considered in these econometric and Gaussian-mixtures contexts.

The outline of the paper is as follows. In Section~\ref{rotsymsec}, we introduce the semiparametric model we will focus on and define the sequences of hypotheses converging to the uniform on the unit sphere we will consider. In Section~\ref{Nullfixedpsec}, we recall the Wald and Watson tests and study their asymptotic null behaviour in the vicinity of uniformity. We derive the corresponding local asymptotic powers in Section~\ref{Alternativefixedpsec}. In Section~\ref{LANfixedpsec}, we show that, irrespective of the rate~$\eta_n$ at which convergence to the uniform takes place, the resulting sequences of statistical experiments converge to some limiting experiments (that depend on~$\eta_n$). There, we also exploit these results to make precise what are the (Le Cam) optimality properties of the Watson test (the lack of robustness of the Wald test, which will follow from the results of Sections~\ref{Nullfixedpsec}-\ref{Alternativefixedpsec}, justifies that we restrict to the Watson test when discussing optimality issues). Throughout, our theoretical findings are confirmed by simulation exercises. In Section~\ref{realsec}, we show the practical relevance of our results on a cosmic ray data set. Finally, Section~\ref{conclusec} summarizes the results and an appendix collects technical proofs.

\section{Rotational symmetry and shrinking neighbourhoods of uniformity} 
\label{rotsymsec}

As announced in the introduction, we will  restrict to a specific, semiparametric, class of distributions on the unit sphere~$\mathcal{S}^{p-1}$. More precisely, we will consider absolutely continuous distributions over~$\mathcal{S}^{p-1}$ (with respect to the surface area measure) that admit densities of the form
\begin{equation}
\label{rotpdf}
\xb \mapsto c_{p,\kappa,f} f(\kappa \xb'\thetab) 
,
\end{equation}
where~$\thetab\in\mathcal{S}^{p-1}$, $\kappa>0$ and~$f$ belongs to the collection~$\mathcal{F}$ of functions from~$\R$ to~$\R^+$ that are monotone increasing, twice differentiable at~$0$, and satisfy~$f(0)=f'(0)=1$. Throughout, the distribution with density~(\ref{rotpdf}) will be said to be \emph{rotationally symmetric about~$\thetab$} and will be denoted as~$\mathcal{R}(\thetab,\kappa,f)$. The restrictions on~$\kappa$ and~$f$ above, under which~$\thetab$ is both the unique mode and the spherical mean of the distribution, ensure identifiability of~$\thetab$, $\kappa$ and~$f$. Clearly,~$\kappa$ measures the strength of the directional signal or its ``concentration" (the larger~$\kappa$, the more concentrated the probability mass is about~$\thetab$). If~$\Xb$ has density~(\ref{rotpdf}), then $\Xb'\thetab$ has density~$c_{p,\kappa,f} (1-t^2)^{(p-3)/2}f(t)$ over~$[-1,1]$ (see, e.g. \citealp[p.\!~136]{Wat1983b}), which shows that the normalization constant in~(\ref{rotpdf}) is given by~$c_{p,\kappa,f}=1/\int_{-1}^1 (1-t^2)^{(p-3)/2}f(t)\,dt$. It is important to note that, irrespective of~$p$ and~$f$, the boundary case~$\kappa=0$ corresponds to the uniform distribution over~$\mathcal{S}^{p-1}$. 
The celebrated Fisher--von Mises--Langevin (FvML) distributions correspond to the particular case~$t\mapsto f(t)=\exp(t)$.

As explained in the introduction, our main focus will be on sequences of hypotheses that are in the vicinity of the uniform distribution. In the present setup, the corresponding ``shrinking neighbourhoods" of uniformity require considering triangular arrays of observations of the form 
$$
\Xb_{ni}, 
\quad i=1,\ldots,n, 
\quad n=1,2,\ldots
,
$$ 
where, for any~$n$, $\Xb_{n1},\ldots,\Xb_{nn}$ form a random sample from~$\mathcal{R}(\thetab_n,\kappa_n,f)$; 
the resulting sequence of hypotheses will be denoted as~${\rm P}^{(n)}_{\thetab_n,\kappa_n,f}$.
Here, $(\thetab_n)$ is a sequence in~$\mathcal{S}^{p-1}$, $(\kappa_n)$ is a sequence in~$\R^+_0$,
and~$f\in\mathcal{F}$ is fixed. The sequence of hypotheses under which, for any~$n$, $\Xb_{n1},\ldots,\Xb_{nn}$ form a random sample from the uniform over~$\mathcal{S}^{p-1}$ will be denoted as~${\rm P}^{(n)}_{0}$ (for convenience, we will also put~${\rm P}^{(n)}_{\thetab,0,f}:={\rm P}^{(n)}_{0}$ for any~$\thetab,f$). Since~$\kappa=0$ corresponds to the uniform distribution over~$\mathcal{S}^{p-1}$, it is natural to adopt the following definition, that  allows~${\rm P}^{(n)}_{\thetab_n,\kappa_n,f}$ to converge to~${\rm P}^{(n)}_{0}$ at an arbitrary rate.

\begin{defin}
\label{definvicf}
Fix a sequence~$(\thetab_n)$ in~$\mathcal{S}^{p-1}$, $f\in\mathcal{F}$, $\xi>0$ and a sequence~$(\eta_n)$ in~$\R^+$ that is~$o(1)$ as~$\to\infty$. Then the sequence of hypotheses~${\rm P}\n_{\thetab_n,\kappa_n,f}$ is in an \emph{$\eta_n$-neighbourhood of uniformity, with locality parameter~$\xi$}, if and only if~$\kappa_n=\sqrt{p}\eta_n\xi+o(\eta_n)$ as~$n\to\infty$.  
\end{defin}




The presence of~$\sqrt{p}$ in the expression~$\kappa_n=\sqrt{p}\eta_n\xi+o(\eta_n)$ may be unexpected at first and will be explained below Definition~\ref{definvicF}. To widen the scope of our results as much as possible, we will often consider more general rotationally symmetric distributions. We will say that the random $p$-vector~$\Xb$, with values on~$\mathcal{S}^{p-1}$, is rotationally symmetric about location~$\thetab(\in\mathcal{S}^{p-1})$ if $\mathbf{O}\Xb$ is equal in distribution to~$\Xb$ for any orthogonal $p\times p$ matrix~$\mathbf{O}$ satisfying $\mathbf{O}\thetab=\thetab$. Such general rotationally symmetric distributions, that do not need be absolutely continuous nor have a concentration that is governed by a parameter~$\kappa$, are characterized by the location parameter~$\thetab$ and the cumulative distribution function~$F$ of~$\Xb'\thetab$. The corresponding distribution will be denoted by~$\mathcal{R}(\thetab,F)$.
Parallel as above,~${\rm P}^{(n)}_{\thetab_n,F_n}$ will then refer to triangular arrays of observations for which~$\Xb_{n1},\ldots,\Xb_{nn}$ form a random sample from~$\mathcal{R}(\thetab_n,F_n)$, where $(\thetab_n)$ is still a sequence in~$\mathcal{S}^{p-1}$ and where~$(F_n)$ is a sequence of cumulative distribution functions on~$[-1,1]$. 

Of course, it is desirable to identify conditions that make sequences of hypotheses~${\rm P}^{(n)}_{\thetab_n,F_n}$ be in $\eta_n$-neighbourhoods of uniformity. It actually follows from~(5.2)-(5.3) in \cite{Cut2015} that, under~${\rm P}\n_{\thetab_n,\kappa_n,f}$, with~$\kappa_n=\sqrt{p}\eta_n\xi+o(\eta_n)$ (where~$\eta_n=o(1)$), one has
\begin{equation}
\label{tatatar}
e_{n1}
:=
{\rm E}[\Xb_{n1}'\thetab_n]
= 
\frac{\eta_n\xi}{\sqrt{p}}  + 
o(\eta_n)
\quad
\textrm{and}
\quad 
\tilde{e}_{n2}
:=
{\rm Var}[\Xb_{n1}'\thetab_n]
=
\frac{1}{p}   + 
o(1)
\end{equation}
as~$n\to\infty$, which is to be compared with the values~$e_{n1}=0$ and~$\tilde{e}_{n2}=1/p$ obtained under~${\rm P}^{(n)}_{0}$. This motivates the following definition.

\begin{defin}
\label{definvicF}
Fix a sequence~$(\thetab_n)$ in~$\mathcal{S}^{p-1}$, a sequence $(F_n)$ of cumulative distribution functions on~$[-1,1]$, $\xi>0$,  and a sequence~$(\eta_n)$ in~$\R^+$ that is~$o(1)$ as~$\to\infty$. Then the sequence of hypotheses~${\rm P}\n_{\thetab_n,F_n}$ is in an \emph{$\eta_n$-neighbourhood of uniformity, with locality parameter~$\xi$}, if and only if
$$
e_{n1}
=
\frac{\eta_n\xi}{\sqrt{p}}  + 
o(\eta_n)
\quad
\textrm{and}
\quad
\tilde{e}_{n2}
=
\frac{1}{p}
+
o(1)
$$
as~$\ny$, where $e_{n1}={\rm E}[\Xb_{n1}'\thetab_n]$
and $\tilde{e}_{n2}={\rm Var}[\Xb_{n1}'\thetab_n]$ are evaluated under~${\rm P}\n_{\thetab_n,F_n}$.  
\end{defin}

 In the present ``low-dimensional" (fixed-$p$) setup, it might have been more natural to define $\eta_n$-neighbourhoods of uniformity with locality parameter~$\xi$ through~$\kappa_n=\eta_n\xi+o(\eta_n)$ in Definition~\ref{definvicf} (which would then translate into~$e_{n1}=\kappa_n=(\eta_n\xi)/p+o(\eta_n)$ in Definition~\ref{definvicF}). Of course, appropriate reparametrization of~$\xi$ into~$p^{\pm 1/2}\xi$ makes the definitions we adopted above and these possible alternative ones perfectly equivalent. The reason why we favour Definitions~\ref{definvicf}-\ref{definvicF} is that they would make easier possible future comparisons between the low- and high-dimensional cases.  
 
Throughout, it will be of interest to compare the results obtained in the vicinity of uniformity to the standard ones obtained away from uniformity. In the framework of Definition~\ref{definvicf}, we will say that the sequence of hypotheses~${\rm P}\n_{\thetab_n,\kappa_n,f}$ stays away from uniformity if and only if~$\kappa_n\to\kappa(>0)$ as~$n\to\infty$. This corresponds to sequences of hypotheses~${\rm P}\n_{\thetab_n,F_n}$ for which $e_{n1}$ and~$\tilde{e}_{n2}$ converge to positive constants (say $e_1$ and~$\tilde{e}_2$, respectively), with the important difference that~$\tilde{e}_2$ here does not need be equal to~$1/p$. For instance, in the FvML case with concentration~$\kappa_n$ converging to~$\kappa(>0)$, one has 
\begin{equation}
\label{e2fvml}
e_1
=
\frac{\mathcal{I}_{p/2}(\kappa)}{\mathcal{I}_{p/2-1}(\kappa)}
\quad\textrm{and}\quad
\tilde{e}_2
=
- 
\frac{p-1}{\kappa}
\frac{\mathcal{I}_{p/2}(\kappa)}{\mathcal{I}_{p/2-1}(\kappa)}
+1
-
\bigg(\frac{\mathcal{I}_{p/2}(\kappa)}{\mathcal{I}_{p/2-1}(\kappa)}\bigg)^2
,
\end{equation}  
where~$\mathcal{I}_r(\cdot)$ denotes the order-$r$ modified Bessel function of the first kind; see, for instance, Lemma~S.2.1 in \cite{Cut2015s}. To present the results in a setup that is closely related to the one we adopted above for neighbourhoods of uniformity, we then have the following definition (that should be compared to Definition~\ref{definvicF}).

\begin{defin}
\label{definvicFfixed}
Fix a sequence~$(\thetab_n)$ in~$\mathcal{S}^{p-1}$, a sequence $(F_n)$ of cumulative distribution functions over~$[-1,1]$, and $\xi,\tilde{e}_{2}>0$. Then the sequence of hypotheses~${\rm P}\n_{\thetab_n,F_n}$  stays away from uniformity (or is in a $1$-neighbourhood of uniformity), with locality parameters~$\xi$ and~$\tilde{e}_{2}$, if and only if
$$
e_{n1}
=
\frac{\xi}{\sqrt{p}}  + 
o(1)
\quad
\textrm{and}
\quad
\tilde{e}_{n2}
=
\tilde{e}_2
+
o(1)
$$
as~$\ny$, where $e_{n1}={\rm E}[\Xb_{n1}'\thetab_n]$
and $\tilde{e}_{n2}={\rm Var}[\Xb_{n1}'\thetab_n]$ are evaluated under~${\rm P}\n_{\thetab_n,F_n}$.  
\end{defin}

As discussed in the introduction, the closer to uniformity the distribution is, the more challenging it should be to perform inference about~$\thetab$. While we will mostly focus on hypothesis testing in the sequel, we present here the following point estimation result, that describes how the most natural estimator for~$\thetab$, namely the (sample) spherical mean, which is the MLE for~$\thetab$ in the FvML parametric submodel, deteriorates when the underlying distribution gets closer to uniformity (throughout, $\stackrel{\mathcal{D}}{\to}$ denotes weak convergence).

\begin{Theor} 
\label{asymsphermean} 
Let~$(\eta_n)$ be either the sequence~$\eta_n\equiv 1$ or a sequence in~$\R^+$ that is~$o(1)$. Assume that~${\rm P}\n_{\thetab,F_n}$ is in an $\eta_n$-neighbourhood of uniformity, with locality parameters~$\xi$ and~$\tilde{e}_2$ if $\eta_n\equiv 1$ and with locality parameter~$\xi$ otherwise. Let 
$\hat{\thetab}_n=\bar{\Xb}_n/\| \bar{\Xb}_n\|$, with~$\bar{\Xb}_n:=\frac{1}{n} \sum_{i=1}^n \Xb_{ni}$.
Then we have the following as~$n\to\infty$ under~${\rm P}\n_{\thetab,F_n}$~:\\
(i) if~$\eta_n\equiv 1$, then 
$$
\sqrt{n}(\hat{\thetab}_n- \thetab)
 \stackrel{\mathcal{D}}{\to}
{\cal N} 
\bigg(
{\bf 0}
\,
, 
\,
{\frac{1-\xi^2/p-\tilde{e}_2}{\xi^2(1-1/p)}}
\,
 ({\bf I}_p- \thetab\thetab\pr)
\bigg) 
;
$$
%
(ii) if~$\eta_n=o(1)$ with~$\sqrt{n}\eta_n\to\infty$, then 
$$
\sqrt{n}\eta_n(\hat{\thetab}_n- \thetab)
 \stackrel{\mathcal{D}}{\to}
{\cal N} 
\bigg(
{\bf 0}
\,
, 
{\frac{1}{\xi^2}} 
\,
({\bf I}_p- \thetab\thetab\pr)
\bigg) 
;
$$
(iii) if~$\sqrt{n}\eta_n\to 1$, then 
$$
\hat{\thetab}_n
 \stackrel{\mathcal{D}}{\to}
\frac{\Zb}{\|\Zb\|},
\quad
\textrm{with }
\Zb
\sim
{\cal N} \big(
\xi \thetab
\,
, 
{\bf I}_p
\big) 
;
$$
(iv) if~$\sqrt{n}\eta_n\to 0$, then 
$
\hat{\thetab}_n
 \stackrel{\mathcal{D}}{\to}
{\rm Unif}(\mathcal{S}^{p-1}),
$
the uniform distribution over~$\mathcal{S}^{p-1}$. 
\end{Theor}

Part~(i) of the result states that standard root-$n$ consistency is obtained away from uniformity. The faster the underlying distribution converges to the uniform in~(ii), the poorer the resulting consistency rate of~$\hat\thetab_n$, that may become arbitrarily slow. In case~(iii),~$\hat\thetab_n$ fails to be consistent, but its asymptotic distribution still depends on the true value of~$\thetab$. Finally, in case~(iv), we are so close to the uniform case that~$\hat\thetab_n$ behaves like uniform noise on the sphere, hence does not bear any information on~$\thetab$. This result therefore confirms that the performance of~$\hat{\thetab}_n$ deteriorates (monotonically) as the speed at  which the underlying distribution converges to the uniform increases. 

Theorem~\ref{asymsphermean} also hints that the rate~$\eta_n\sim 1/\sqrt{n}$ will play a special role in this paper. This rate is actually the slowest one for which~${\rm P}\n_{0}$ and~${\rm P}\n_{\thetab_n,\kappa_n,f}$, with~$\kappa_n=\sqrt{p}\eta_n\xi+o(\eta_n)$, are mutually contiguous (this is a corollary of Theorem~2.2 in \cite{Cut2015}, that states that, for any fixed~$f\in\mathcal{F}$, the sequence of (concentration) parametric models~$\big\{ {\rm P}\n_{\thetab_n,\kappa,f}:\kappa\geq 0 \big\}$ is \emph{locally and asymptotically normal (LAN)} at~$\kappa=0$, with contiguity rate~$1/\sqrt{n}$). In line with Theorem~\ref{asymsphermean}, most results in the sequel will discriminate between the following regimes~: \emph{away from uniformity} ($\eta_n\equiv 1$), \emph{beyond contiguity} ($\eta_n=o(1)$ with~$\sqrt{n}\eta_n\to\infty$), \emph{under contiguity} ($\eta_n\sim 1/\sqrt{n}$) and \emph{under strict contiguity} ($\sqrt{n}\eta_n\to 0$).


\section{Contiguity-robust testing} 
\label{Nullfixedpsec}

In the most general rotationally symmetric setup introduced in the previous section, we consider the problem of testing the null hypothesis that the modal location~$\thetab$ is equal to some given location~$\thetab_0$, under unspecified cumulative distribution function~$F$. More precisely, using the notation introduced in Section~\ref{rotsymsec}, we consider the testing problem
\begin{equation}
\label{fixedpproblem}
\mathcal{H}_0^{(n)}: 
\cup_{F}
\,
\big\{ {\rm P}^{(n)}_{\thetab_{0},F} \big\} 
\quad
\textrm{ against }
\quad
\mathcal{H}_1^{(n)}: 
\cup_{\thetab\neq \thetab_{0}} \cup_{F}
\big\{ {\rm P}^{(n)}_{\thetab,F} \big\} 
,
\end{equation}
where $\thetab_0(\in\mathcal{S}^{p-1})$ is fixed and the unions in~$F$ are over the collection of cumulative distribution functions on~$[-1,1]$. In this section, we investigate whether or not the two most classical tests for this problem remain valid (in the sense that they still meet asymptotically the nominal level constraint) in the vicinity of uniformity. 

These classical tests are based on the sample average~$\bar{\Xb}_n:=\frac{1}{n}\sum_{i=1}^n \Xb_{ni}$ 
of the observations~$\Xb_{ni}$, $i=1,\ldots,n$ at hand, and take the following form~:
\begin{itemize}
\item[(i)] the Watson score test~$\phi_n^W$ (\citealp[p.\!~140]{Wat1983b}) rejects the null~$\mathcal{H}_0^{(n)}$ at asymptotic level~$\alpha$ whenever
\begin{equation}
\label{Watsontest}
W_n 
:=
\frac{n (p-1)\bar{\Xb}_n\pr ({\bf I}_p- \thetab_{0} \thetab_{0}\pr) \bar{\Xb}_n}{1-\frac{1}{n} \sum_{i=1}^n(\Xb_{ni}\pr\thetab_{0})^2}
>
\chi^2_{p-1,1-\alpha}
,
\end{equation}
where~$\mathbf{I}_p$ denotes the $p$-dimensional identity matrix  and~$\chi^2_{\ell,1-\alpha}$ stands for the $\alpha$-upper quantile of the $\chi^2_\ell$ distribution. 
\vspace{2mm}
\item[(ii)]  
The Wald test~$\phi_n^S$ (\citealp{Hay1990,HayPur1985}) rejects the null~$\mathcal{H}_0^{(n)}$ at asymptotic level~$\alpha$ if
\begin{equation}
\label{Waldtest}
S_n
:=
\frac{n(p-1) (\bar{\Xb}_{n}\pr\thetab_{0})^2
\,
 \hat{\thetab}_n\pr({\bf I}_p- \thetab_{0}\thetab_{0}\pr)\hat{\thetab}_n}{1-\frac{1}{n} \sum_{i=1}^n(\Xb_{ni}\pr\thetab_{0})^2}
>
\chi^2_{p-1,1-\alpha}
,
\vspace{2mm}
\end{equation}
where $\hat{\thetab}_n=\bar{\Xb}_n/\| \bar{\Xb}_n\|$ is the estimator of~$\thetab_n$ considered in Theorem~\ref{asymsphermean}. 
\end{itemize}

The test statistics~$W_n$ and~$S_n$ are known to be asymptotically equivalent in probability away from uniformity, which is confirmed in Part~(i) of Theorem~\ref{LDnullasymp} below. The main goal of this theorem, however, is to describe the asymptotic null behaviour of these test statistics in the vicinity of uniformity (see the appendix for a proof).

\begin{Theor} 
\label{LDnullasymp} 
Let~$(\eta_n)$ be either the sequence~$\eta_n\equiv 1$ or a sequence in~$\R^+$ that is~$o(1)$. 
Assume that~${\rm P}\n_{\thetab_n,F_n}$ is in an $\eta_n$-neighbourhood of uniformity, with locality parameters~$\xi$ and~$\tilde{e}_2$ if $\eta_n\equiv 1$ and with locality parameter~$\xi$ otherwise, for some~$\xi,\tilde{e}_{2}>0$.
Then we have the following as~$n\to\infty$ under~${\rm P}\n_{\thetab_0,F_n}$~: 
(i) if~$\eta_n\equiv 1$ or (ii) if~$\eta_n=o(1)$ with~$\sqrt{n}\eta_n\to\infty$, then 
$$
W_n\stackrel{\mathcal{D}}{\to}\chi^2_{p-1}
\quad
\textrm{and}
\quad
S_n\stackrel{\mathcal{D}}{\to}\chi^2_{p-1}
$$
$($and one actually then has $S_n=W_n+o_{\rm P}(1))$;
(iii) if~$\sqrt{n}\eta_n\to 1$, then 
$$
W_n\stackrel{\mathcal{D}}{\to}\chi^2_{p-1}
\quad
\textrm{and}
\quad
S_n 
\stackrel{\mathcal{D}}{\to}
\bigg(
1+
\frac{Q}{(Z+\xi)^2}
\bigg)^{-1}
Q
,
\vspace{1mm}
$$
where~$Z\sim\mathcal{N}(0,1)$ and~$Q\sim\chi^2_{p-1}$ are independent;
(iv) if~$\sqrt{n}\eta_n\to 0$, then 
$$
W_n\stackrel{\mathcal{D}}{\to}\chi^2_{p-1}
\quad
\textrm{and}
\quad
S_n\stackrel{\mathcal{D}}{\to}
\bigg(
1+
\frac{Q}{Z^2}
\bigg)^{-1}
Q
,
$$
still where~$Z\sim\mathcal{N}(0,1)$ and~$Q\sim\chi^2_{p-1}$ are independent.
\end{Theor}

This result shows that the asymptotic equivalence in probability between the Watson and Wald test statistics survives beyond contiguity (case~(ii)), but does not under (strict) contiguity. Also, we see that the Watson test~$\phi^W_n$ remains asymptotically valid in the vicinity of uniformity, irrespective of the rate~$\eta_n$ at which the convergence to the uniform takes place. In contrast, the Wald test~$\phi^S_n$ fails to be asymptotically valid under (strict) contiguity, hence is not robust. In the contiguous case (case~(iii)), the asymptotic null distribution of the Wald statistic depends on the locality parameter~$\xi$, which, even in the unrealistic case in which it would be known that the contiguous regime is the ``true" one, would jeopardise implementation of the Wald test. 

\begin{figure}[htbp!]
\begin{center} 
\includegraphics[height=5.2cm, width=132mm]{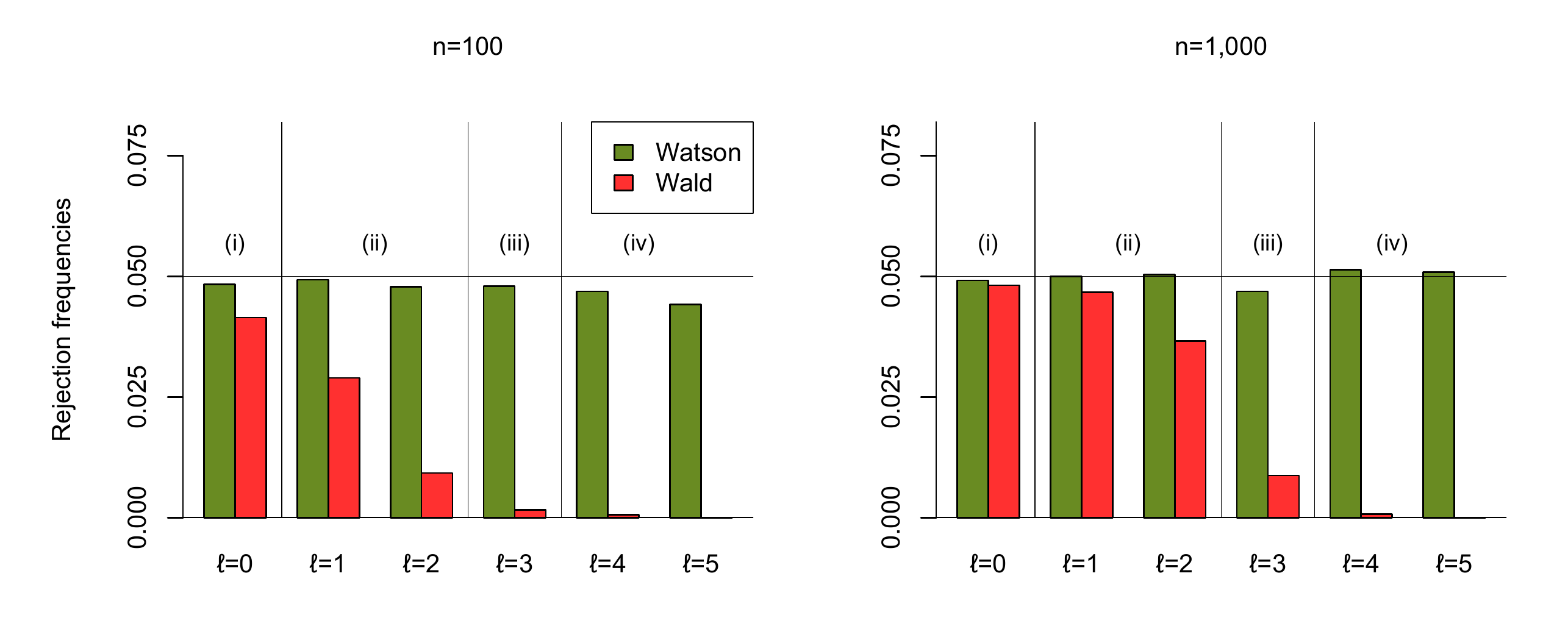}
\vspace{-5mm}
\caption{Rejection frequencies of the Watson (green) and Wald (red) tests for~$\mathcal{H}_0\n:\thetab_n=\thetab_0(=(0,0,1)')$, when performed, at nominal level~$\alpha=5\%$, on $M=10,\!000$ independent random samples of size~$n=100$ (left) and size~$n=1,\!000$ (right) from the FvML distribution on~$\mathcal{S}^2$ with modal location~$\thetab_0$ and a concentration~$\kappa_\ell$ such that~$e^{(\ell)}_{n1}={\rm E}[\Xb^{(\ell)\prime}_{ni}\thetab_0]=n^{-\ell/6}/\sqrt{p}$, for 
(i)~$\ell=0$ (away from uniformity), 
{\protect\linebreak}
(ii) $\ell=1,2$ (beyond contiguity),
(iii) $\ell=3$ (under contiguity), and
(iv) $\ell=4,5$ (under strict contiguity).}
\label{hist1}
\end{center}
\end{figure}

To illustrate Theorem~\ref{LDnullasymp} numerically, we generated, for each value of~$\ell=0,1, \ldots,5$, a collection of~$M=10,\!000$ random samples~$\Xb_{n1}^{(\ell)},\ldots,\Xb_{nn}^{(\ell)}$ from the FvML distribution on~$\mathcal{S}^2$ with modal location~$\thetab_0=(0,0,1)'$ and a concentration~$\kappa_\ell$ that is such that~$e^{(\ell)}_{n1}={\rm E}[\Xb^{(\ell)\prime}_{ni}\thetab_0]=n^{-\ell/6}/\sqrt{p}$; this yields, for~$\ell>0$, $n^{-\ell/6}$-neighbourhoods of uniformity with locality parameter~$\xi=1$, and, for~$\ell=0$, 1-neighbourhoods of uniformity with locality parameters~$\xi=1$ and~$\tilde{e}_2=\tilde{e}_2(\kappa_\ell)$ from~(\ref{e2fvml}). The various values of~$\ell$ clearly allow us to consider all regimes considered in Theorem~\ref{LDnullasymp}~: 
(i) away from uniformity ($\ell=0$),
(ii) beyond contiguity ($\ell=1,2$),
(iii) under contiguity ($\ell=3$), and
(iv) under strict contiguity ($\ell=4,5$).
Figure~\ref{hist1} reports, for sample sizes~$n=100$ and~$n=1,\!000$, the resulting empirical rejection frequencies of the Watson and Wald tests for~$\mathcal{H}_0\n:\thetab_n=\thetab_0(=(0,0,1)')$, performed at nominal level~$5\%$. Clearly, this confirms the robustness of the Watson test and reveals that the Wald test becomes extremely conservative close to uniformity. 

\section{Local asymptotic powers} 
\label{Alternativefixedpsec}

Theorem~\ref{LDnullasymp} above shows that the classical Watson test~$\phi_n^W$ remains valid in the vicinity of uniformity. But of course, it is desirable that this validity-robustness is not obtained at the expense of efficiency, that is, it is desirable that the Watson test still exhibits high local asymptotic powers in the vicinity of uniformity. We now investigate whether this is the case or not. 

Consider a local perturbation~$\thetab_{0}+ \nu_n \taub_n$ of the null value~$\thetab_{0}$, where the sequence~$(\taub_n)$ in~$\R^p$ converges to~$\taub(\neq 0)$, so that the severity, in terms of rate, of such local alternatives is measured by the sequence~$\nu_n$. Of course, it is assumed that~$\thetab_{0}+ \nu_n \taub_n\in{\mathcal S}^{p-1}$ for any~$n$, which imposes that
$
1
=
(\thetab_{0}+ \nu_n \taub_n)\pr (\thetab_{0}+ \nu_n \taub_n) 
=
1 +2 \nu_n \thetab_{0}\pr \taub_n+ \nu_n^2 \| \taub_n \|^2
,
$
or equivalently, that
\begin{equation}
\label{spherconstraint}
\thetab_{0}\pr \taub_n=-\frac{1}{2} \nu_n \| \taub_n \|^2
\
\big(\!=O(\nu_n)\big)
.
\end{equation}
If~$\nu_n=o(1)$, then this leads to~$\thetab_{0}\pr \taub=0$. If~$\nu_n\equiv 1$, then we must rather have~$\thetab_{0}\pr \taub=-\| \taub \|^2/2$. The following result derives the asymptotic distributions of the Watson and Wald test statistics under appropriate alternatives of this form.

\begin{Theor}
\label{LDnonnullasymp} 
Let~$(\eta_n)$ be either the sequence~$\eta_n\equiv 1$ or a sequence in~$\R^+$ that is~$o(1)$. Let~$(\taub_n)$ be a sequence in~$\R^p$ converging to~$\taub(\neq 0)$ and that is such that~$\thetab_{0}+ \nu_n \taub_n\in{\mathcal S}^{p-1}$ for any~$n$, where~$\nu_n=1/(\sqrt{n}\eta_n)$ if~$\sqrt{n}\eta_n\to \infty$ (away from uniformity or beyond contiguity) and~$\nu_n\equiv 1$ if~$\sqrt{n}\eta_n=O(1)$ (under contiguity or under strict contiguity). Assume that~${\rm P}\n_{\thetab_0+\nu_n\taub_n,F_n}$ is in an $\eta_n$-neighbourhood of uniformity, with locality parameters~$\xi$ and~$\tilde{e}_2$ if $\eta_n\equiv 1$ and with locality parameter~$\xi$ otherwise, for some~$\xi,\tilde{e}_2>0$.
Then we have the following as~$n\to\infty$ under~${\rm P}\n_{\thetab_0+\nu_n\taub_n,F_n}$~: 
(i) if~$\eta_n\equiv 1$, then 
$$
W_n\stackrel{\mathcal{D}}{\to}
\chi^2_{p-1}
\bigg(
\,
\frac{1-1/p}{1-\xi^2/p-\tilde{e}_2} \, \xi^2 \|\taub\|^2
\bigg)
\quad
\textrm{and}
\quad
S_n\stackrel{\mathcal{D}}{\to}
\chi^2_{p-1}
\bigg(
\,
\frac{1-1/p}{1-\xi^2/p-\tilde{e}_2} \, \xi^2 \|\taub\|^2
\bigg)
$$
$($and one actually then has $S_n=W_n+o_{\rm P}(1))$;
(ii) if~$\eta_n=o(1)$ with~$\sqrt{n}\eta_n\to\infty$, then 
$$
W_n\stackrel{\mathcal{D}}{\to}
\chi^2_{p-1}
\big(
\xi^2 \|\taub\|^2
\big)
\quad
\textrm{and}
\quad
S_n\stackrel{\mathcal{D}}{\to}
\chi^2_{p-1}
\big(
\xi^2 \|\taub\|^2
\big)
$$
$($and one then still has $S_n=W_n+o_{\rm P}(1))$;
(iii) if~$\sqrt{n}\eta_n\to 1$, then 
$$
W_n\stackrel{\mathcal{D}}{\to}
\chi^2_{p-1}
\Big(
\,\frac{1}{4}\xi^2 \|\taub\|^2 (4-\|\taub\|^2)
\Big)
\quad
\textrm{and}
\quad
S_n
\stackrel{\mathcal{D}}{\to}
\bigg(
1+
\frac{Q_{\xi,\taub}}{(Z+\xi-\xi\|\taub\|^2/2)^2}
\bigg)^{-1}
Q_{\xi,\taub}
,
$$
where~$Z\sim\mathcal{N}(0,1)$ and~$Q_{\xi,\taub}\sim\chi^2_{p-1}\big(
\xi^2 \|\taub\|^2
\big)
$ are independent;
%
(iv) if~$\sqrt{n}\eta_n\to 0$, then 
$$
W_n\stackrel{\mathcal{D}}{\to}\chi^2_{p-1}
\quad
\textrm{and}
\quad
S_n
\stackrel{\mathcal{D}}{\to}
\bigg(
1+
\frac{Q}{Z^2}
\bigg)^{-1}
Q
,
$$
where~$Z\sim\mathcal{N}(0,1)$ and~$Q\sim\chi^2_{p-1}$ are independent.
\end{Theor}

This result shows that the asymptotic equivalence in probability, away from uniformity and beyond contiguity, between the Watson and Wald tests not only holds under the null but also extends to the local alternatives considered. Both tests there exhibit non-trivial asymptotic powers against alternatives that are increasingly severe when the rate at which the underlying distribution converges to the uniform gets faster; note that, in line with Theorem~\ref{asymsphermean}, the consistency rate goes from the standard~$\nu_n=1/\sqrt{n}$ rate away from uniformity to rates that are arbitrarily slow close to contiguity. Under contiguity, the Watson test detects alternatives at a constant rate~$\nu_n\equiv 1$, yet fails to be consistent there, irrespective of the fixed alternative~$\thetab_n\equiv \thetab_0+\taub(\in\mathcal{S}^{p-1})$ considered.
Finally, under strict contiguity, both the Watson and Wald tests are blind to such fixed alternatives, hence cannot show non-trivial asymptotic powers against any alternative there. 

The non-centrality parameter in the asymptotic distribution of~$W_n$ in Theorem~\ref{LDnonnullasymp}(iii) may seem puzzling at first sight, compared to the more standard ones in~(i)-(ii). Note that the Watson test essentially rejects the null for large values of~$\|({\bf I}_p- \thetab_{0} \thetab_{0}\pr) \bar{\Xb}_n\|$, that is, for large values of the norm of the projection of~$\bar{\Xb}_n$ onto the orthogonal complement to~$\thetab_0$. It therefore makes sense that the non-centrality parameter in Theorem~\ref{LDnonnullasymp}(iii) (resp., the corresponding asymptotic power of the Watson test) increases from its minimum value zero (resp., its minimum value~$\alpha$) to its maximum value when~$\|\taub\|$ increases from~$0$ ($\thetab$ equal to the ``north pole"~$\thetab_0$) to~$\sqrt{2}$ ($\thetab$ belongs to the ``equator" with respect to~$\thetab_0$) and decreases from its maximum value to its minimum value zero (resp., its minimum value~$\alpha$) when~$\|\taub\|$ increases from~$\sqrt{2}$  ($\thetab$ belongs to the equator) to~$2$ ($\thetab$ equal to the ``south pole"~$-\thetab_0$). 

We performed the following simulation exercise to see how well the finite-sample behaviours of the Watson and Wald tests actually reflect the theoretical results of Theorem~\ref{LDnonnullasymp}. For each combination of~$\ell=0,1,2,3$ and~$r=0,1,\ldots,6$, we generated~$M=10,\!000$ independent FvML random samples~$\Xb_{n1}^{(\ell,r)},\ldots,\Xb_{nn}^{(\ell,r)}$ of size~$n=200$ on~$\mathcal{S}^2$ with a modal location~$\thetab_r^{(\ell)}$ given in~(\ref{localter}) below and a concentration~$\kappa_\ell$ that is such that~$e^{(\ell)}_{n1}={\rm E}[\Xb^{(\ell)\prime}_{ni}\thetab_r^{(\ell)}]=n^{-\ell/4}/\sqrt{p}$. The integer~$\ell$ allows us to consider the various asymptotic regimes, namely (i) away from uniformity ($\ell=0$), (ii) beyond contiguity ($\ell=1$), under contiguity ($\ell=2$), and under strict contiguity ($\ell=3$). Alternatives were chosen according to the rates in Theorem~\ref{LDnonnullasymp}, and are associated with
\begin{equation}
\label{localter}
\thetab_r^{(\ell)}
:=
\left\{
\begin{array}{ll}
\frac{\thetab_0+ n^{\frac{\ell}{4}-\frac{1}{2}} (\frac{r}{6}\taub_{\rm max})}{\|\thetab_0+ n^{\frac{\ell}{4}-\frac{1}{2}} (\frac{r}{6}\taub_{\rm max})\|}
& \textrm{ for } \ell=0,1
\\[4mm]
\thetab_0+ \taub_r=\vb_r
& \textrm{ for } \ell=2,3,
\end{array} 
\right.
\end{equation}
with~$\thetab_0=(0,0,1)'$, $\taub_{\rm max}=(2,0,0)'$, $r=0,1,\ldots,6$, and~$\vb_r=(\sin(r\pi/6),0,\cos(r\pi/6))'$; clearly, irrespective of~$\ell$, the value~$r=0$ corresponds to the null hypothesis, whereas $r=1,2,\ldots,6$ provide increasingly severe alternatives. 

\begin{figure}[htbp!] 
\begin{center}
\includegraphics[width=152mm]{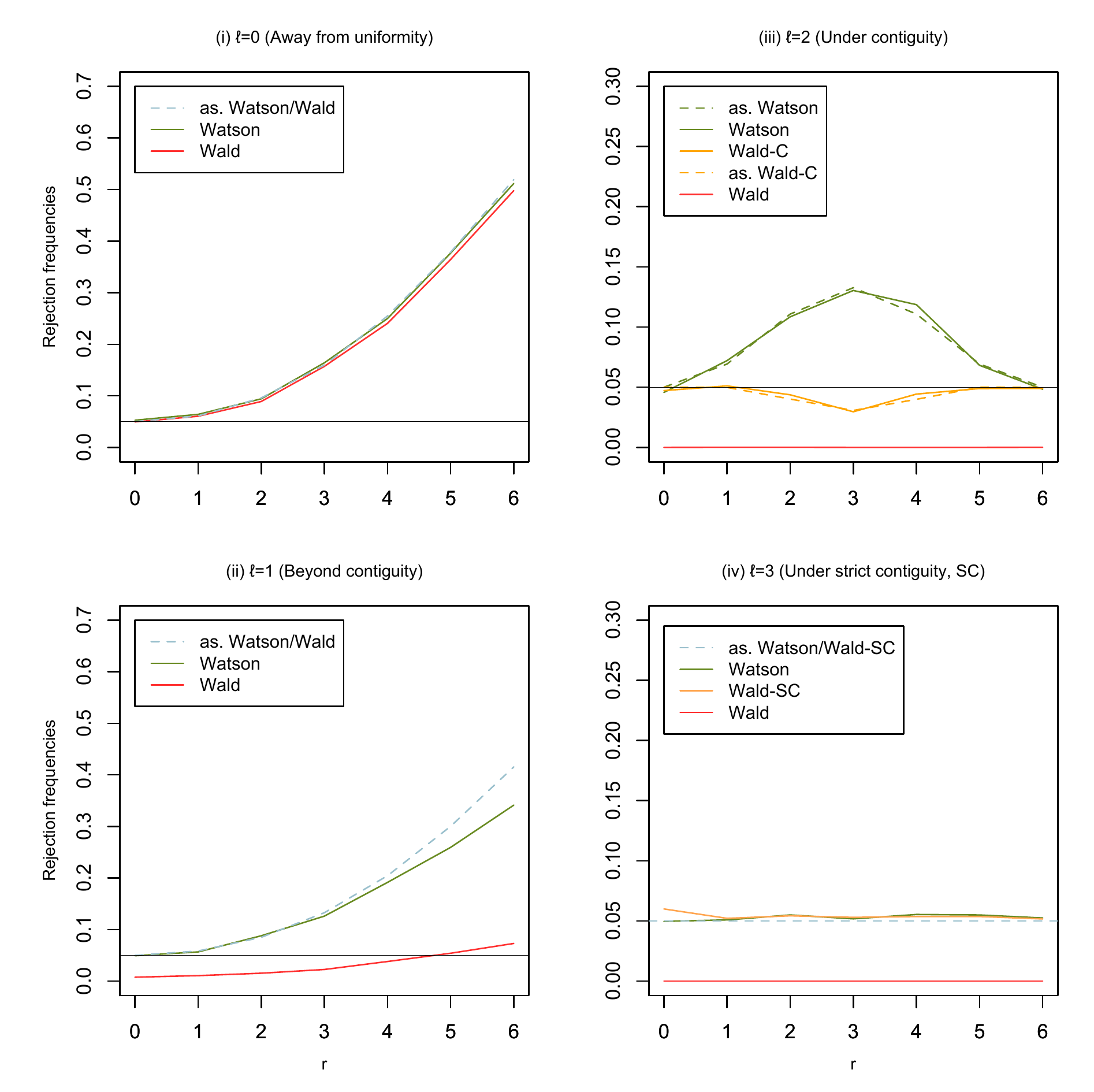}
\caption{
Rejection frequencies of various tests for~$\mathcal{H}_0\n:\thetab_n=\thetab_0(=(0,0,1)')$, when performed, at nominal level~$\alpha=5\%$, on $M=10,\!000$ independent random samples of size~$n=200$ from the FvML distribution on~$\mathcal{S}^2$ with modal location~$\thetab_r^{(\ell)}$ in~(\ref{localter}) and a concentration~$\kappa_\ell$ such that~$e^{(\ell)}_{n1}={\rm E}[\Xb^{(\ell)\prime}_{ni}\thetab_r^{(\ell)}]=n^{-\ell/4}/\sqrt{p}$, for  
(i)~$\ell=0$ (away from uniformity), 
(ii) $\ell=1$ (beyond contiguity),
(iii) $\ell=2$ (under contiguity), and
{\protect\linebreak}
(iv) $\ell=3$ (under strict contiguity). In each case, the value~$r=0$ corresponds to the null hypothesis, whereas $r=1,2,\ldots,6$ provide increasingly severe alternatives. Some asymptotic power curves are plotted in dashed lines.
}
\label{hist2} 
\end{center}
\end{figure}

The resulting rejection frequencies of the following tests for~$\mathcal{H}_0\n:\thetab_n=\thetab_0$, all performed at nominal level~$5\%$, are plotted in Figure~\ref{hist2}~: (1) the Watson test~$\phi_n^W$ in~(\ref{Watsontest}), (2) the Wald test~$\phi_n^S$ in~(\ref{Waldtest}), (3) the ``contiguity-Wald" test~$\phi_{n,\xi;{\rm C}}^S$ (resp., (4) the ``strict-contiguity-Wald" test~$\phi_{n;{\rm SC}}^S$) rejecting the null when the Wald test statistic~$S_n$ exceeds the upper-$\alpha$ quantile of the asymptotic null distribution in Theorem~\ref{LDnullasymp}(iii) (resp., in Theorem~\ref{LDnullasymp}(iv)). For~(3)-(4), these critical values were estimated from a random sample of size~$10^6$ drawn from the corresponding asymptotic null distribution (note that, for~$\phi_{n,\xi;{\rm C}}^S$, estimating the critical value, hence conducting this test, not only requires assuming that we are in the contiguous regime, but further requires knowing the true value of the corresponding locality parameter~$\xi$, which is of course unrealistic). In none of the asymptotic regimes are all four tests considered. Asymptotic powers of the Watson test are also plotted (for~$\ell=0,1$ and for~$\ell=3$, these asymptotic powers coincide with those of the Wald test and of the ``strict-contiguity-Wald test", respectively). 

Results are in a very good agreement with Theorem~\ref{LDnonnullasymp}. While the Watson and Wald tests provide essentially the same empirical powers away from uniformity ($\ell=0$), they show opposite non-null behaviours under contiguity~($\ell=2$). There, the Watson test detects alternatives of the form~$\thetab_n=\thetab_0+\taub$ (with the non-monotonic power pattern described when commenting on the non-centrality parameter in Theorem~\ref{LDnonnullasymp}(iii) above), while the Wald test basically never rejects such alternatives. Interestingly, the contiguity-oracle test~$\phi_{n,\xi;{\rm C}}^S$ behaves very poorly as well, since its empirical rejection frequencies, in line with the  corresponding asymptotic powers, are \emph{uniformly smaller} than the nominal level~$\alpha$. Finally, under strict contiguity, the Watson and ``strict-contiguity-Wald" tests, in accordance with Theorems~\ref{LDnullasymp}-\ref{LDnonnullasymp}, provide empirical rejection frequencies virtually equal to the nominal level, while the standard Wald test basically never rejects the null there. 

\section{Adaptive Le Cam optimality} 
\label{LANfixedpsec}

Away from uniformity (that is, for~$\eta_n\equiv 1$), the Watson test shows non-trivial asymptotic powers against alternatives of the form~$\thetab_n=\thetab_0+n^{-1/2}\taub_n$, where the sequence~$\taub_n$ is~$O(1)$ but not~$o(1)$ and is such that~$\thetab_0+n^{-1/2}\taub_n\in\mathcal{S}^{p-1}$ for any~$n$. No tests can improve on this consistency rate, which is a consequence of the \emph{local asymptotic normality (LAN)} result derived in \cite{PaiVer2015a}. Better~: the same LAN result shows that, in the FvML case, the Watson test is locally asymptotically maximin, hence provides, in the Le Cam maximin sense, the best asymptotic powers that can be achieved in the FvML case. However, from the results in \cite{PaiVer2015a}, it is easy to conclude that, still away from uniformity, the optimality of the Watson test does not extend beyond the FvML setup. 

This raises natural questions in the vicinity of uniformity~: in the corresponding regimes (beyond contiguity, under contiguity, under strict contiguity), is the Watson test still rate-optimal? If it is, does it still enjoy optimality properties that are parallel to those stated above? To answer these questions, we derive the following LAN result (see the appendix for a proof).

\begin{Theor} 
\label{vicinityfixedpLAN}
Consider the sequence of parametric models~$\mathcal{P}\n_{\eta_n,\xi,f}=\{{\rm P}\n_{\thetab,\kappa_n,f}:\thetab\in\mathcal{S}^{p-1}\}$, with~$\kappa_n=\sqrt{p}\eta_n\xi+o(\eta_n)$ as~$n\to\infty$, where $(\eta_n)$ is a sequence in~$\R^+$ that is~$o(1)$, $\xi>0$ is fixed, and~$f:\R\to\R^+$ is monotone increasing, twice differentiable at~$0$, and satisfies~$f(0)=f'(0)=1$. 
If~$\sqrt{n}\eta_n\to \infty$ (beyond contiguity), let 
$$
\nu_n
=
\frac{1}{\sqrt{n}\eta_n} \, ,
\quad
\Deltab_{\thetab,\xi}\n:=
\xi \sqrt{np} (\mathbf{I}_p - \thetab\thetab\pr) \bar{\Xb}_n
\quad
\textrm{ and }
\quad
\Gammab_{\thetab,\xi}:=
\xi^2   (\mathbf{I}_p - \thetab\thetab\pr )
;
$$
if~$\sqrt{n}\eta_n\to 1$ (under contiguity), let 
$$
\nu_n\equiv 1,
\quad
\Deltab_{\thetab,\xi}\n:=
\xi \sqrt{np} \bar{\Xb}_n-\xi^2  \thetab
\quad
\textrm{ and }
\quad
\Gammab_{\thetab,\xi}:=
\xi^2  \mathbf{I}_p
;
$$
if~$\sqrt{n}\eta_n=o(1)$ (under strict contiguity), let 
$$
\nu_n\equiv 1,
\quad
\Deltab_{\thetab,\xi}\n:=\mathbf{0}
\quad
\textrm{ and }
\quad
\Gammab_{\thetab,\xi}:=
\mathbf{0}
.
$$
Let further~$(\taub_n)$ be a bounded sequence in~$\R^p$ that is not~$o(1)$ and that is such that~$\thetab+ \nu_n \taub_n\in{\mathcal S}^{p-1}$ for any~$n$.
Then, for any~$\thetab\in\mathcal{S}^{p-1}$, we have that, 
as~$\ny$ under~${\rm P}\n_{\thetab,\kappa_n,f}$, 
$$
\log \frac{d{\rm P}\n_{\thetab+\nu_n\taub_n,\kappa_n,f}}{d{\rm P}\n_{\thetab,\kappa_n,f}} 
=
\taub_n\pr 
\Deltab_{\thetab,\xi}\n
 - \frac{1}{2}
\taub_n\pr \Gammab_{\thetab,\xi} \taub_n
+
o_{\rm P}(1)
\quad
\textrm{and}
\quad
\Deltab_{\thetab,\xi} \n\stackrel{\mathcal{D}}{\to}\mathcal{N}_p(\mathbf{0},\Gammab_{\thetab,\xi})
.
$$
In other words, $\mathcal{P}\n_{\eta_n,\xi,f}$ is LAN, with central sequence~$\Deltab_{\thetab,\xi}\n$, Fisher information matrix~$\Gammab_{\thetab,\xi}$, and contiguity rate~$\nu_n$.     
\end{Theor}

Beyond contiguity, a locally asymptotically maximin test for~$\mathcal{H}\n_0:\thetab_n=\thetab_0$ is therefore rejecting the null at asymptotic level~$\alpha$ whenever
$$
Q_{\rm BC}\n
:=
(\Deltab_{\thetab_0,\xi}\n)\pr \,
\Gammab_{\thetab_0,\xi}^-
\Deltab_{\thetab_0,\xi} \n
=
np \bar{\Xb}_n'
(\mathbf{I}_p - \thetab_0\thetab_0\pr) 
\bar{\Xb}_n
>\chi^2_{p-1,1-\alpha}
,
$$
where $\Ab^-$ denotes the Moore--Penrose pseudoinverse of~$\Ab$.
Under~${\rm P}\n_{\thetab_0,\kappa_n,f}$, with~$\kappa_n=\sqrt{p}\eta_n\xi+o(\eta_n)$, where~$\eta_n=o(1)$ and~$\sqrt{n}\eta_n\to\infty$, we have
$
{\rm E}[(\Xb_{n1}'\thetab_0)^2]
=
1/p + 
o(1)
$; see~(\ref{tatatar}). Lemma~\ref{step1} (see the appendix) thus implies that 
$$
W_n=Q_{\rm BC}\n+o_{\rm P}(1)
$$
as~$n\to\infty$ under the same sequence of hypotheses --- hence also, from contiguity, under sequences of local alternatives of the form~${\rm P}\n_{\thetab_0+\taub_n/(\sqrt{n}\eta_n),\kappa_n,f}$. We conclude that, beyond contiguity, the Watson test remains locally asymptotically maximin. Actually, this optimality property,  quite remarkably, holds \emph{at virtually any~$f$} (that is, at any~$f$ meeting the conditions of Theorem~\ref{vicinityfixedpLAN}), which is in contrast with the fact that, away from uniformity, the Watson test is locally asymptotically maximin \emph{at the FvML only}. Now, by applying the Le Cam third lemma, the LAN result above allows us to derive the asymptotic distribution of the Watson test statistic under the sequences of local alternatives~${\rm P}\n_{\thetab_0+\taub_n/(\sqrt{n}\eta_n),\kappa_n,f}$; doing so actually confirms, in the present absolutely continuous setup, the non-null result obtained for~$W_n$ in Theorem~\ref{LDnonnullasymp}(ii).  

The story is different under contiguity. Proceeding as above, it may be tempting there to consider the test rejecting the null~$\mathcal{H}\n_0:\thetab_n=\thetab_0$ at asymptotic level~$\alpha$ whenever
\begin{equation}
\label{Oracletest}
Q_{\rm C;oracle}\n
:=
\Deltab_{\thetab_0,\xi}\n
(\Gammab_{\thetab_0,\xi})^{-}
\Deltab_{\thetab_0,\xi}\n
=
\|\sqrt{np} \bar{\Xb}_n-\xi  \thetab_0\|^2
>\chi^2_{p,1-\alpha},
\end{equation} 
based on the central sequence and information matrix obtained under contiguity (see Theorem~\ref{vicinityfixedpLAN}). This test, however, is much less satisfactory than the optimal test we just considered beyond contiguity. The reason is two-fold. First, as hinted by the notation in~(\ref{Oracletest}), this test is an oracle test, in the sense that it requires knowing the underlying value of the locality parameter~$\xi$. Second, the optimality properties of this test (if any) are unclear, due to the non-standard nature of the limiting experiment at hand. 

To comment on the latter point, note that the LAN result above, under contiguity, leads, for any fixed~$\thetab$, to a limiting experiment of the form
\begin{eqnarray}
\lefteqn{
\Big( 
\R^p, \mathcal{B}^p, \mathcal{P}_{\thetab,\xi}
=
\Big\{
{\rm P}^{\mathcal{N}_p(\xi^2\taub,\xi^2\mathbf{I}_p)}: \taub \in \R^p \textrm{ such that } \thetab+\taub\in\mathcal{S}^{p-1}
\Big\}
\Big)
}
\nonumber
\\[2mm]
& & \hspace{25mm}
=
\Big( 
\R^p, \mathcal{B}^p, \mathcal{P}_{\thetab,\xi}
=
\Big\{
{\rm P}^{\mathcal{N}_p(\xi^2\taub,\xi^2\mathbf{I}_p)}: \taub \in -\thetab+\mathcal{S}^{p-1}
\Big\}
\Big)
.
\label{limitexpns}
\end{eqnarray}
The problem of testing~$\mathcal{H}_0\n:\thetab_n=\thetab_0$ against~$\mathcal{H}_1\n:\thetab_n\neq \thetab_0$ translates, in the corresponding $\thetab_0$-limiting experiment, into the testing problem
\begin{equation}
\label{limprob}
\bigg\{
\begin{array}{l}
\mathcal{H}_0:\taub=\mathbf{0} \\[0mm]
\mathcal{H}_1:\taub\in (-\thetab_0+\mathcal{S}^{p-1})\setminus \{\mathbf{0}\},
\end{array}
\end{equation}
based on a single observation~$\Deltab$ from the $p$-variate normal distribution with mean~$\xi^2\taub$ and covariance matrix~$\xi^2\mathbf{I}_p$. 
While the limiting experiment in~(\ref{limitexpns}) is, as always in the LAN framework, a Gaussian shift experiment, the non-linear constraint $\taub \in -\thetab_0+\mathcal{S}^{p-1}$ on its location parameter makes this limiting experiment non-standard. And to the best of our knowledge, no \emph{globally} optimal test is known for the problem~(\ref{limprob}), irrespective of the optimality concept considered (leading to most powerful tests, maximin tests, most stringent tests, etc). This prevents the construction of locally asymptotically optimal tests in the corresponding sequence of experiments~$\mathcal{P}\n_{\eta_n,\xi,f}$ and makes unclear whether or not the test in~(\ref{Oracletest}) is optimal in some sense.

\begin{figure}[htbp!] 
\begin{center}
\includegraphics[width=82mm]{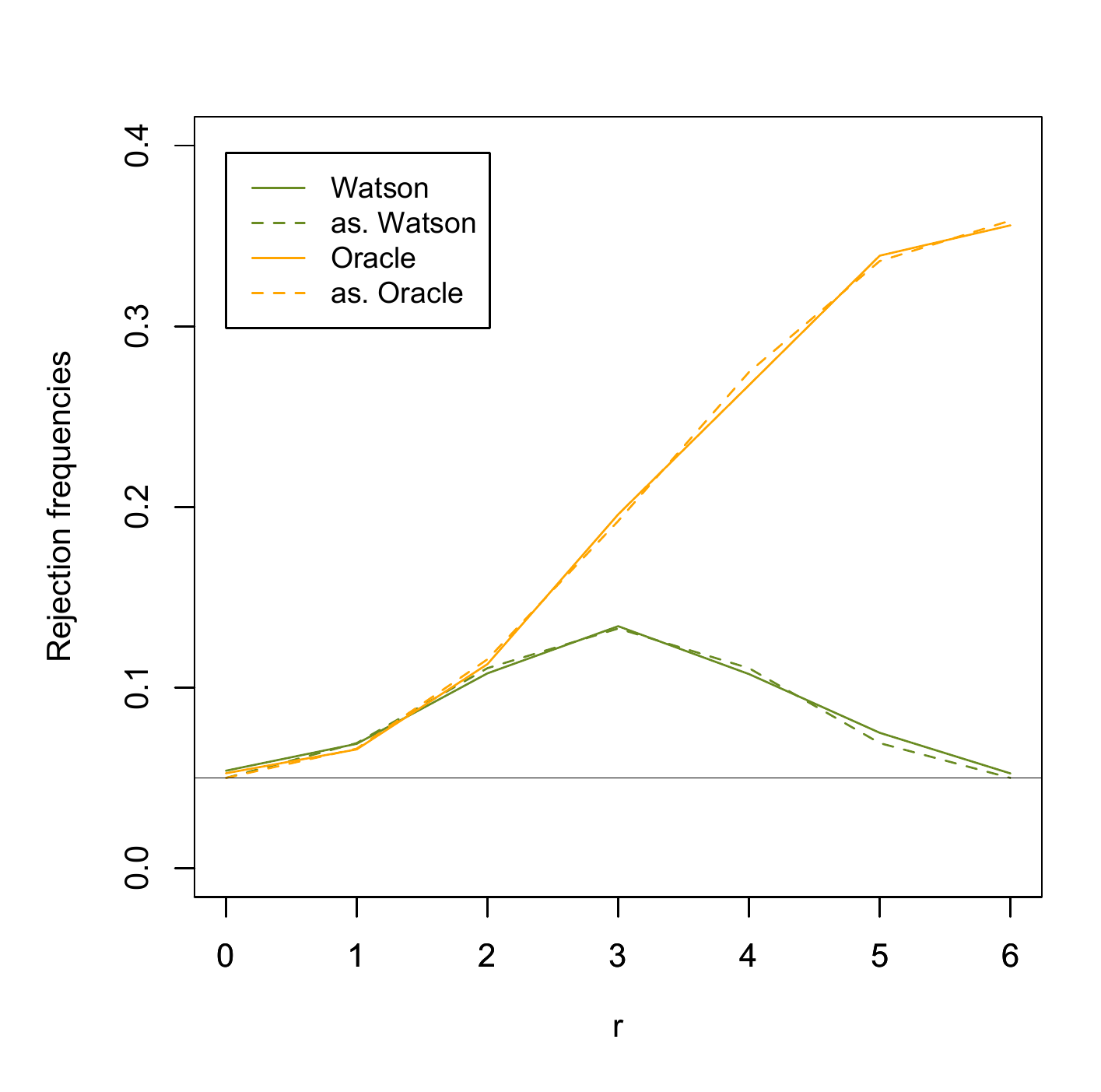}
\caption{
Rejection frequencies of the Watson test in~(\ref{Watsontest}) and of the oracle test in~(\ref{Oracletest}) when testing~$\mathcal{H}_0\n:\thetab_n=\thetab_0(=(0,0,1)')$, at nominal level~$\alpha=5\%$, on $M=10,\!000$ independent random samples of size~$n=200$ from the FvML distribution on~$\mathcal{S}^2$ with modal location~$\thetab_r=\thetab_0+\vb_r$ (see~(\ref{localter})) and a concentration~$\kappa$ such that~$e_{n1}={\rm E}[\Xb^{\prime}_{ni}\thetab_r]=1/\sqrt{np}$ (corresponding to the contiguous regime), for~$r=0$ (null hypothesis) and $r=1,2,\ldots,6$ (increasingly severe alternatives). The distributional setup is therefore the same as in the upper-right panel of Figure~\ref{hist2}. The dashed lines are the corresponding asymptotic power curves.
}
\label{figureLAN} 
\end{center}
\end{figure}

However, it is easy to show that the test rejecting the null of~(\ref{limprob}) whenever
$$
\frac{1}{\xi^2}
\,
\Deltab\pr
(\mathbf{I}_p - \thetab_0\thetab_0\pr) 
\Deltab
>\chi^2_{p-1,1-\alpha}
$$
is \emph{locally} maximin at level~$\alpha$ for the testing problem~(\ref{limprob}). As a corollary, under contiguity, the test rejecting the null~$\mathcal{H}_0\n:\thetab_n=\thetab_0$ at asymptotic level~$\alpha$ whenever 
$$
Q_{\rm C}\n
:=
\frac{1}{\xi^2}
\,
\Deltab_{\thetab_0,\xi}^{(n)\prime}
(\mathbf{I}_p - \thetab_0\thetab_0\pr) 
\Deltab_{\thetab_0,\xi}\n
=
np \bar{\Xb}_n'
(\mathbf{I}_p - \thetab_0\thetab_0\pr) 
\bar{\Xb}_n
>\chi^2_{p-1,1-\alpha}
$$
is \emph{bilocally asymptotically maximin}, where the term ``bilocally" refers to local-in-$\thetab$ and local-in-$\taub$ optimality (standard locally asymptotically optimal tests are associated with local-in-$\thetab$ optimality only). Since Lemma~\ref{step1} still ensures that, under contiguity,~$W_n=Q_{\rm C}\n+o_{\rm P}(1)$ under the null (hence also under sequences of contiguous local alternatives), the Watson test also enjoys this bilocal asymptotic optimality property. Interestingly, the oracle test in~(\ref{Oracletest}) does not enjoy the same optimality property, which is easily seen by comparing, as~$\|\taub\|\to 0$, the asymptotic powers of the Watson test (resulting from Theorem~\ref{LDnonnullasymp}(iii)) with those of the oracle test (obtained from the fact that, in view of the Le Cam third lemma,
\begin{equation}
\label{Oracletestpower}
Q_{\rm C;oracle}\n
\stackrel{\mathcal{D}}{\to}
\chi^2_{p}(\xi^2\|\taub\|^2)
\end{equation} 
under alternatives of the form~${\rm P}\n_{\thetab_0+\taub,\kappa_n,f}$, with~$\kappa_n=\sqrt{p/n}\,\xi+o(n^{-1/2})$). Of course, the oracle test might still outperform the Watson test for more severe alternatives, that is, for larger values of~$\|\taub\|$. This is actually the case, as we show in Figure~\ref{figureLAN} by comparing the corresponding asymptotic powers as well as the respective empirical rejection frequencies, obtained in a simulation exercise similar to the one in the upper-right panel of Figure~\ref{hist2} (see the caption of Figure~\ref{figureLAN} for details).


Beyond these Le Cam optimality issues, the LAN result  in Theorem~\ref{vicinityfixedpLAN} guarantees that the Watson test is at least rate-optimal under contiguity~: in the contiguous regime, the Watson test shows non-trivial asymptotic powers against $\thetab$-fixed alternatives of the form~${\rm P}\n_{\thetab_0+\taub,\kappa_n,f}$, with~$\kappa_n=\sqrt{p/n}\,\xi+o(n^{-1/2})$ (Theorem~\ref{LDnonnullasymp}(iii)), and no tests can detect less severe local alternatives of the form~${\rm P}\n_{\thetab_0+\nu_n\taub_n,\kappa_n,f}$, with~$\nu_n=o(1)$, $(\taub_n)$ bounded, and~$\kappa_n=\sqrt{p/n}\,\xi+o(n^{-1/2})$. Moreover, it is still so that, under contiguity, the local asymptotic powers of the Watson test in Theorem~\ref{LDnonnullasymp}(iii) can be obtained by applying the Le Cam third lemma with the central sequence and Fisher information matrix from (the contiguous-regime part of) Theorem~\ref{vicinityfixedpLAN} above (under contiguity, the same result can actually also be obtained by applying the third lemma to the LAN result in Theorem~3.1 of \cite{Cut2015}, which
is quite remarkable since this LAN result is with respect to~$\kappa$, unlike the one in Theorem~\ref{vicinityfixedpLAN} that is with respect to~$\thetab$). 

Finally, under strict contiguity, Theorem~\ref{vicinityfixedpLAN} implies that no asymptotic $\alpha$-level tests can detect even the most severe alternatives of the form~$\thetab_n=\thetab_0+\taub$, so that the Watson test may be considered optimal in this case, too. Of course, the optimality here is somewhat degenerate since the trivial $\alpha$-test, that randomly rejects the null with probability~$\alpha$, is also optimal under strict contiguity.

\section{Real data example}
\label{realsec}

\begin{figure}[htbp!]  
\begin{center}
\vspace{5mm}  
\includegraphics[width=49mm]{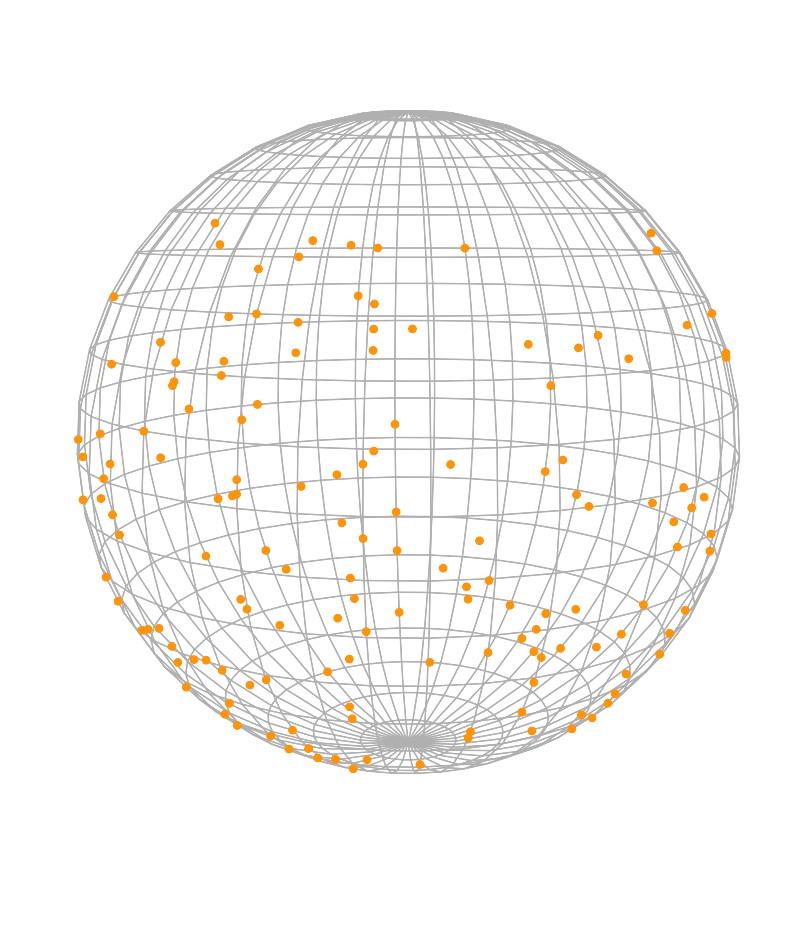}
\includegraphics[width=49mm]{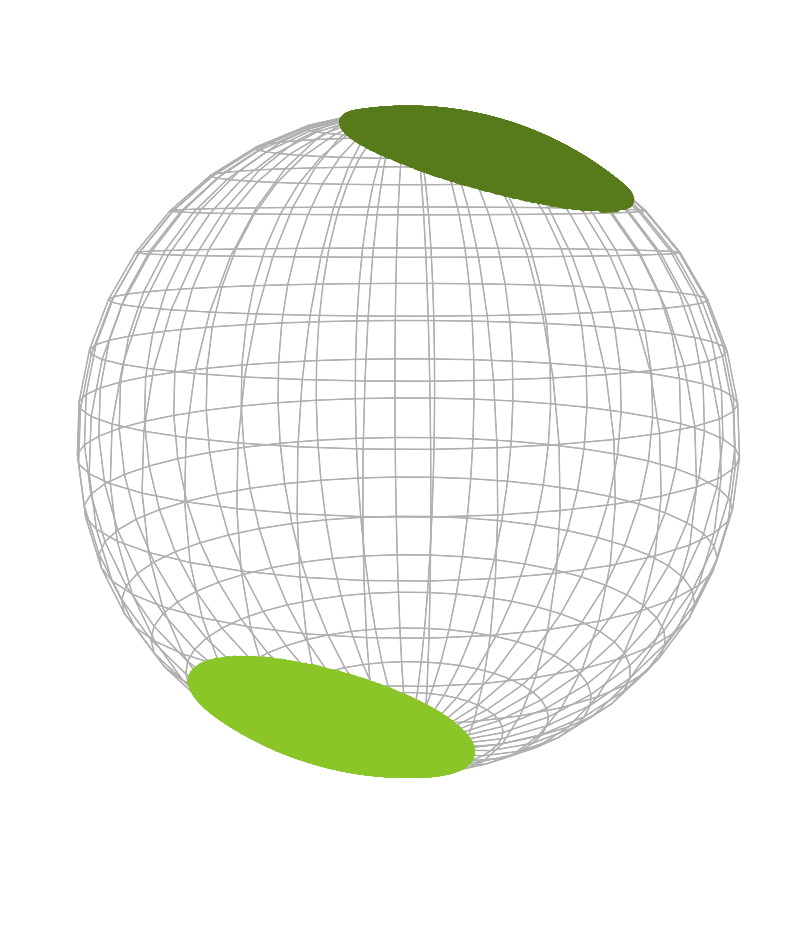}
\includegraphics[width=49mm]{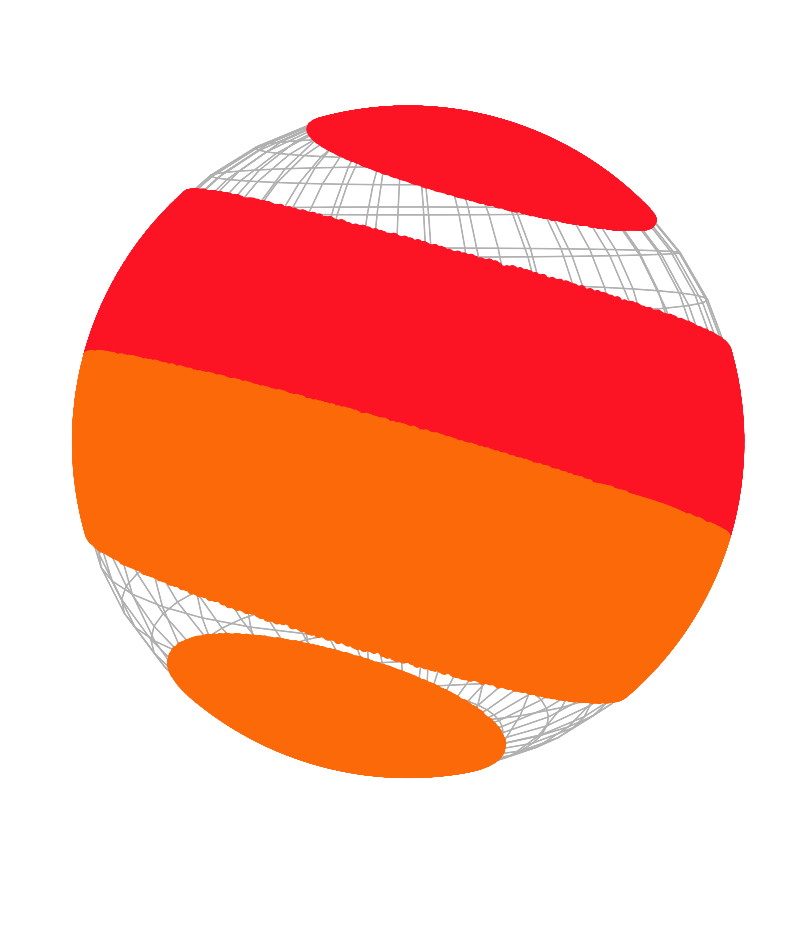}
\vspace{-3mm}
\caption{ 
(Left:) 
The $n=148$ measurements of cosmic ray directions from  \cite{Toyetal1965}.
(Middle:)
the asymptotic $95\%$-confidence zone~${\cal C}^{W}_{n;.95}$ obtained by inverting the Watson test.
(Right:)
the corresponding confidence zone, ${\cal C}^{S}_{n;.95}$, obtained from the Wald test. Both for the Watson and the Wald confidence zones, the symmetric component containing the point estimate~$\hat{\thetab}_n$ is shown with lighter colors (light green and orange, respectively).
}
\label{Zones} 
\end{center}
\end{figure}

In this section, we illustrate the practical relevance of our results on a cosmic ray data set. This data set, that was first used in \cite{Toyetal1965} to study primary cosmic rays in certain energy regions, has also been analysed, among others, in \citet[p.~102]{Fish87} and \cite{LSV14}. When applied to the $n=148$ arrival directions of cosmic rays at hand, the classical Rayleigh test of uniformity over~$\mathcal{S}^2$ rejects the null at asymptotic level~$5\%$; yet visual inspection of the left panel of Figure~\ref{Zones} below suggests that concentration is quite moderate, so that inference on the modal location~$\thetab$ may be delicate. We will compare, in the light of the results derived in the previous sections, the confidence zones for~$\thetab$ obtained by inverting the Watson and Wald tests.

Letting again~$\hat{\thetab}_n=\bar{\Xb}_n/\| \bar{\Xb}_n\|$, the Watson and Wald tests lead to the confidence zones (at asymptotic confidence level~$1-\alpha$)
\begin{equation}
\label{WatsonCZ}
{\cal C}^{W}_{n;1-\alpha}
:= 
\bigg\{ 
\thetab \in {\cal S}^{p-1}: 
W_n(\thetab)
:=
\frac{n (p-1)\bar{\Xb}_n\pr ({\bf I}_p- \thetab \thetab\pr) \bar{\Xb}_n}{1-\frac{1}{n} \sum_{i=1}^n(\Xb_{ni}\pr\thetab)^2}
\leq \chi^2_{p-1,1-\alpha} 
\bigg\}  
\end{equation}
and
\begin{equation}
\label{WaldCZ}
{\cal C}^{S}_{n;1-\alpha}
:=
\bigg\{ 
 \thetab \in {\cal S}^{p-1}:
S_n(\thetab)
:=
\frac{n(p-1) (\bar{\Xb}_{n}\pr\thetab)^2
\,
 \hat{\thetab}_n\pr({\bf I}_p- \thetab\thetab\pr)\hat{\thetab}_n}{1-\frac{1}{n} \sum_{i=1}^n(\Xb_{ni}\pr\thetab)^2}
\leq \chi^2_{p-1,1-\alpha} 
\bigg\} 
,
\end{equation}
respectively. Note that $W_n(\thetab)$ and~$S_n(\thetab)$ are respectively obtained from~(\ref{Watsontest}) and~(\ref{Waldtest}) by substituting~$\thetab$ for~$\thetab_0$, hence are the Watson and Wald test statistics to be used when testing that the modal location is equal to~$\thetab$.

Since~$p$ is small for the cosmic ray data set, it is computationally feasible to evaluate these confidence zones by simply considering a sufficiently fine grid over~$S^{p-1}$. The resulting confidence zones (at asymptotic confidence level~$1-\alpha=95\%$) are plotted in Figure~\ref{Zones}. Clearly, the Wald confidence zone is much larger than the Watson one. This arguably results from the fact that the Wald test is overly conservative in the vicinity of uniformity. In contrast, the Watson test, that was proved to be robust to arbitrarily mild departures from uniformity, provides more accurate confidence zones. 
We conclude that, in the present example showing little deviation from uniformity, the Watson and Wald procedures behave in perfect agreement with our asymptotic results in Theorem~\ref{LDnullasymp} and with their finite-sample illustration in Figure~\ref{hist1}.

Two further comments are in order~: 

\begin{itemize}
\item[(i)] The bipolar nature of the Watson/Wald confidence zones may be puzzling at first. However, the invariance of~$W_n(\thetab)$ and~$S_n(\thetab)$ under reflections of~$\thetab$ about the centre of~$\mathcal{S}^{p-1}$ directly implies that the confidence zones in~(\ref{WatsonCZ})-(\ref{WaldCZ}) are always symmetric with respect to this centre. In practice, of course, the ``symmetric components" of these confidence zones do play very different roles and it is natural to favour the one containing the spherical mean~$\hat{\thetab}_n$ (that is plotted in light colors in Figure~\ref{Zones}), even though this symmetric component alone is not a $95\%$-confidence zone for~$\thetab$.
\item[(ii)] The symmetric component of the Watson confidence zone containing the point estimate~$\hat{\thetab}_n$, namely the intersection between~${\cal C}^{\rm W}_{n;.95}$ and the hemisphere with pole~$\hat{\thetab}_n$, is made of a well-behaved connected region. In contrast, the corresponding Wald symmetric component is \emph{not} connected but rather is the union of a zone containing~$\hat{\thetab}_n$ and a zone containing the great circle orthogonal to~$\hat{\thetab}_n$. Inspection of~(\ref{WaldCZ}) makes it clear that this great circle will \emph{always} be part of the Wald confidence zone, which is of course undesirable (incidentally, this is also at the origin of the uniformly biased ``contiguity-Wald" power curve in the upper-right panel of Figure~\ref{hist2}). 
\end{itemize}

When the underlying distribution does not deviate much from uniformity, Watson confidence zones therefore outperform their Wald counterparts on all counts. It is remarkable that these two procedures that so far have been perceived as perfectly interchangeable (due to their asymptotic equivalence away from uniformity) behave so differently in the vicinity of uniformity.



\section{Summary}
\label{conclusec}

In the spherical location problem, the classical Watson test, unlike the Wald test based on the spherical mean, is robust to asymptotic scenarios in which the underlying distribution converges to the uniform distribution. Irrespective of the rate of this convergence (leading to the beyond contiguity, under contiguity, and under strict contiguity regimes), the Watson test exhibits the same asymptotic $(\chi^2_{p-1}$) distribution as under distributions that are fixed away from uniformity. The Watson test is also \emph{rate-adaptive}, in the sense that, irrespective of the regime considered, no tests can show non-trivial asymptotic powers against less severe alternatives than those detected by the Watson test. 

This test further enjoys excellent, Le Cam-type, optimality properties that can be summarized as follows~:
(i)
for distributions that are fixed away from uniformity, the Watson test is optimal under FvML densities. 
(ii)
Beyond contiguity, the Watson test is optimal under virtually any distribution, which is of course a stronger optimality property.
(iii)
Under contiguity, the Watson test is, uniformly in the underlying distribution, locally-in-$\taub$ optimal.
(iv)
Finally, under strict contiguity, the Watson test is optimal, but in a degenerate way, since, so close to uniformity, the trivial $\alpha$-test is also optimal. 
We conclude that, interestingly, the Watson test shows a ``non-monotonic" optimality pattern as one gets closer to uniformity. 

Throughout, Monte-Carlo studies showed that, irrespective of the regime considered, our asymptotic results  actually provide very accurate descriptions of the finite-sample behaviours of the Watson and Wald tests, even for moderate sample sizes (of the order of 100 or 200). Finally, the practical relevance of our results was illustrated on a real data set that shows little deviation from uniformity.

\appendix

\section{Proofs}
\label{appA}

Most proofs in this technical appendix are based on the so-called  \emph{tangent-normal decomposition} of~$\Xb_{ni}$, that is, on the expression~$\Xb_{ni}=u_{ni}\thetab_n+v_{ni} {\bf  S}_{ni}$,
where
\vspace{1mm}
$$
u_{ni}
:=
\Xb_{ni}'\thetab_n
,
\quad
v_{ni}
:=
 \sqrt{1- u_{ni}^2}
,
\quad
\textrm{and}
\quad
\Sb_{ni}:=
\left\{
\begin{array}{ll}
\displaystyle \frac{\Xb_{ni}-(\Xb_{ni}'\thetab_n)\thetab_n}{\|\Xb_{ni}-(\Xb_{ni}'\thetab_n)\thetab_n\|}
& 
\textrm{if } \Xb_{ni}\neq \thetab_n
\\[3mm]
\mathbf{0}
&
\textrm{otherwise}
.
\end{array}
\right.
\vspace{1mm}
$$


We start with the proof of Theorem~\ref{asymsphermean}. 
\vspace{2mm}

\noindent {\sc Proof of Theorem~\ref{asymsphermean}.} 
Fix a sequence of hypotheses~${\rm P}\n_{\thetab,F_n}$ such that
$$
e_{n1}
=
\frac{\eta_n\xi}{\sqrt{p}}  + 
o(\eta_n)
\quad
\textrm{and}
\quad
\tilde{e}_{n2}
=
\tilde{e}_{2}
+
o(1)
,
$$
with~$\xi,\tilde{e}_{2}>0$. This covers all cases considered in the statement of the theorem (if~$\eta_n=o(1)$, then we work with~$\tilde{e}_2=1/p$).  
Letting $e_{n1}=\eta_n\xi_n/\sqrt{p}$, where~$(\xi_n)\to \xi$, write
$$
\sqrt{n} \bar{\Xb}_n
-
\frac{\sqrt{n}\eta_n\xi_n}{\sqrt{p}}
\,
\thetab
=
\frac{1}{\sqrt{n}} \sum_{i=1}^n 
(u_{ni} - e_{n1}) \thetab
+
 \frac{1}{\sqrt{n}} \sum_{i=1}^n v_{ni} {\bf S}_{ni}
=:
V_{n1}+V_{n2}
,
$$
say. The Lindeberg CLT for triangular arrays yields
$$
\bigg(
\!
\begin{array}{c}
V_{n1}\\[-1mm]
V_{n2}
\end{array}
\!
\bigg)
\stackrel{\mathcal{D}}{\to}
\mathcal{N}
\Bigg(
\bigg(
\!
\begin{array}{c}
\mathbf{0} \\[-1mm]
\mathbf{0}
\end{array}
\!
\bigg)
,
\bigg(
\begin{array}{cc}
\tilde{e}_2 \thetab \thetab \pr & \mathbf{0} \\[-1mm]
\mathbf{0} & \frac{d}{p} ({\bf I}_p- \thetab\thetab\pr)
\end{array}
\bigg)
\Bigg)
,
$$
where we let
\begin{equation}
\label{defd}
d:=
\frac{1-e_2}{1-1/p}
,
\ \ \textrm{ with } 
e_2:=\tilde{e}_2+\Big(\lim_{n\to\infty} e_{n1}\Big)^2
.
\end{equation}
Consequently,  
\begin{equation} 
\label{theone1}
\Yb_n^{\xi,\eta}
:=
\sqrt{n}\eta_n
(
\eta_n^{-1}\sqrt{p} \bar{\Xb}_n- 
\xi_n
 \thetab
)
 \stackrel{\mathcal{D}}{\to}
{\cal N} \big(
{\bf 0}
, 
p\tilde{e}_2 \, \thetab\thetab \pr 
+
d
 ({\bf I}_p- \thetab\thetab\pr)
\big) 
.
\end{equation}
Parts~(iii)-(iv) of the result directly follow (note that we have~$\tilde{e}_2=e_2=1/p$ in these cases), and we may thus focus on Parts~(i)-(ii). Applying the uniform delta method (see Theorem~3.8 in \citealp{van1998}) to~(\ref{theone1}), with the mapping $\xb \mapsto \xb/ \| \xb \|$, then yields
\begin{eqnarray}
\sqrt{n}\eta_n
(
\hat{\thetab}_n- \thetab
)
&=&
\sqrt{n}\eta_n
\xi^{-1}({\bf I}_p- \thetab\thetab\pr)
(
\eta_n^{-1}\sqrt{p} \bar{\Xb}_n- 
\xi_n
 \thetab
)
+o_{\rm P}(1)
\nonumber
\\[2mm]
&=&
\xi^{-1}({\bf I}_p- \thetab\thetab\pr)
\Yb_n^{\xi,\eta}
+o_{\rm P}(1)
\label{delt}
,
\end{eqnarray}
which, by using~(\ref{theone1}) again, establishes the result.
\cqfd



 \begin{Lem} 
\label{step1}
Let $(\thetab_n)$ be an arbitrary sequence in~$\mathcal{S}^{p-1}$ and~$(F_n$) be an arbitrary sequence of cumulative distribution functions on~$[-1,1]$. Then 
$\frac{1}{n} \sum_{i=1}^n u_{ni}^2={\rm E}[(\Xb_{n1}'\thetab_n)^2]+ o_{\rm P}(1)$ as~$n\to\infty$ under~${\rm P}^{(n)}_{\thetab_n,F_n}$ $($where the expectation is evaluated under~${\rm P}^{(n)}_{\thetab_n,F_n})$. 
\end{Lem}

\noindent {\sc Proof.}
Since~$\sup_n {\rm E}[|u_{ni}|^2]\leq 1$, the result readily follows from the weak law of large numbers for triangular arrays.
\cqfd
\vspace{3mm}

\noindent {\sc Proof of Theorem~\ref{LDnullasymp}.} 
Fix a sequence of hypotheses~${\rm P}\n_{\thetab_0,F_n}$ such that
$$
e_{n1}
=
\frac{\eta_n\xi}{\sqrt{p}}  + 
o(\eta_n)
\quad
\textrm{and}
\quad
\tilde{e}_{n2}
=
\tilde{e}_2
+
o(1)
,
$$
with~$\xi,\tilde{e}_2>0$. As in the proof of Theorem~\ref{asymsphermean}, we restrict to~$\tilde{e}_2=1/p$ whenever~$\eta_n=o(1)$. Write then  $e_{n1}=\eta_n\xi_n/\sqrt{p}$, where~$(\xi_n)\to \xi$. All derivations in the proof of Theorem~\ref{asymsphermean} then hold, with~$\thetab$ replaced with~$\thetab_0$ everywhere. In particular, (\ref{theone1}) yields
\begin{equation} 
\label{theone4}
\Bigg(
\begin{array}{c}
\thetab_0\pr \Yb_n^{\xi,\eta}
\\[-1mm]
({\bf I}_p- \thetab_0\thetab_0\pr) \Yb_n^{\xi,\eta}
\end{array}
\Bigg)
 \stackrel{\mathcal{D}}{\to}
{\cal N} \bigg(
\bigg(
\!
\begin{array}{c}
\mathbf{0} \\[-1mm]
\mathbf{0}
\end{array}
\!
\bigg)
, 
\bigg(
\begin{array}{cc}
p\tilde{e}_2 &  \mathbf{0} 
\\[-1mm] 
 \mathbf{0}  & d ({\bf I}_p- \thetab_0\thetab_0\pr)
\end{array}
\bigg)
\bigg)
,
\end{equation}
where~$d$ is as in~(\ref{defd}), still with $e_2:=\tilde{e}_2+(\lim_{n\to\infty} e_{n1})^2$. Note then that the Watson statistic satisfies
\vspace{1mm}
\begin{equation}
\label{asWn}
W_n 
=
 \frac{(p-1)
(\Yb_n^{\xi,\eta})' 
({\bf I}_p- \thetab_0 \thetab_0\pr) 
\Yb_n^{\xi,\eta}
}{p(1-\frac{1}{n} \sum_{i=1}^n u_{ni}^2)} 
=
\frac{1}{d}
(\Yb_n^{\xi,\eta})' 
({\bf I}_p- \thetab_0 \thetab_0\pr) 
\Yb_n^{\xi,\eta}
+o_{\rm P}(1)
,
\end{equation}
where we used Lemma~\ref{step1}.
It follows that~$W_n \stackrel{\mathcal{D}}{\to} \chi^2_{p-1}$ in all four cases~(i)-(iv).

We then turn to the Wald statistic~$S_n$ and consider first the cases~(i)-(ii). Note that~(\ref{delt}) entails 
\begin{eqnarray*}
n\eta_n^2
\,
\hat{\thetab}_n\pr ({\bf I}_p- \thetab_0\thetab_0\pr) \hat{\thetab}_n
&=& 
n\eta_n^2
 (\hat{\thetab}_n- \thetab_0)\pr({\bf I}_p- \thetab_0\thetab_0\pr) (\hat{\thetab}_n- \thetab_0)
\\[1mm]
&=& 
\xi^{-2}
(\Yb_n^{\xi,\eta})'
 ({\bf I}_p- \thetab_0\thetab_0\pr) 
\Yb_n^{\xi,\eta} + o_{\rm P}(1)
.
\end{eqnarray*}
This leads to
\begin{eqnarray*}
\lefteqn{
S_n
=
 \frac{n(p-1)(\bar{\Xb}_n'\thetab_0)^2}{1-\frac{1}{n} \sum_{i=1}^n u_{ni}^2} 
\frac{1}{n\eta_n^2\xi^2}
(\Yb_n^{\xi,\eta})'
 ({\bf I}_p- \thetab_0\thetab_0\pr) 
\Yb_n^{\xi,\eta} + o_{\rm P}(1)
}
\\[1mm]
& &
\hspace{-2mm}
=
\frac{1}{n\eta_n^2\xi^2}
(\sqrt{np}\, \bar{\Xb}_n'\thetab_0)^2
W_n+ o_{\rm P}(1)
=
\frac{1}{n\eta_n^2\xi^2}
(\thetab_0' \Yb_n^{\xi,\eta}+\sqrt{n}\eta_n\xi_n)^2
W_n+ o_{\rm P}(1)
,
\end{eqnarray*}
which shows that~$S_n=W_n+ o_{\rm P}(1)$, hence proves~(i)-(ii) for~$S_n$. 

Turning to cases~(iii)-(iv), note that (\ref{theone1}) rewrites
$
\sqrt{np}\, \bar{\Xb} 
 \stackrel{\mathcal{D}}{\to}
{\cal N} 
\big(
\lambda_\xi
 \thetab_{0}
, 
{\bf I}_p
\big) 
,
$
with $\lambda_\xi=\xi$ and $\lambda_\xi=0$ in case~(iii) and in case~(iv), respectively (recall that~$\tilde{e}_2=e_2=1/p$ in these cases). Hence,
\begin{eqnarray*}
S_n
&=&
 \frac{(p-1)(\sqrt{np}\,\thetab_0'\bar{\Xb}_n)^2}{p(1-\frac{1}{n} \sum_{i=1}^n u_{ni}^2)}
\,
\hat{\thetab}_n\pr ({\bf I}_p- \thetab_0\thetab_0\pr) \hat{\thetab}_n
\\[1mm]
&=&
(\sqrt{np}\,\thetab_0'\bar{\Xb}_n)^2
\,
\frac{[\sqrt{np} \bar{\Xb}_n]\pr ({\bf I}_p- \thetab_0\thetab_0\pr) [\sqrt{np} \bar{\Xb}_n ]}{\|\sqrt{np} \bar{\Xb}_n \|^2}
+o_{\rm P}(1)
\\[1mm]
&=&
(\sqrt{np}\,\thetab_0'\bar{\Xb}_n)^2
\,
\frac{\|\sqrt{np}({\bf I}_p- \thetab_{0}\thetab_{0}\pr)\bar{\Xb}_n \|^2}{\|(\sqrt{np} \, \thetab_0\pr\bar{\Xb}_n)\thetab_0+\sqrt{np}({\bf I}_p- \thetab_{0}\thetab_{0}\pr)\bar{\Xb}_n \|^2}
+o_{\rm P}(1)
\\[1mm]
&=&
\frac{(\sqrt{np}\,\thetab_0'\bar{\Xb}_n)^2\|\sqrt{np}({\bf I}_p- \thetab_{0}\thetab_{0}\pr)\bar{\Xb}_n \|^2}{(\sqrt{np} \, \thetab_0\pr\bar{\Xb}_n)^2+\|\sqrt{np}({\bf I}_p- \thetab_{0}\thetab_{0}\pr)\bar{\Xb}_n \|^2}
+o_{\rm P}(1)
.
\end{eqnarray*}
By combining~(\ref{asWn}) and~(\ref{theone1}), we then obtain   
\begin{equation}
\label{Waldecomp}
S_n
=
\frac{(\sqrt{np}\,\thetab_0'\bar{\Xb}_n)^2 W_n}{(\sqrt{np} \, \thetab_0\pr\bar{\Xb}_n)^2+W_n}
+o_{\rm P}(1)
=
\frac{(\thetab_0' \Yb_n^{\xi,\eta}+\sqrt{n}\eta_n\xi_n)^2 W_n}{(\thetab_0' \Yb_n^{\xi,\eta}+\sqrt{n}\eta_n \xi_n)^2+W_n}
+o_{\rm P}(1)
.
\end{equation}
From~(\ref{theone4})-(\ref{asWn}), it is seen that $Z_n:=\thetab_0' \Yb_n^{\xi,\eta}$ is asymptotically standard normal, $W_n$ is asymptotically $\chi^2_{p-1}$, and that~$Z_n$ and~$W_n$  are asymptotically mutually independent. This provides
$$
S_n
=
\frac{(Z_n+\lambda_\xi)^2}{(Z_n+\lambda_\xi)^2+W_n}
\,
 W_n
+o_{\rm P}(1)
=
\bigg(
1+\frac{W_n}{(Z_n+\lambda_\xi)^2}
\bigg)^{-1}
\,
 W_n
+o_{\rm P}(1)
,
$$
which establishes the result.
\cqfd
\vspace{3mm}


\noindent {\sc Proof of Theorem~\ref{LDnonnullasymp}.} 
Fix a sequence of hypotheses~${\rm P}\n_{\thetab_n,F_n}$ such that
$$
e_{n1}
=
\frac{\eta_n\xi}{\sqrt{p}}  + 
o(\eta_n)
\quad
\textrm{and}
\quad
\tilde{e}_{n2}
=
\tilde{e}_2
+
o(1)
,
$$
with~$\xi,\tilde{e}_2>0$ and $\thetab_n:=\thetab_0+\nu_n\taub_n$, where $(\nu_n)$ is as in the statement of Theorem~\ref{LDnonnullasymp} (we still restrict to~$\tilde{e}_2=1/p$ whenever~$\eta_n=o(1)$). Letting~$e_{n1}=\eta_n\xi_n/\sqrt{p}$, where~$(\xi_n)\to \xi$, and proceeding as in the proof of Theorem~\ref{asymsphermean}, 
we can write
$$
\sqrt{n} \bar{\Xb}_n
-
\frac{\sqrt{n}\eta_n \xi_n}{\sqrt{p}} \,\thetab_n
=
\frac{1}{\sqrt{n}} \sum_{i=1}^n 
(u_{ni} - e_{n1}) \thetab_n
+
 \frac{1}{\sqrt{n}} \sum_{i=1}^n v_{ni} {\bf S}_{ni}
=:
W_{n1}+W_{n2}
,
$$
say, where~${\bf S}_{ni}$ is now based on~$\thetab_n$. Under~${\rm P}\n_{\thetab_n,F_n}$, 
$$
\bigg(
\!
\begin{array}{c}
W_{n1}\\[-1mm]
W_{n2}
\end{array}
\!
\bigg)
\stackrel{\mathcal{D}}{\to}
\mathcal{N}
\Bigg(
\bigg(
\!
\begin{array}{c}
\mathbf{0} \\[-1mm]
\mathbf{0}
\end{array}
\!
\bigg)
,
\bigg(
\begin{array}{cc}
\tilde{e}_2 \thetab_{0\alpha}\thetab_{0\alpha} \pr & \mathbf{0} \\[-1mm]
\mathbf{0} & \frac{d}{p}
({\bf I}_p- \thetab_{0\alpha}\thetab_{0\alpha}\pr)
\end{array}
\bigg)
\Bigg)
,
$$
where we let~$\thetab_{0\alpha}:=\thetab_{0}+\delta\taub$, with~$\delta:=1$ if~$\nu_n\equiv 1$ (under contiguity or under strict contiguity) and $\delta:=0$ otherwise (away from contiguity or beyond contiguity). 
Parallel to~(\ref{theone1}), we obtain
\begin{equation} 
\label{theone1nonull}
\Yb_n^{\xi,\eta}
:=
\sqrt{n}\eta_n
(
\eta_n^{-1}\sqrt{p} \bar{\Xb}- 
\xi_n
 \thetab_{0}
)
 \stackrel{\mathcal{D}}{\to}
{\cal N} \big(
\lambda_\xi
\taub
, 
p\tilde{e}_2 \, \thetab_{0\alpha}\thetab_{0\alpha} \pr 
+
d
 ({\bf I}_p- \thetab_{0\alpha}\thetab_{0\alpha}\pr)
\big) 
,
\end{equation}
where~$\lambda_\xi$ is as in the proof of Theorem~\ref{LDnullasymp}. Letting $\taub_\thetab:=({\bf I}_p- \thetab_0\thetab_0\pr)\taub
=\taub-\delta(\thetab_0'\taub)\thetab_0$, this provides
\begin{equation} 
\label{theone4nonnull}
\Bigg(
\begin{array}{c}
\thetab_0\pr \Yb_n^{\xi,\eta}
\\[-2mm]
({\bf I}_p- \thetab_0\thetab_0\pr) \Yb_n^{\xi,\eta}
\end{array}
\Bigg)
 \stackrel{\mathcal{D}}{\to}
{\cal N} \bigg(
\lambda_\xi
\bigg(
\!
\begin{array}{c}
 \thetab_0\pr\taub \\[-2mm]
\taub_\thetab
\end{array}
\!
\bigg)
, 
\Sigb
\bigg)
,
\end{equation}
with
$$
\Sigb
:=
\bigg(
\begin{array}{cc}
d+(p\tilde{e}_2-d)(1+\delta\thetab_0'\taub)^2 &  \mathbf{0} 
\\[-2mm] 
 \mathbf{0}  & d({\bf I}_p- \thetab_0\thetab_0\pr)+(p\tilde{e}_2-d)\delta\taub_\thetab\taub_\thetab'
\end{array}
\bigg)
\bigg)
.
$$

Note that the Watson statistic still satisfies 
\vspace{1mm}
$
W_n 
=
\Tb_n\pr\Tb_n+o_{\rm P}(1)
$
(see~(\ref{asWn})), where $\Tb_n:=d^{-1/2}
({\bf I}_p- \thetab_0 \thetab_0\pr) 
\Yb_n^{\xi,\eta}$. 
We then readily obtain
$
W_n
\stackrel{\mathcal{D}}{\to} 
\chi^2_{p-1}
\big(
\lambda_\xi^2 \|\taub_\thetab\|^2/d
\big)
,
$
where the asymptotic distribution rewrites
$
\chi^2_{p-1}
\big(
\xi^2 \|\taub\|^2/d
\big)
,
$
$
\chi^2_{p-1}
\big(
\xi^2 \|\taub\|^2
\big)
,
$
$
\chi^2_{p-1}
\big(
\xi^2 \|\taub_\thetab\|^2
\big)
,
$
and
$
\chi^2_{p-1}
,
$
in cases~(i), (ii), (iii), and~(iv), respectively. Since a direct computation shows that~$\xi^2 \|\taub_\thetab\|^2$ coincides with the non-centrality parameter in Part~(ii) of the result, this completes the proof for the Watson test.

Turning to the Wald test statistic~$S_n$. For cases~(i)-(ii), the exact same reasoning as in the proof of Theorem~\ref{LDnullasymp}, this time applied to~(\ref{theone1nonull}), yields that $S_n=W_n+ o_{\rm P}(1)$ under the sequence of hypotheses considered, which yields the result. Now, in cases~(iii)-(iv), the result in~(\ref{Waldecomp}), or equivalently,
$$
S_n
=
\bigg(
1+\frac{W_n}{(\thetab_0' \Yb_n^{\xi,\eta}+\lambda_\xi)^2}
\bigg)^{-1}
\,
 W_n
+o_{\rm P}(1)
,
$$
still holds under the sequence of hypotheses considered. The results in cases~(iii)-(iv) then directly follow from~(\ref{theone4nonnull}) --- in case~(iii), recall indeed that~$\thetab_{0}\pr \taub=-\| \taub \|^2/2$.
\cqfd
\vspace{3mm}

\noindent {\sc Proof of Theorem~\ref{vicinityfixedpLAN}.} 
\emph{(Under contiguity).}
Note that, in the contiguous regime, $\kappa_n=\sqrt{p/n}\,\xi+o(n^{-1/2})$ as~$n\to\infty$. Writing~$\Delta_{\thetab}\n:=\sqrt{np}\,\bar{\Xb}_n'\thetab$, Theorem~3.1 in~\cite{Cut2015} then implies that, for any sequence~$(\thetab_n)$ in~$\mathcal{S}^{p-1}$, 
$$
\log \frac{d{\rm P}\n_{\thetab_n,\kappa_n,f}}{d{\rm P}\n_{0}} 
=
\xi \Delta_{\thetab_n}\n - \frac{\xi^2}{2}
+
o_{\rm P}(1)
, 
$$
as~$n\to\infty$ under~${\rm P}\n_{0}$. Using~(\ref{spherconstraint}), it then readily follows that 
\begin{eqnarray*}
\lefteqn{
\log \frac{d{\rm P}\n_{\thetab+\taub_n,\kappa_n,f}}{d{\rm P}\n_{\thetab,\kappa_n,f}} 
=
\log \frac{d{\rm P}\n_{\thetab+\taub_n,\kappa_n,f}}{d{\rm P}\n_{0}} 
-
\log \frac{d{\rm P}\n_{\thetab,\kappa_n,f}}{d{\rm P}\n_{0}} 
}
\\[2mm]
& & 
\hspace{3mm}
=
\
\xi
(
\Delta_{\thetab+\taub_n}\n
-
\Delta_{\thetab}\n
)
+
o_{\rm P}(1)
=
\xi  \sqrt{np}\,\bar{\Xb}_n'\taub_n
+
o_{\rm P}(1)
\\[2mm]
& & 
\hspace{3mm}
=
\
\taub_n\pr \Deltab_{\thetab,\xi}\n
+
\xi^2 \taub_n\pr\thetab
+ 
o_{\rm P}(1)
=
\taub_n\pr \Deltab_{\thetab,\xi}\n
 - \frac{1}{2}
\taub_n\pr \Gammab_{\thetab,\xi} \taub_n
+ 
o_{\rm P}(1)
\end{eqnarray*}
as~$\ny$ under~${\rm P}_0\n$, hence, from contiguity, also under~${\rm P}\n_{\thetab,\kappa_n,f}$. 
Now, (5.2)-(5.3) in~\cite{Cut2015} show that, under~${\rm P}\n_{\thetab_n,\kappa_n,f}$, with~$\kappa_n=\sqrt{p}\,\eta_n\xi+o(\eta_n)$ (where~$\eta_n=o(1)$), one has
$$
e_{n1}
=
\frac{\eta_n\xi}{\sqrt{p}}  + 
o(\eta_n)
\quad
\textrm{and}
\quad 
\tilde{e}_{n2}
=
\frac{1}{p}   + 
o(1)
.
$$
Consequently, (\ref{theone1}) applies and provides $\sqrt{np} \bar{\Xb}- \xi \thetab \stackrel{\mathcal{D}}{\to}{\cal N} (\mathbf{0}, {\bf I}_p)$ as~$\ny$ under~${\rm P}\n_{\thetab,\kappa_n,f}$. This establishes the result in the contiguous regime.  
\vspace{2mm}


\emph{(Under strict contiguity).}
From~(A.7) in~\cite{Cut2015} implies that, for any sequence~$(\thetab_n)$ in~$\mathcal{S}^{p-1}$, we learn that
$$
{\rm E}\bigg[
\bigg(
\log \frac{d{\rm P}\n_{\thetab_n,\kappa_n,f}}{d{\rm P}\n_{0}} 
\bigg)^2
\bigg]
=
O(n^2\kappa_n^4)
+
O(n\kappa_n^2)
$$
as~$n\to\infty$ under~${\rm P}\n_{0}$. In the strictly contiguous case, this readily provides
$$
\log \frac{d{\rm P}\n_{\thetab+\taub_n,\kappa_n,f}}{d{\rm P}\n_{\thetab,\kappa_n,f}} 
=
\log \frac{d{\rm P}\n_{\thetab+\taub_n,\kappa_n,f}}{d{\rm P}\n_{0}} 
-
\log \frac{d{\rm P}\n_{\thetab,\kappa_n,f}}{d{\rm P}\n_{0}} 
=
o_{\rm P}(1)
$$
as~$n\to\infty$ under~${\rm P}\n_{0}$. The result then follows from the mutual contiguity of~${\rm P}\n_{0}$ and~${\rm P}\n_{\thetab,\kappa_n,f}$. 
\vspace{-4mm}


\emph{(Beyond contiguity).} 
Write
$$
\log \frac{d{\rm P}\n_{\thetab+\nu_n\taub_n,\kappa_n,f}}{d{\rm P}\n_{\thetab,\kappa_n,f}} 
=
\sum_{i=1}^n 
\big(
\log f(\kappa_n u_{ni}+\kappa_n \nu_n\taub_n'\Xb_{ni})
-
\log f(\kappa_n u_{ni})
\big)
=
 L_{n1}+L_{n2}+L_{n3},
$$
with
$$
L_{n1}:=n \kappa_n\nu_n \taub_n' \bar{\Xb}_n,
\qquad
L_{n2}:=\kappa_n\nu_n \taub_n' 
\sum_{i=1}^n (\varphi_f(\kappa_n u_{ni})-1)\Xb_{ni}
$$
and
$$
L_{n3}:=
\sum_{i=1}^n 
\big(
\log f(\kappa_n u_{ni}+\kappa_n \nu_n\taub_n'\Xb_{ni})
-
\log f(\kappa_n u_{ni})
-
\kappa_n\nu_n  \varphi_f(\kappa_n u_{ni}) \taub_n' \Xb_{ni}
\big)
.
$$

The result then follows from the following lemma.
\cqfd

\begin{Lem}
\label{lemdedingue}
Let the assumptions of Theorem~\ref{vicinityfixedpLAN} hold and restrict to the case where $\sqrt{n}\eta_n\to\infty$. Then, as~$\ny$ under~${\rm P}\n_{\thetab,\kappa_n,f}$, 
(i) $L_{n1}= \taub_n' 
\Deltab\n_{\thetab,\xi}
-
\frac{1}{2} 
\taub_n\pr\Gammab_{\thetab,\xi}
\taub_n+o_{\rm P}(1)
$, where $\Deltab\n_{\thetab,\xi}
 \stackrel{\mathcal{D}}{\to}
{\cal N} \big(
{\bf 0}
, 
\Gammab_{\thetab,\xi} \big)$;
(ii) 
$
L_{n2}
=
-\frac{\xi^2}{2}
\varphi_f'(0)
\|\taub_n\|^2
+o_{\rm P}(1)
;
$
(iii)
$
L_{n3}
=
\frac{\xi^2}{2}
\varphi_f'(0)
\|\taub_n\|^2
+o_{\rm P}(1)
.
$
\end{Lem}

{\sc Proof of Lemma~\ref{lemdedingue}.} 
Throughout this proof, we write~$\kappa_n=\sqrt{p}\eta_n\xi_n$, where~$\xi_n\to\xi$. All expectations, variances, and stochastic convergence statements will be under~${\rm P}\n_{\thetab,\kappa_n,f}$.

(i) Since~$\eta_n$ is~$o(1)$, we still have that
$
e_{n1}
=
\kappa_n/p + 
o(\kappa_n)
$
and 
$
\tilde{e}_{n2}
=
1/p + 
o(1)
$
under~${\rm P}\n_{\thetab,\kappa_n,f}$. Consequently, (\ref{theone1}) applies and provides
$
\sqrt{np}\bar{\Xb}_n- 
\sqrt{n}\eta_n\xi_n
 \thetab
 \stackrel{\mathcal{D}}{\to}
{\cal N} \big(
{\bf 0}
, 
{\bf I}_p
\big),  
$
which implies that 
\begin{equation}
\label{agains}
\sqrt{n}(\bar{\Xb}_n- 
e_{n1} \thetab)
 \stackrel{\mathcal{D}}{\to}
{\cal N} \big(
{\bf 0}
, 
{\textstyle\frac{1}{p}}{\bf I}_p
\big).  
\end{equation}
Jointly with the fact that~$\taub_n\pr\thetab=o(1)$ in the present setup (see~(\ref{spherconstraint})), this implies that 
$$
\taub_n\pr
[
\Deltab\n_{\thetab,\xi}
-
\xi_n\sqrt{np}(\bar{\Xb}_n- 
e_{n1} \thetab)
]
=
-\xi_n \sqrt{np}\, \thetab' (\bar{\Xb}_n- 
e_{n1} \thetab) (\taub_n\pr\thetab)+o_{\rm P}(1)=o_{\rm P}(1).
$$
Using~(\ref{spherconstraint}), this readily yields  
$$
L_{n1}
=
n \kappa_n\nu_n  \taub_n' (\bar{\Xb}_n - e_{n1}\thetab)
+
n \kappa_n  e_{n1}  \nu_n \taub_n' \thetab
=
n \kappa_n\nu_n  \taub_n' (\bar{\Xb}_n - e_{n1}\thetab)
-
\frac{1}{2}
n \kappa_n  e_{n1}   \nu_n^2 \|\taub_n\|^2
$$
$$
=
 \taub_n'[ \xi_n   \sqrt{np} (\bar{\Xb}_n - e_{n1}\thetab) ]
-
\frac{1}{2}
n \kappa_n^2   \nu_n^2 \Big(\, \frac{1}{p} +o(1)\Big)   \|\taub_n\|^2
=
\taub_n' 
\Deltab\n_{\thetab,\xi}
-
\frac{1}{2} 
\taub_n\pr\Gammab_{\thetab,\xi}
\taub_n+o_{\rm P}(1)
.
$$
Finally, the asymptotic normality result for~$\Deltab\n_{\thetab,\xi}$ readily follows by premultiplying~(\ref{agains}) with~$\xi\sqrt{p} (\mathbf{I}_p-\thetab\thetab')$.  
\vspace{3mm}

(ii) Using the tangent-normal decomposition of~$\Xb_{ni}$, split~$L_{n2}$ into
\begin{eqnarray*}
L_{n2}
&=&
\kappa_n \nu_n 
\sum_{i=1}^n  (\varphi_f(\kappa_n u_{ni})-1)
u_{ni} (\taub_n\pr\thetab)
+
\taub_n\pr
\bigg(
\kappa_n \nu_n 
\sum_{i=1}^n  (\varphi_f(\kappa_n u_{ni})-1)
v_{ni} \Sb_{ni}
\bigg)
\\[2mm]
&=:& 
L_{n2a}+\taub_n\pr {\bf L}_{n2b},
\end{eqnarray*}
say. 
Since
\begin{eqnarray*}
{\rm E}[L_{n2a}]
&=&
n\kappa_n\nu_n
\,
{\rm E}[  (\varphi_f(\kappa_n u_{n1})-1)
u_{n1} ] \,(\taub_n\pr\thetab)
\\[1mm]
&= & 
-\frac{1}{2}n\kappa_n\nu_n^2 \|\taub_n\|^2
\,
c_{p,\kappa_n,f}
\int_{-1}^1 
(1-s^2)^{(p-3)/2}
(\varphi_f(\kappa_n s) -1)
 s f(\kappa_n s) 
\,ds
,
\end{eqnarray*}
we obtain
\begin{eqnarray}
{\rm E}[L_{n2a}]
&=&
-\frac{1}{2}n\kappa_n^2 \nu_n^2 \|\taub_n\|^2
\,
\bigg(
\varphi_f'(0)
c_{p}
\int_{-1}^1 
(1-s^2)^{(p-3)/2}
 s^2
\,ds
+o(1)
\bigg)
\nonumber
\\[1mm]
&=&
-\frac{1}{2p}n\kappa_n^2 \nu_n^2 \varphi_f'(0) \|\taub_n\|^2
+o(n\kappa_n^2 \nu_n^2)
=
-\frac{\xi^2}{2}\varphi_f'(0)\|\taub_n\|^2
+o(1)
.
\label{auboulot}
\end{eqnarray}
Now,
\begin{eqnarray*}
{\rm Var}[L_{n2a}]
&=&
n\kappa_n^2 \nu_n^2 
\,
{\rm Var}[  (\varphi_f(\kappa_n u_{n1})-1)
u_{n1} ] (\taub_n\pr\thetab)^2
\\[1mm]
&\leq &
n\kappa_n^2 \nu_n^2 
\,
{\rm E}[  (\varphi_f(\kappa_n u_{n1})-1)^2
u_{n1}^2 ] 
(\taub_n\pr\thetab)^2
\\[1mm]
&\leq &
\frac{1}{4}n\kappa_n^2 \nu_n ^4
\|\taub_n\|^4
\,
{\rm E}[  (\varphi_f(\kappa_n u_{n1})-1)^2 ] 
,
\end{eqnarray*}
where 
\begin{eqnarray}
\lefteqn{
{\rm E}[  (\varphi_f(\kappa_n u_{n1})-1)^2] 
=
\kappa_n^2
c_{p,\kappa_n,f}
\int_{-1}^1 
(1-s^2)^{(p-3)/2}
\Big(\frac{\varphi_f(\kappa_n s) -1}{\kappa_n s}\Big)^2 
s^2
f(\kappa_n s) 
\,ds
}
\nonumber
\\[1mm]
& & 
\hspace{2mm}
=
\kappa_n^2
(\varphi_f'(0))^2  
c_{p}
\int_{-1}^1 
(1-s^2)^{(p-3)/2}
s^2 
\,ds
+o(\kappa_n^2)
=
\frac{\kappa_n^2}{p}
(\varphi_f'(0))^2  
+o(\kappa_n^2)
=
O(\kappa_n^2)
.
\label{auboulot2}
\end{eqnarray}
Thus, 
$
{\rm Var}[L_{n2a}]=O(n\kappa_n^4\nu_n^4)=O(n^{-1}),
$
which, jointly with~(\ref{auboulot}), implies that
\begin{equation}
\label{steste1}
L_{n2a}=-\frac{\xi^2}{2}\varphi_f'(0)\|\taub_n\|^2+o_{\rm P}(1).
\end{equation}

Now, using~(\ref{auboulot2}) again, we obtain
\begin{eqnarray*}
{\rm E}[\|{\bf L}_{n2b}\|^2]
&=&
n \kappa_n^2 \nu_n^2\,
{\rm E}[
(\varphi_f(\kappa_n u_{n1})-1)^2
v_{n1}^2
]
\\[1mm]
&\leq&
n \kappa_n^2 \nu_n^2\,
{\rm E}[
(\varphi_f(\kappa_n u_{n1})-1)^2
]
=
O(n \kappa_n^4 \nu_n^2)
=
o(1)
.
\end{eqnarray*}
Consequently, ${\bf L}_{n2b}=o_{\rm P}(1)$, which, jointly with~(\ref{steste1}), establishes Part~(ii) of the result. 
\vspace{3mm}

(iii) Decomposing~$L_{n3}$ into~$\sum_{i=1}^n T_{ni}$, write
\begin{eqnarray*}
{\rm E}[T_{n1}]
&=&
{\rm E}[
\log f(\kappa_n u_{n1}+\kappa_n \nu_n(\taub_n'\thetab u_{n1}+\taub_n'\Sb_{n1}v_{n1}))
-
\log f(\kappa_n u_{n1})
\\[0mm]
& &
\hspace{50mm}
-
\kappa_n\nu_n  \varphi_f(\kappa_n u_{ni}) (\taub_n'\thetab u_{n1}+\taub_n'\Sb_{n1}v_{n1})
]
\\[2mm]
&=&
\frac{c_{p,\kappa_n,f}}{\mu_p}
\int_{-1}^1
\int_{\mathcal{S}^\perp_{\thetab}}
(1-s^2)^{(p-3)/2}
f(\kappa_n s)
\Big\{
\log f(\kappa_n s+\kappa_n \nu_n(\taub_n'\thetab s+\taub_n' \ub \sqrt{1-s^2}))
\\[2mm]
& &
\hspace{30mm}
-
\log f(\kappa_n s)
-
\kappa_n\nu_n  \varphi_f(\kappa_n s) (\taub_n'\thetab s+\taub_n'\ub\sqrt{1-s^2})
\Big\}
\,d\sigma(\ub)\,ds
,
\end{eqnarray*}
where~$\sigma(\cdot)$ stands for the surface area measure on~$\mathcal{S}^\perp_{\thetab}:=\{\xb\in\mathcal{S}^{p-1}:\xb'\thetab=0\}$ and~$\mu_p:=\sigma(\mathcal{S}^\perp_{\thetab})$. Letting~$t=\kappa_n s$ then provides
\begin{eqnarray*}
\lefteqn{
{\rm E}[T_{n1}]
=
\frac{c_{p,\kappa_n,f}}{\mu_p\kappa_n}
\int_{-\kappa_n}^{\kappa_n}
\int_{\mathcal{S}^\perp_{\thetab}}
(1-({\textstyle\frac{t}{\kappa_n}})^2)^{(p-3)/2}
f(t)
\Big\{
\log f(t+\kappa_n \nu_n(\taub_n'\thetab ({\textstyle\frac{t}{\kappa_n}})+\taub_n' \ub \sqrt{1-({\textstyle\frac{t}{\kappa_n}})^2}))
}
\\[2mm]
& &
\hspace{30mm}
-
\log f(t)
-
\kappa_n\nu_n  \varphi_f(t) (\taub_n'\thetab ({\textstyle\frac{t}{\kappa_n}})+\taub_n'\ub\sqrt{1-({\textstyle\frac{t}{\kappa_n}})^2})
\Big\}
\,d\sigma(\ub)\,dt
\\[2mm]
& &
\hspace{0mm}
=
\frac{c_{p,\kappa_n,f} \kappa_n\nu_n^2}{2\mu_p}
\int_{-\kappa_n}^{\kappa_n}
\int_{\mathcal{S}^\perp_{\thetab}}
(1-({\textstyle\frac{t}{\kappa_n}})^2)^{(p-3)/2}
 (\taub_n'\thetab ({\textstyle\frac{t}{\kappa_n}})+\taub_n'\ub\sqrt{1-({\textstyle\frac{t}{\kappa_n}})^2})^2
g_n(t,\ub)
f(t)
\,d\sigma(\ub)\,dt
,
\end{eqnarray*}
where 
\begin{eqnarray*}
\lefteqn{
g_n(t,\ub)
:=
}
\\[2mm]
& & 
\hspace{-6mm}
\frac{
\log f(t+\kappa_n \nu_n(\taub_n'\thetab ({\textstyle\frac{t}{\kappa_n}})+\taub_n' \ub \sqrt{1-({\textstyle\frac{t}{\kappa_n}})^2}))
-
\log f(t)
-
\kappa_n\nu_n  \varphi_f(t) (\taub_n'\thetab ({\textstyle\frac{t}{\kappa_n}})+\taub_n'\ub\sqrt{1-({\textstyle\frac{t}{\kappa_n}})^2})
}
{
\frac{1}{2}
\kappa_n^2\nu_n^2  (\taub_n'\thetab ({\textstyle\frac{t}{\kappa_n}})+\taub_n'\ub\sqrt{1-({\textstyle\frac{t}{\kappa_n}})^2})^2
}
\cdot
\end{eqnarray*}

Hence,
\begin{eqnarray*}
{\rm E}[T_{n1}]
&=&
\frac{c_{p,\kappa_n,f}  \kappa_n\nu_n^2}{2\mu_p}
 (\taub_n'\thetab)^2 
\int_{\mathcal{S}^\perp_{\thetab}}
\bigg(
\int_{-\kappa_n}^{\kappa_n}
({\textstyle\frac{t}{\kappa_n}})^2
(1-({\textstyle\frac{t}{\kappa_n}})^2)^{(p-3)/2}
g_n(t,\ub)
f(t)
\,dt
\bigg)
\,d\sigma(\ub)
\\[2mm]
& &
+
\frac{c_{p,\kappa_n,f}  \kappa_n\nu_n^2}{2\mu_p}
\,
\taub_n\pr
\bigg(
\int_{\mathcal{S}^\perp_{\thetab}}
\bigg(
\int_{-\kappa_n}^{\kappa_n}
(1-({\textstyle\frac{t}{\kappa_n}})^2)^{(p-1)/2}
g_n(t,\ub)
f(t)
\,dt
\bigg)
\,
\ub\ub\pr
\,d\sigma(\ub)
\bigg)
\taub_n
\\[2mm]
& &
+
\frac{c_{p,\kappa_n,f}  \kappa_n\nu_n^2}{\mu_p}
\,
(\taub_n'\thetab)
\taub_n\pr
\int_{\mathcal{S}^\perp_{\thetab}}
\bigg(
\int_{-\kappa_n}^{\kappa_n}
 ({\textstyle\frac{t}{\kappa_n}})
(1-({\textstyle\frac{t}{\kappa_n}})^2)^{(p-2)/2}
g_n(t,\ub)
f(t)
\,dt
\bigg)
\ub
\,d\sigma(\ub)
.
\end{eqnarray*}
By using the identities
$$
\int_{-\kappa_n}^{\kappa_n}
({\textstyle\frac{t}{\kappa_n}})^2
(1-({\textstyle\frac{t}{\kappa_n}})^2)^{(p-3)/2}
\,dt
=
\frac{\kappa_n}{c_p}
c_p
\int_{-1}^{1}
s^2
(1-s^2)^{(p-3)/2}
\,ds
=
\frac{\kappa_n}{c_p p}
,
$$
$$
\int_{-\kappa_n}^{\kappa_n}
(1-({\textstyle\frac{t}{\kappa_n}})^2)^{(p-1)/2}
\,dt
=
\kappa_n
\int_{-1}^{1}
(1-s^2)^{(p+2-3)/2}
\,ds
=
\frac{\kappa_n}{c_{p+2}}
=
\frac{\kappa_n(p-1)}{c_p p}
,
$$
and
$$
\int_{-\kappa_n}^{\kappa_n}
|{\textstyle\frac{t}{\kappa_n}}|
(1-({\textstyle\frac{t}{\kappa_n}})^2)^{(p-2)/2}
\,dt
=
\kappa_n
\int_{-1}^{1}
|s|(1-s^2)^{(p-2)/2}
\,ds
=
\frac{2\kappa_n}{p}
,
$$
we obtain
\begin{eqnarray*}
{\rm E}[T_{n1}]
&=&
\frac{c_{p,\kappa_n,f}  \kappa_n^2\nu_n^2}{2c_p p \mu_p}
 (\taub_n'\thetab)^2 
\int_{\mathcal{S}^\perp_{\thetab}}
\bigg(
\int_{-\kappa_n}^{\kappa_n}
h_n(t)
g_n(t,\ub)
f(t)
\,dt
\bigg)
\,d\sigma(\ub)
\\[2mm]
& &
+
\frac{c_{p,\kappa_n,f} (p-1) \kappa_n^2\nu_n^2}{2c_pp\mu_p}
\,
\taub_n\pr
\bigg(
\int_{\mathcal{S}^\perp_{\thetab}}
\bigg(
\int_{-\kappa_n}^{\kappa_n}
k_n(t)
g_n(t,\ub)
f(t)
\,dt
\bigg)
\ub\ub\pr
\,d\sigma(\ub)
\bigg)
\taub_n
\\[2mm]
& &
+
\frac{2c_{p,\kappa_n,f}  \kappa_n^2\nu_n^2}{p\mu_p}
\,
(\taub_n'\thetab)
\taub_n\pr
\int_{\mathcal{S}^\perp_{\thetab}}
\bigg(
\int_{-\kappa_n}^{\kappa_n}
\ell_n(t)
g_n(t,\ub)
f(t)
\,dt
\bigg)
\ub
\,d\sigma(\ub)
,
\end{eqnarray*}
where we let
$$
h_n(t)
:=
\frac{
({\textstyle\frac{t}{\kappa_n}})^2
(1-({\textstyle\frac{t}{\kappa_n}})^2)^{(p-3)/2}
}
{
\int_{-\kappa_n}^{\kappa_n}
({\textstyle\frac{t}{\kappa_n}})^2
(1-({\textstyle\frac{t}{\kappa_n}})^2)^{(p-3)/2}
\,dt
}
,
$$
$$
k_n(t)
:=
\frac{
(1-({\textstyle\frac{t}{\kappa_n}})^2)^{(p-1)/2}
}
{
\int_{-\kappa_n}^{\kappa_n}
(1-({\textstyle\frac{t}{\kappa_n}})^2)^{(p-1)/2}
\,dt
}
$$
and
$$
\ell_n(t)
:=
\frac{
{\textstyle\frac{t}{\kappa_n}}
(1-({\textstyle\frac{t}{\kappa_n}})^2)^{(p-2)/2}
}
{
\int_{-\kappa_n}^{\kappa_n}
|{\textstyle\frac{t}{\kappa_n}}|
(1-({\textstyle\frac{t}{\kappa_n}})^2)^{(p-2)/2}
\,dt
}
\cdot
$$
Splitting the third term of~${\rm E}[T_{n1}]$ according to~$\int_{-\kappa_n}^0+\int_0^{\kappa_n}$, we then obtain
\begin{eqnarray*}
{\rm E}[T_{n1}]
&=&
\frac{c_{p,\kappa_n,f}\kappa_n^2\nu_n^2}{2c_p p}
(\taub_n'\thetab)^2 
\bigg(
\frac{1}{\mu_p}
\int_{\mathcal{S}^\perp_{\thetab}}
\Big( (\log f)''(0) +o(1)\Big)
\,d\sigma(\ub)
\bigg)
\\[2mm]
& &
+
\frac{c_{p,\kappa_n,f}(p-1) \kappa_n^2\nu_n^2}{2c_p p}
\,
\taub_n\pr
\bigg(
\frac{1}{\mu_p}
\int_{\mathcal{S}^\perp_{\thetab}}
\Big( (\log f)''(0) +o(1)\Big)
\ub\ub\pr
\,d\sigma(\ub)
\bigg)
\taub_n
\\[2mm]
& &
-
\frac{c_{p,\kappa_n,f}  \kappa_n^2\nu_n^2}{c_p p}
\,
(\taub_n'\thetab)
\taub_n\pr
\bigg(
\frac{1}{\mu_p}
\int_{\mathcal{S}^\perp_{\thetab}}
\Big( (\log f)''(0) +o(1)\Big)
\ub
\,d\sigma(\ub)
\bigg)
\\[2mm]
& &
+
\frac{c_{p,\kappa_n,f}  \kappa_n^2\nu_n^2}{c_p p}
\,
(\taub_n'\thetab)
\taub_n\pr
\bigg(
\frac{1}{\mu_p}
\int_{\mathcal{S}^\perp_{\thetab}}
\Big( (\log f)''(0) +o(1)\Big)
\ub
\,d\sigma(\ub)
\bigg)
.
\end{eqnarray*}
Since the four $o(1)$'s in this expression are uniform in~$\ub$ and~$\mathcal{S}^\perp_{\thetab}$ is compact, it follows (by using~(\ref{spherconstraint}) that  
\begin{eqnarray*}
{\rm E}[T_{n1}]
&\!\!\!=\!\!\!&
O( \kappa_n^2\nu_n^4)
+
\frac{
(1+o(1))(p-1) \kappa_n^2\nu_n^2}{2p}
\,
\varphi_f'(0)
\,\taub_n\pr
\Big[\frac{1}{p-1}(\mathbf{I}_p-\thetab\thetab')\Big]
\taub_n
+
o( \kappa_n^2\nu_n^3)
\\[2mm]
&\!\!\!=\!\!\!&
\frac{
\kappa_n^2\nu_n^2}{2p}
\,
\varphi_f'(0)
\|\taub_n\|^2
+
o( \kappa_n^2\nu_n^2)
=
\frac{\xi^2}{2n}
\varphi_f'(0)
\|\taub_n\|^2 
+
o(n^{-1})
.
\end{eqnarray*}
Therefore, 
$$
{\rm E}[L_{n3}]
=
\frac{\xi^2}{2}
\varphi_f'(0)
\|\taub_n\|^2 
+
o(1)
.
$$
Thus it only remains to show that~${\rm Var}[L_{n3}]=o(1)$.

To do so, write
\begin{eqnarray*}
{\rm Var}[T_{n1}]
&\leq &
{\rm E}[
(\log f(\kappa_n u_{n1}+\kappa_n \nu_n(\taub_n'\thetab u_{n1}+\taub_n'\Sb_{n1}v_{n1}))
-
\log f(\kappa_n u_{n1})
\\[2mm]
& & 
\hspace{40mm}
-
\kappa_n\nu_n  \varphi_f(\kappa_n u_{ni}) (\taub_n'\thetab u_{n1}+\taub_n'\Sb_{n1}v_{n1})
)^2
]
\\[2mm]
&=&
\frac{c_{p,\kappa_n,f}}{\mu_p}
\int_{-1}^1
\int_{\mathcal{S}^\perp_{\thetab}}
(1-s^2)^{(p-3)/2}
f(\kappa_n s)
\Big\{
\log f(\kappa_n s+\kappa_n \nu_n(\taub_n'\thetab s+\taub_n' \ub \sqrt{1-s^2}))
\\[2mm]
& &
\hspace{20mm}
-
\log f(\kappa_n s)
-
\kappa_n\nu_n  \varphi_f(\kappa_n s) (\taub_n'\thetab s+\taub_n'\ub\sqrt{1-s^2})
\Big\}^2
\,d\sigma(\ub)\,ds
.
\end{eqnarray*}
Letting again~$t=\kappa_n s$ yields
\begin{eqnarray*}
{\rm Var}[T_{n1}]
&\leq &
\frac{c_{p,\kappa_n,f}}{\mu_p \kappa_n}
\int_{-\kappa_n}^{\kappa_n}
\int_{\mathcal{S}^\perp_{\thetab}}
(1-({\textstyle \frac{t}{\kappa_n}})^2)^{(p-3)/2}
f(t)
\Big\{
\log f(t+\kappa_n \nu_n(\taub_n'\thetab ({\textstyle \frac{t}{\kappa_n}})+\taub_n' \ub \sqrt{1-({\textstyle \frac{t}{\kappa_n}})^2}))
\\[2mm]
& &
\hspace{30mm}
-
\log f(t)
-
\kappa_n\nu_n  \varphi_f(t) (\taub_n'\thetab ({\textstyle \frac{t}{\kappa_n}})+\taub_n'\ub\sqrt{1-({\textstyle \frac{t}{\kappa_n}})^2})
\Big\}^2
\,d\sigma(\ub)\,dt
.
\end{eqnarray*}
Proceeding as for the expectation, we may then write
\begin{eqnarray*}
\hspace{-6mm}
{\rm Var}[T_{n1}]
&\!\!\!\leq \!\!\!&
\frac{c_{p,\kappa_n,f}\kappa_n^3\nu_n^4}{4\mu_p}
\int_{-\kappa_n}^{\kappa_n}
\int_{\mathcal{S}^\perp_{\thetab}}
(1-({\textstyle \frac{t}{\kappa_n}})^3)^{(p-3)/2}
 (\taub_n'\thetab ({\textstyle \frac{t}{\kappa_n}})+\taub_n'\ub\sqrt{1-({\textstyle \frac{t}{\kappa_n}})^2})^4
(g_n(t,\ub))^2
f(t)
\,d\sigma(\ub)\,dt
\\[2mm]
&\!\!\!\leq \!\!\!&
C\frac{c_{p,\kappa_n,f}\kappa_n^3\nu_n^4}{\mu_p}
\int_{-\kappa_n}^{\kappa_n}
\int_{\mathcal{S}^\perp_{\thetab}}
(1-({\textstyle \frac{t}{\kappa_n}})^2)^{(p-3)/2}
(g_n(t,\ub))^2
f(t)
\,d\sigma(\ub)\,dt
\\[2mm]
&\!\!\!\leq \!\!\!&
C\frac{c_{p,\kappa_n,f}\kappa_n^4\nu_n^4}{c_p \mu_p}
\int_{-\kappa_n}^{\kappa_n}
\int_{\mathcal{S}^\perp_{\thetab}}
m_n(t)
(g_n(t,\ub))^2
f(t)
\,d\sigma(\ub)\,dt
,
\end{eqnarray*}
where~$C$ is some positive constant and
$$
m_n(t)
:=
\frac{
(1-({\textstyle\frac{t}{\kappa_n}})^2)^{(p-3)/2}
}
{
\int_{-\kappa_n}^{\kappa_n}
(1-({\textstyle\frac{t}{\kappa_n}})^2)^{(p-3)/2}
\,dt
}
=
\frac{c_p}{\kappa_n}
(1-({\textstyle\frac{t}{\kappa_n}})^2)^{(p-3)/2}
.
$$
Hence, along the same lines as above, we obtain
\begin{eqnarray*}
\hspace{-6mm}
{\rm Var}[T_{n1}]
&\leq &
C\frac{c_{p,\kappa_n,f}\kappa_n^4\nu_n^4}{c_p \mu_p}
\int_{\mathcal{S}^\perp_{\thetab}}
\int_{-\kappa_n}^{\kappa_n}
m_n(t)
(g_n(t,\ub))^2
f(t)
\,dt
\,d\sigma(\ub)
\\[2mm]
&= &
C\frac{c_{p,\kappa_n,f}\kappa_n^4\nu_n^4}{c_p}
\bigg(
\frac{1}{\mu_p}
\int_{\mathcal{S}^\perp_{\thetab}}
\Big(\Big( (\log f)''(0) \Big)^2+o(1)\Big)
\,d\sigma(\ub)
\bigg)
\\[2mm]
&= &
C\kappa_n^4\nu_n^4(\varphi_f'(0))^2
+
o(\kappa_n^4\nu_n^4)
=
O(\kappa_n^4\nu_n^4)
=
O(n^{-2})
.
\end{eqnarray*}
Therefore, ${\rm Var}[L_{n3}]=n {\rm Var}[T_{n1}]=o(1)$, as was to be shown. 
\cqfd
\vspace{2mm}

\bibliographystyle{imsart-nameyear.bst} 
\bibliography{Vicinity.bib}           
\vspace{3mm} 


\end{document}